\documentclass[final]{siamltex}
\usepackage{xcolor}
\usepackage{amsmath}
\usepackage{amssymb}
\usepackage{mathrsfs}
\usepackage{graphicx}
\usepackage{bm}
\usepackage{ucs}
\usepackage[utf8x]{inputenc}
\usepackage{wrapfig}
\usepackage{color}
\usepackage{float}
\usepackage{epstopdf}
\usepackage{algorithm,algorithmic}
\usepackage[breaklinks,colorlinks=true,linkcolor=blue,citecolor=red,
  backref=page]{hyperref}

\def\btheta{{\boldsymbol{\theta}}}
\def\bbeta{{\boldsymbol{\beta}}}
\def\calP{{\bf P}}

\def\debproof{ {\bf Proof.} }
\def\finproof{\hfill {\small $\Box$}}

\DeclareMathAlphabet{\itbf}{OML}{cmm}{b}{it}
 \DeclareMathAlphabet\mathbfcal{OMS}{cmsy}{b}{n}

\renewcommand{\hat}{\widehat}
\def\EE{\mathbb{E}}
\def\RR{\mathbb{R}}

\def\bx{{{\itbf x}}}

\newcommand{\la}{\lambda}

\newcommand{\ep}{\varepsilon}

\newcommand{\om}{\omega}
\newcommand{\btau}{{\bm{\tau}}}
\newcommand{\bn}{{{\itbf n}}}
\newcommand{\xpar}{\bx_{\parallel}}
\newcommand{\cR}{\mathcal{R}}
\newcommand{\lin}{\left<}
\newcommand{\rin}{\right>}

\begin{document}
%%%%%%%%%%%%%%%%%%%%%%%%%%%%%
\title{Transport of power in random waveguides with turning points} 
\author{Liliana Borcea\footnotemark[1] \and Josselin
  Garnier\footnotemark[2] \and Derek Wood\footnotemark[1]} \maketitle
  
 %%%%%%%%%%%%%%%%%%%%%%%%%%%%%
 % footnotes for the addresses

\renewcommand{\thefootnote}{\fnsymbol{footnote}}
\footnotetext[1]{Department of Mathematics, University of Michigan,
  Ann Arbor, MI 48109. {{\tt borcea@umich.edu} and {\tt
      derekmw@umich.edu}}} \footnotetext[2]{Centre de
  Math\'{e}matiques Appliqu\'{e}es, Ecole Polytechnique, 91128
  Palaiseau Cedex, France. {\tt josselin.garnier@polytechnique.edu}}
 % headings \pagestyle{myheadings} \thispagestyle{plain}
\markboth{L. BORCEA, J. GARNIER, D. WOOD}{Transport in random waveguide}

%%%%%%%%%%%%%%%%%%%%%%%%%%%%%
\begin{abstract}
We present a mathematical theory of time-harmonic  wave propagation and reflection in a two-dimensional
random acoustic waveguide with sound soft boundary and turning points. The  boundary has small fluctuations on the scale of the wavelength, modeled as random. The waveguide supports multiple 
propagating modes. The number of these modes changes due to slow variations of  the waveguide cross-section. The changes 
occur at turning points, where waves transition from propagating to evanescent or the other way around.  We consider
a regime where scattering at the random boundary has
significant effect on the wave traveling from one turning point to another.
This effect is described by the coupling of its components, the modes. We derive the mode coupling theory from first principles,  and quantify the  randomization of the wave and the transport and reflection of power in the waveguide. 
We show in particular that scattering at the random boundary may increase or decrease the net power transmitted
through the waveguide depending on the source. 
\end{abstract}

%%%%%%%%%%%%%%%%%%%%%%%%%%%%%%%%%%%%%%%%%%%%%%%%%%%
\section{Introduction}
%%%%%%%%%%%%%%%%%%%%%%%%%%%%%%%%%%%%%%%%%%%%%%%%%%%
Guided waves have been studied extensively due to their numerous applications in electromagnetics 
\cite{collin1960field}, optics and communications \cite{marcuse2013theory,snyder1983optical}, quantum waveguides
\cite{exner2015quantum}, ocean acoustics \cite{tolstoy1966ocean}, etc.  
The classic waveguides, with straight boundaries and filled with media that are either homogeneous or do not vary along 
the direction of propagation of the  waves, are  well understood.  The wave equation in such waveguides is separable
and the solution is a superposition of independent modes, which are  propagating and evanescent modes and, in the case 
of penetrable boundary, radiation modes. 
When the geometry of the waveguide varies,  or the medium filling it is heterogeneous, the modes are coupled. Examples of mathematical 
studies that account for mode coupling and lead to numerical methods that model  such waveguides can be found in \cite{bruno1993numerical,ciraolo2008method,dinesen2001fast,dhia2000generalized,hazard2008improved}. We are interested 
in mode coupling due to small random perturbations of the waveguide, which can be quantified more explicitly using asymptotic analysis. 

Random models are useful for waveguides with rough boundaries and filled with composite media that vary at small scale, 
comparable  to the wavelength. Such variations are typically small and  unknown, so they introduce uncertainty in the model of wave propagation.  The random models of the  boundary and the wave speed are used to quantify how this uncertainty propagates 
to the uncertainty of the solution of the wave equation.   This is useful information in applications like imaging \cite{BKLU,borcea2014paraxial,borcea2013quantitative,acosta2015source}. The mode coupling theory in waveguides  filled 
with random media has been developed in \cite{kohler77,dozier1978statistics,garnier_papa,gomez,garnier2008effective} 
for sound waves and in  \cite{marcuse2013theory,alonso2015electromagnetic} for electromagnetic waves. The theory has also 
been extended to waveguides with random perturbations of straight  boundaries \cite{alonso2011wave,borcea2014paraxial,gomez2011wave}.

In this paper we develop  a mode coupling theory in random waveguides with turning points. 
We consider for simplicity a two-dimensional acoustic waveguide with sound soft boundary, but the theory can 
be extended to three dimensions and to electromagnetics.  The waveguide has a slowly bending 
axis, a slowly changing opening and a randomly perturbed boundary. The slow variations occur on 
a long scale with respect to the wavelength,  whereas the random perturbations are at a scale similar to the wavelength.

In the absence of the random perturbations,  waves in slowly changing waveguides can be analyzed with a  
local decomposition in modes that are approximately independent  \cite{ahluwalia1974asymptotic,ting1983wave}. This is known as the adiabatic approximation \cite[Section 19-2]{snyder1983optical}, and the result differs from that in waveguides 
with straight boundaries in one important aspect: The change in the opening of the waveguide causes the 
number of propagating modes to increase or decrease by $1$ at  locations called turning points. Modes transition 
there from propagating to evanescent or the other way around, and due to energy conservation, the impinging propagating 
mode is turned back i.e., is reflected. Turning waves in slowly changing waveguides are studied mathematically in \cite{anyanwu1978asymptotic}. A recent study of their interaction with a random boundary is given in \cite{borcea2016turning}, 
for the case of weak random fluctuations that affect only the turning modes.
Here we extend the results in \cite{borcea2016turning} to stronger  fluctuations, that couple all the waveguide modes. 

Starting from the wave equation, and using the separation of scales of variation of the waveguide, we derive an asymptotic 
model for the wave field that accounts for coupling of all the modes, propagating and evanescent. This coupling is
described by a stochastic system of differential equations for the random mode amplitudes, endowed with initial conditions that 
model the source excitation and radiation conditions.  The excitation is due to a point source, but due to the linearity of the 
equation, other  source excitations  can be handled by superposition. We obtain an extension of the  diffusion 
 approximation theorem proved in \cite{papanicolaou1974asymptotic} to carry out the asymptotic analysis of the mode 
 amplitudes. The result simplifies when the random fluctuations are smooth, because the forward and backward 
 going components of the propagating modes become independent. This is known as the forward scattering approximation,  and 
 applies to propagation between the turning points. At the turning points there is strong coupling of the components of the turning 
 waves, described by random reflection coefficients.  With the diffusion approximation theorem, and in the forward scattering 
 approximation regime, we quantify the net effect  of  the random boundary on the transmitted 
 and reflected power in the waveguide. This is the main result of the paper, and shows that  the random boundary can be useful for increasing the transmitted power.
 
The paper is organized as follows: We begin in section \ref{sect:formulation} with the formulation of the problem, and the derivation 
of the asymptotic model for the wave field. Then we give in section \ref{sect:wavedec} the mode decomposition and derive a 
closed system of stochastic differential equations for the random amplitudes of the propagating modes between turning points. These 
equations are complemented by source excitation conditions, radiation conditions, as well as continuity and reflection conditions 
at the turning points. The asymptotic limit of the solution of these equations and the forward scattering approximation are in section \ref{sect:results}. We use these results to quantify the transmitted and reflected power in the waveguide in section \ref{sect:netscat}. The diffusion approximation theorem used to carry the asymptotic limit is stated and proved in section 
\ref{app:adif}. We end with a summary in section \ref{sect:sum}.
% ------------------------
\section{Formulation of the problem}
\label{sect:formulation}
In this section we give the mathematical model for time-harmonic waves
in a random waveguide with variable cross-section and bending axis. We
begin in section \ref{sect:form1} with the setup, and describe the
scaling in section \ref{sect:form2} in terms of a small, dimensionless
parameter $\ep$. We use it in section \ref{sect:form3} to write the
wave problem in a form that can be analyzed in the asymptotic limit
$\ep \to 0$.
\subsection{Setup}
\label{sect:form1}

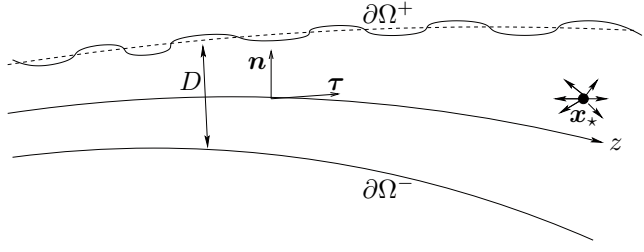
\begin{figure}[t]
\begin{picture}(0,0)%
\hspace{1in}\includegraphics[width = 0.65\textwidth]{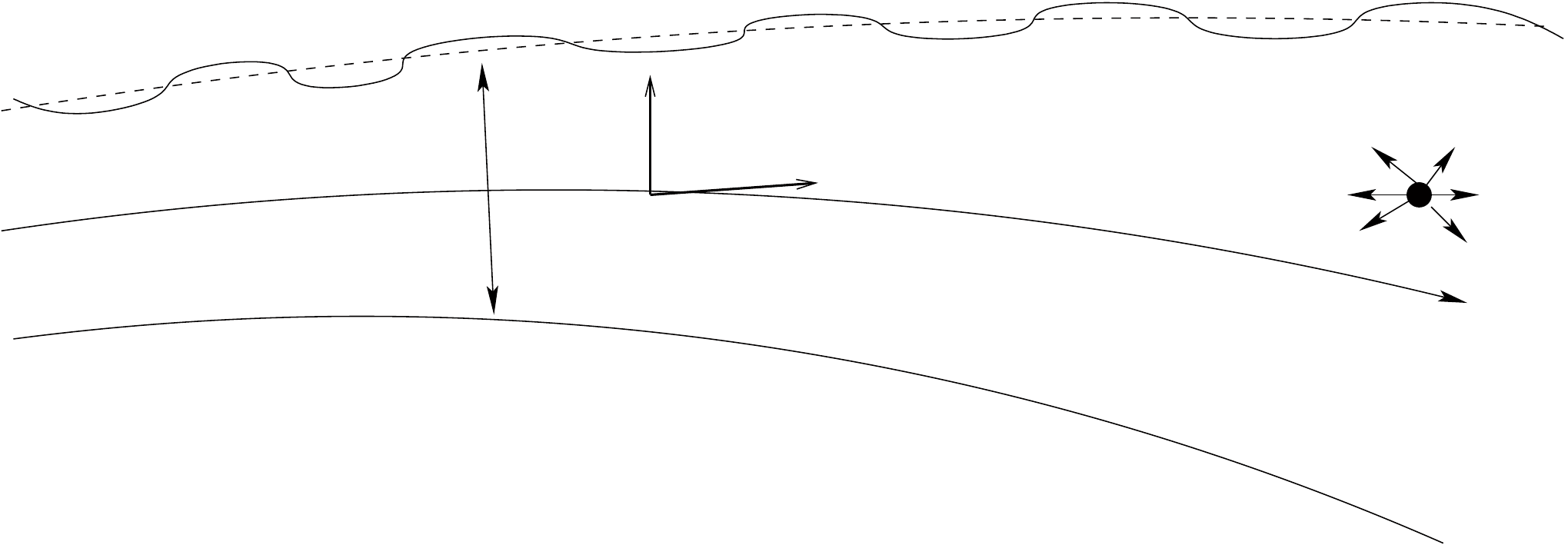}%
\end{picture}%
\setlength{\unitlength}{3947sp}%
\begingroup\makeatletter\ifx\SetFigFont\undefined%
\gdef\SetFigFont#1#2#3#4#5{%
  \reset@font\fontsize{#1}{#2pt}%
  \fontfamily{#3}\fontseries{#4}\fontshape{#5}%
  \selectfont}%
\fi\endgroup%
\begin{picture}(9770,1095)(1568,-1569)
\put(6550,-1000){\makebox(0,0)[lb]{\smash{{\SetFigFont{7}{8.4}{\familydefault}{\mddefault}{\updefault}{\color[rgb]{0,0,0}{\normalsize $z$}}%
}}}}
\put(6300,-820){\makebox(0,0)[lb]{\smash{{\SetFigFont{7}{8.4}{\familydefault}{\mddefault}{\updefault}{\color[rgb]{0,0,0}{\normalsize $\bx_\star$}}%
}}}}
\put(3850,-620){\makebox(0,0)[lb]{\smash{{\SetFigFont{7}{8.4}{\familydefault}{\mddefault}{\updefault}{\color[rgb]{0,0,0}{\normalsize $D$}}%
}}}}
\put(4780,-620){\makebox(0,0)[lb]{\smash{{\SetFigFont{7}{8.4}{\familydefault}{\mddefault}{\updefault}{\color[rgb]{0,0,0}{\normalsize $\btau$}}%
}}}}
\put(4280,-490){\makebox(0,0)[lb]{\smash{{\SetFigFont{7}{8.4}{\familydefault}{\mddefault}{\updefault}{\color[rgb]{0,0,0}{\normalsize ${\itbf n}$}}%
}}}}
\put(5000,-1300){\makebox(0,0)[lb]{\smash{{\SetFigFont{7}{8.4}{\familydefault}{\mddefault}{\updefault}{\color[rgb]{0,0,0}{\normalsize $\partial \Omega^-$}}%
}}}}
\put(5000,-200){\makebox(0,0)[lb]{\smash{{\SetFigFont{7}{8.4}{\familydefault}{\mddefault}{\updefault}{\color[rgb]{0,0,0}{\normalsize $\partial \Omega^+$}}%
}}}}
\end{picture}%
\caption{Illustration of a waveguide with slowly varying width $D$ and
  bending axis parametrized by the arc length $z$.  The boundary
  $\partial \Omega$ is the union of the curves $\partial \Omega^-$
  (the bottom boundary) and $\partial \Omega^+$ (the top boundary).
  The top boundary is perturbed by small random fluctuations. The unit
  tangent to the axis of the waveguide is denoted by $\btau$ and the
  unit normal $\bn$ points toward the upper boundary. The source of
  waves is at $\bx_\star$.  }
\label{fig:setup}
\end{figure}
We consider a two-dimensional acoustic waveguide with sound-soft
boundary.  The waveguide occupies the semi-infinite domain $\Omega$,
bounded above and below by two curves $\partial \Omega^+$ and
$\partial \Omega^-$, as shown in Figure \ref{fig:setup}. The top
boundary $\partial \Omega^+$ is perturbed by small random fluctuations
about the curve $\partial \Omega^+_o$ shown in the figure with the
dotted line. The axis of the waveguide is at half the distance $D$
between $\partial \Omega_0^+$ and $\partial \Omega^-$. It is a smooth
curve parametrized by the arc length $z \in \mathbb{R}$, that bends
slowly, meaning that its tangent $\btau(z/L)$ and curvature
$\kappa(z/L)$ vary on a scale $L$ which is large with respect to the
waveguide width $D(z/L)$. The width function $D$  has bounded first two derivatives, 
and  to avoid complications in the analysis of scattering of the waves at the random boundary, we
also assume that it is monotonically increasing.

Because of the changing geometry, it is convenient to use orthogonal
curvilinear coordinates with axes along $\btau(z/L)$ and $\bn(z/L)$,
where $\bn$ is the unit vector orthogonal to $\btau$, pointing toward
the upper boundary. For any $\bx \in \Omega$, written henceforth as
$\bx = (r,z)$, we have
\begin{equation}
\bx = \xpar(z) + r \bn \Big(\frac{z}{L}\Big),
\label{eq:coord1}
\end{equation}
where $\xpar(z)$ is along the waveguide axis at arc length $z$, satisfying 
\begin{equation}
\partial_z \xpar(z) = \btau\Big(\frac{z}{L}\Big),
\end{equation}
and $r$ is the coordinate in the normal direction. The domain $\Omega$
is the set
\begin{equation}
\Omega = \{(r,z): ~  z \in \mathbb{R}, ~ r \in (r^-(z),r^+(z))\},
\label{eq:coord1.p}
\end{equation}
where
\begin{equation}
r^-(z) = -\frac{D(z/L)}{2},
\label{eq:coord2}
\end{equation}
is at the bottom boundary $\partial \Omega^-$ and 
\begin{equation}
r^+(z) = \frac{D(z/L)}{2} \left[ 1 + 1_{(-Z_M,Z_M)}(z) \, \sigma \nu
  \Big(\frac{z}{\ell}\Big) \right],
\label{eq:coord3}
\end{equation}
is at the randomly perturbed top boundary $\partial \Omega^+$. The
perturbation is modeled by the random process $\nu$ and it extends
over the interval $(-Z_M,Z_M)$, the support of the indicator function
$1_{(-Z_M,Z_M)}(z)$, where $Z_M > L$ is a long scale needed to
impose outgoing boundary conditions on the waves.\footnote{In practice
  $Z_M$ may be chosen based on the duration of the observation time of
  the wave, using the hyperbolicity of the wave equation in the time
  domain.  The single frequency wave analyzed in this paper is the
  Fourier transform of the time dependent wave field.}  We let the
boundaries of the waveguide be straight and parallel for $|z| > Z_M$.

The random process $\nu$ is stationary with zero mean
\begin{equation}
\EE \big[\nu(\zeta)\big] = 0,
\label{eq:rp1}
\end{equation}
and auto-correlation function 
\begin{equation}
\cR(\zeta) = \EE \big[ \nu(0) \nu(\zeta)\big].
\label{eq:rp2}
\end{equation}
We assume that $\nu$  is mixing, with rapidly decaying mixing rate, as
defined for example in \cite[section 2]{papanicolaou1974asymptotic},
and it is bounded, with bounded first two derivatives, almost
surely. This implies in particular that $\cR$ is integrable and has at
least four bounded derivatives. We normalize $\nu$ by 
\begin{equation}
\cR(0) = 1, \qquad \int_{-\infty}^\infty d \zeta \, \cR(\zeta) =  O(1),
\label{eq:rp3}
\end{equation}
so that $\sigma$ in \eqref{eq:coord3} is the standard deviation of the
fluctuations of $\partial \Omega^+$, and $\ell$ quantifies their
correlation length.

The waves are generated by a point source at $\bx_\star = (r_\star,z_\star = 0)
\in \Omega$, which emits a complex signal $f(\om)$ at frequency
$\om$. We take the origin of $z$ at the source, so that $z_\star = 0$.
The waveguide is filled with a homogeneous medium with wave speed $c$,
and the wave field $p(\om,\bx)$ satisfies the Helmholtz equation
\begin{equation}
\Delta p(\om,\bx) + k^2 p(\om,\bx) = f(\om) \delta(\bx-\bx_\star),
\quad \bx = (r,z) \in \Omega,
\label{eq:we1}
\end{equation}
where $k = \om/c$ is the wavenumber. In curvilinear coordinates this
takes the form
\begin{align}
\left[\partial_r^2 -
  \frac{\frac{1}{L}\kappa\big(\frac{z}{L}\big)\partial_r}{1 -
    \frac{r}{L} \kappa\big(\frac{z}{L}\big)} +
  \frac{\partial_z^2}{\Big[ 1 - \frac{r}{L}
      \kappa\big(\frac{z}{L}\big)\Big]^2} + \frac{\frac{r}{L^2}
    \kappa'\big(\frac{z}{L}\big)\partial_z}{\Big[ 1 - \frac{r}{L}
      \kappa\big(\frac{z}{L}\big)\Big]^3} +k^2 \right] p(\om,r,z)
\nonumber \\= \left|1 - \frac{r_\star}{L} \kappa(0)\right|^{-1}
f(\om)\delta(z) \delta(r-r_\star),\label{eq:we3}
\end{align}
as shown in appendix \ref{ap:CurvCoord}, where $\kappa'$ is the
derivative of the curvature $\kappa$.  The sound soft boundary 
$\partial \Omega^+ \cup \partial \Omega^-$  is modeled by the homogeneous Dirichlet boundary conditions
\begin{equation}
p(\om,r^+(z),z) = p(\om,r^-(z),z) = 0,
\label{eq:we2}
\end{equation}
and at points $\bx = (r,z)$ with $|z| > Z_M$ we have
radiation conditions that state that $p(\om,r,z)$ is outgoing and
bounded.

\subsection{Scaling}
\label{sect:form2}
There are four length scales in the problem: The wavelength
$\la = 2 \pi/k$, the width of the waveguide $D$, the scale $L$ of the
slow variations of the waveguide, and the correlation length $\ell$ of
the random fluctuations of the boundary $\partial \Omega^+$. They
satisfy
\begin{equation} 
L \gg D \sim \la \sim \ell, 
\label{eq:sc1}
\end{equation}
where $\sim$ denotes ``of the same order as'', and we model the
separation of scales using the dimensionless parameter
\begin{equation}
\ep = \frac{\ell}{L}, \quad 0 < \ep \ll 1.
\label{eq:sc2}
\end{equation}
Our analysis of the wave field $p(\om,r,z)$ is in the asymptotic limit $\ep \to 0$.

As shown in section \ref{sect:wavedec}, the ratio of $D$ and $\la/2$
defines $ N(z) = \lfloor {2 D(z/L)}/{\la} \rfloor $, the number of
propagating components of the wave, called modes, where $\lfloor ~
\rfloor$ denotes the integer part. The assumption $D \sim \la$ in
\eqref{eq:sc1} means that
\begin{equation}
\label{eq:ASN}
N_{\rm min} \le N(z) \le N_{\rm max},
\end{equation}
for all $z$, 
where $N_{\rm min}$ and $N_{\rm max}$ are natural numbers, independent of $\ep$.

The scales $\la$ and $\ell$ are of the same order in \eqref{eq:sc1} so
that the waves interact efficiently with the random fluctuations of
the boundary. This interaction, called cumulative scattering,
randomizes the wave field as it propagates in the waveguide. The
distance from the source at which the randomization occurs depends on
the standard deviation $\sigma$ of the fluctuations. We scale $\sigma$
as
\begin{equation}
\sigma = \sqrt{\ep} \tilde{\sigma}, \quad \tilde{\sigma} = O(1),
\label{eq:sc4}
\end{equation}
so that we observe the randomization at distances $z \sim L$.

The scaled variables are defined as follows: The arc length $z$ is
scaled by $L$,
\begin{equation}
\tilde{z} = \frac{z}{L},
\label{eq:sc5}
\end{equation}
and the similar lengths
$D$, $r$ and $\la$ are scaled by $\ell$, to obtain
\begin{equation}
\tilde{D}(\tilde z) = \frac{D(z/L)}{\ell},
~~ \tilde r = \frac{r}{\ell}, ~ ~ \tilde k = k \ell.
\label{eq:sc5p}
\end{equation}
We also introduce the scaled bound $\tilde{Z}_M= {Z_M}/{L}$ of the support of the random
fluctuations,  which is a large number, independent of $\ep$.
\subsection{Asymptotic model}
\label{sect:form3}
Let us multiply equation \eqref{eq:we3} by $L^2 [ 1- r\kappa/L]^2$ and
use the scaling relations \eqref{eq:sc4}-\eqref{eq:sc5p}. Dropping the 
tilde to simplify notation, because all variables are scaled henceforth,  we obtain
\begin{align}
&\hspace{-0.15in}\left[\partial_z^2 + \frac{(1-\ep r
      \kappa(z))^2}{\ep^2}(\partial_r^2+k^2) - \frac{\kappa(z)(1-\ep r
      \kappa(z))}{\ep} \partial_r+ \frac{\ep r \kappa'(z) }{(1-\ep r
      \kappa(z))} \partial_z\right] p(\om,r,z) \nonumber
  \\ &\hspace{2in}= \frac{f(\om)[1-\ep r_\star \kappa(0)]}{\ep}
  \delta(r-r_\star) \delta(z),
\label{eq:AM1}
\end{align}
with homogeneous Dirichlet boundary conditions \eqref{eq:we2} at 
\begin{equation}
r^-(z) = -\frac{D(z)}{2}, ~ ~ r^+(z) = \frac{D(z)}{2} \left[ 1+
  1_{(-Z_M,Z_M)}(z) \sqrt{\ep} \sigma \nu\Big(\frac{z}{\ep}\Big) \right],
\label{eq:AM2}
\end{equation}
and appropriate radiation conditions for $|z| > Z_M$. These
equations define the asymptotic model for the wave field, and we
wish to analyze it in the limit $\ep \to 0$.

The boundary has $\ep$ dependent fluctuations, so to ensure that the
boundary conditions are satisfied at all orders of $\ep$, we change
variables to
\begin{equation}
r = \rho + \frac{[2 \rho + D(z)]}{4} \sqrt{\ep} \sigma \nu
\Big(\frac{z}{\ep}\Big),
\label{eq:AM3}
\end{equation}
for $|z| < Z_M$, and denote the transformed wave field by
\begin{equation}
p^\ep(\om,\rho,z) = p\Big(\om,\rho+ \frac{(2 \rho + D(z))}{4}
\sqrt{\ep} \sigma \nu \Big(\frac{z}{\ep}\Big),z \Big).
\label{eq:AM4}
\end{equation}
At $|z| > Z_M$ there are no fluctuations so the transformation is the
identity $r = \rho$. We use the same notation $p^\ep$ for the wave
field at all $z \in \mathbb{R}$, and analyze it separately in the regions
with the random fluctuations and without. The results are connected by
continuity at $z = \pm Z_M$.

The change of variables \eqref{eq:AM3} makes the boundary conditions 
independent of $\ep$,
\begin{equation}
p^\ep\Big(\om,\pm\frac{D(z)}{2},z\Big) = 0,
\label{eq:AM6}
\end{equation}  
and maps the random fluctuations to the differential operator in the
wave equation. Explicitly, we show in appendix \ref{ap:DerAs} that the
wave equation becomes
\begin{align}
\sum_{j=0}^\infty \ep^{j/2 - 1} \mathcal{L}_j^\ep \, p^\ep(\om,\rho,z)
=\hat f(\om)\big[ 1 + O(\sqrt{\ep})\big] \delta(\rho-\rho_\star)
\delta(z), \label{eq:AM5}
\end{align}
for $|\rho| < D(z)/2$ and $|z| < Z_M$, with the leading term in the
expansion of the operator 
\begin{equation}
\mathcal{L}_0^\ep = \big(\ep \partial_z \big)^2 + \partial_\rho^2 + k^2.
\label{eq:AM7}
\end{equation}
This is the Helmholtz operator in a perfect waveguide, with straight
and parallel boundaries. The random fluctuations appear in the first
perturbation operator,
\begin{equation}
\mathcal{L}_1^\ep = -\sigma\Big\{ \nu\Big(\frac{z}{\ep}\Big)
\partial_\rho^2 + \frac{[2 \rho + D(z)]}{4} \Big[
  \nu''\Big(\frac{z}{\ep}\Big) \partial_\rho + 2
  \nu'\Big(\frac{z}{\ep}\Big) \ep \partial_{\rho
    z}^2\Big]\Big\}.\label{eq:AM8}
\end{equation}
The second perturbation operator has a deterministic part, due to the
curvature of the axis of the waveguide, and a random part, quadratic
in the random fluctuations,
\begin{align}
\mathcal{L}_2^\ep = &-\kappa(z)\big[2 \rho  \big(\partial_\rho^2 + k^2) +
\partial_\rho\big] + \frac{\sigma^2}{4}\Big\{3
\nu^2\Big(\frac{z}{\ep}\Big)+ \frac{\big[2 \rho + D(z)\big]^2}{4}
{\nu'}^{\,2}\Big(\frac{z}{\ep}\Big) \Big\}\partial_\rho^2
\nonumber \\ &+ \frac{[2 \rho + D(z)]\sigma^2}{4} \Big\{
\nu\Big(\frac{z}{\ep}\Big)\nu'\Big(\frac{z}{\ep}\Big) \ep
\partial_{\rho z}^2 + \Big[{\nu'}^{\,2}\Big(\frac{z}{\ep}\Big) +
\frac{1}{2}
\nu''\Big(\frac{z}{\ep}\Big)\nu\Big(\frac{z}{\ep}\Big)\Big]
\partial_\rho\Big\}.
\label{eq:AM9}
\end{align}
The remaining operators in the asymptotic series in \eqref{eq:AM5}
depend on higher powers of the fluctuations $\nu$, but play no role in
the limit $\ep \to 0$.

By assumption, there are no variations of the waveguide at $|z| >
Z_M$, so the operator in the left hand side of \eqref{eq:AM5} reduces
to $\mathcal{L}_0^\ep$ in this region.

\section{Mode decomposition and coupling}
\label{sect:wavedec}
To analyze the solution of the wave equation \eqref{eq:AM5} with
boundary conditions \eqref{eq:AM6} in the limit $\ep \to 0$, we begin
in section \ref{sect:dec1} with the mode decomposition of the wave
field. The modes are special solutions of the wave equation with 
operator \eqref{eq:AM7}. They represent propagating and
evanescent waves which are coupled by the perturbation operators
\eqref{eq:AM8}-\eqref{eq:AM9}, as explained in section
\ref{sect:dec2}. We are interested in the propagating modes, which are
left and right going waves with random amplitudes satisfying a
stochastic system of equations derived in section \ref{sect:dec3}.  It
is this system that we analyze in the asymptotic limit $\ep \to 0$ to
quantify the cumulative scattering effects in the waveguide.
\subsection{Mode decomposition}
\label{sect:dec1}
The second term in \eqref{eq:AM7} is the Sturm-Liouville operator
$\partial_\rho^2 +k^2$ acting on functions that vanish at $ \rho =
\pm D(z)/2$, for any given $z$. Its eigenvalues $\lambda_j$ are real
and distinct
\begin{equation}
\lambda_j(z) = k^2 - \mu_j^2(z), \qquad \mu_j(z) = \frac{\pi
  j}{D(z)}, ~~ j = 1, 2, \ldots
\label{eq:eigval}
\end{equation}
and the eigenfunctions
\begin{equation}
y_j(\rho,z) = \sqrt{\frac{2}{D(z)}} \sin \left[\frac{(2
    \rho + D(z))}{2} \mu_j(z) \right],
\label{eq:eigf}
\end{equation}
form an orthonormal $L^2$ basis in $[-D(z)/2,D(z)/2]$. We use this
basis to decompose the solution of \eqref{eq:AM7} in one dimensional
waves $p_j^\ep(\om,z)$ called modes, for each $z$,
\begin{equation}
p^\ep(\om,\rho,z) = \sum_{j=1}^\infty p_j^\ep(\om,z) y_j(\rho,z).
\label{eq:DC1}
\end{equation}
As shown in appendix \ref{ap:ModeCoup}, the modes can be written as
\begin{equation}
p_j^\ep(\om,z) = u_j^\ep(\om,z)[1+O(\sqrt{\ep})],
\label{eq:DC2}
\end{equation}
with $u_j^\ep(\om,z)$ satisfying a coupled system of one dimensional
wave equations that we now describe:

In the perturbed section $|z| < Z_M$ of the waveguide, $u_j^\ep$ satisfies 
\begin{align} 
\frac{1}{\ep} \Big[ (\ep \partial_z)^2 + k^2 -
  \mu_j^2(z)\Big]u_j^\ep(\om,z) + \frac{\sigma}{\sqrt{\ep}}
\mu_j^2(z) \nu\Big(\frac{z}{\ep}\Big)u_j^\ep(\om,z) + \sigma^2
g_j^\ep(\om,z) u_j^\ep(\om,z) \nonumber \\\approx \mathcal{C}_j^\ep(\om,z) + f(\om)
y_j(\rho_\star,0)\delta(z),
\label{eq:DC3}
\end{align}
where the approximation indicates that we dropped lower order terms that 
have no contribution in the limit $\ep \to 0$. The coefficient 
$g_j^\ep$ in the left hand side is 
\begin{equation}
g_j^\ep(\om,z) = -\frac{3}{4} \mu_j^2(z) \nu^2\Big(\frac{z}{\ep}\Big) -
\Big[ \frac{(\pi j)^2}{12} + \frac{1}{16}\Big]
    {\nu'}^{\,2}\Big(\frac{z}{\ep}\Big) ,
\label{eq:DC4}
\end{equation}
and 
\begin{align}
\mathcal{C}_{j}^\ep(\om,z) &= \hspace{-0.05in}\sum_{q=1, q\ne
  j}^\infty \Big[ \frac{\sigma\Gamma_{jq}}{\sqrt{\ep}}
  \nu''\Big(\frac{z}{\ep}\Big) + \sigma^2 \gamma_{jq}\Big(\frac{z}{\ep}\Big) +
  \gamma_{jq}^o(z)\Big]u_q^\ep(\om,z) \nonumber
\\ &+\hspace{-0.05in}\sum_{q=1, q\ne j}^\infty \Big[
  \frac{\sigma\Theta_{jq}}{\sqrt{\ep}} \nu'\Big(\frac{z}{\ep}\Big) +
  \sigma^2 \theta_{jq}\Big(\frac{z}{\ep}\Big) + \theta_{jq}^o(z)\Big] \ep \partial_z
u_q^\ep(\om,z),
\label{eq:DC5}
\end{align}
models the coupling between the modes. The leading coupling
coefficients $\Gamma_{jq}$ and $\Theta_{jq}$ are constants, given in
equation \eqref{eq:DC6} in appendix \ref{ap:ModeCoup}.  The second
order coefficients $\gamma_{jq}(z/\ep)$ and $\theta_{jq}(z/\ep)$ are
quadratic in $\nu(z/\ep)$, as described in equations
\eqref{eq:DC7}-\eqref{eq:DC8}, and the coefficients $\gamma_{jq}^o(z)$
and $\theta_{jq}^o(z)$, given in equations
\eqref{eq:DC9}-\eqref{eq:DC10}, are due to the slow changes in the
waveguide.

In the region $|z| > Z_M$, where the waveguide has straight and
parallel boundaries, the wave equation simplifies to
\begin{align} 
\hspace{-0.05in}\frac{1}{\ep} \Big[ (\ep \partial_z)^2 + k^2 -
  \mu_j^2(z)\Big]u_j^\ep(\om,z) = 0.
\label{eq:DC11}
\end{align}
Depending on the index $j$, its solution is either an outgoing
propagating wave or a decaying evanescent wave. This wave is connected
to the solution of \eqref{eq:DC3} by the continuity of $u_j^\ep$ and
$\partial_z u_j^\ep$ at $z = \pm Z_M$.
\subsection{Random mode amplitudes}
\label{sect:dec2}
Equations \eqref{eq:DC3} are perturbations of the wave equation with
operator $(\ep \partial_z)^2 + k^2 - \mu_j^2(z)$, where the
perturbation term models the coupling of the modes.  This coupling is
similar to that in waveguides with randomly perturbed parallel
boundaries, studied in \cite{alonso2011wave,BKLU}, but the slow
variation of the waveguide introduces two differences: The first is
the presence of the extra terms $\gamma_{jq}^o(z)$
and $\theta_{jq}^o(z)$ in \eqref{eq:DC5}, given by 
\eqref{eq:DC9}-\eqref{eq:DC10}, which turn out to play no role in the
limit $\ep \to 0$. The second difference is important, as it gives a
$z$ dependent number
\begin{equation}
N(z) = \left \lfloor \frac{k D(z)}{\pi} \right\rfloor
\label{eq:DC12}
\end{equation}
of mode indexes $j = 1, \ldots, N(z)$ for which $k^2 - \mu_j^2(z) >
0$. These modes are oscillatory functions in $z$, and represent left
and right going waves. For indexes $j > N(z)$ the modes are decaying
evanescent waves.

\subsubsection{Turning points}
\label{sect:TP}
The function \eqref{eq:DC12} that defines the number of propagating
modes is piecewise constant. Starting from the origin, where we denote
the number of propagating modes by $N^{(0)} = N(0)$, the function
\eqref{eq:DC12} increases by $1$ at arc lengths $z_+^{(t)} >0$, for $t
= 1, \ldots , t_M^+$, and decreases by $1$ at $z_{-}^{(t)} <0$, for
$t= 1, \ldots , {t}_{-M}$. The jump locations $z_{\pm}^{(t)}$, ordered
as
\[
-Z_M < \ldots < z_{-}^{(2)} < z_{-}^{(1)} < 0 < z_{+}^{(1)} <
z_{+}^{(2)} < \ldots < Z_M,
\] 
are the zeroes of the eigenfunctions \eqref{eq:eigval}, and are called
turning points \cite{lynn1970uniform,anyanwu1978asymptotic}. We
assume henceforth that the monotonically increasing $D(z)$ satisfies
\begin{equation}
D' \big(z_{\pm}^{(t)}\big) >0, \quad \forall \, t \ge 1,
\label{eq:DC13p}
\end{equation}
so that the turning points are simple and isolated.  Consistent with
our scaling, they are spaced at order one scaled distances.

Between any two consecutive turning points $z_\pm^{(t-1)}$ and
$z_\pm^{(t)}$, where we set by convention $ z_\pm^{(0)} = 0, $ the number
of propagating modes equals the constant
\begin{equation}
N_{\pm}^{(t-1)} = N^{(0)} \pm (t-1).
\label{eq:DC13pp}
\end{equation}
This number is bounded above and below as in \eqref{eq:ASN}, with
$N_{\rm min} = N(-Z_M)$ and $N_{\rm max} = N(Z_M),$ so the bounds $t_M^+$ and
$t_M^{-}$ on the indexes $t$ are
\begin{equation}
t_M^{-} = N^{(0)} - N_{\rm min}+1 ~ ~ \mbox{and} ~ ~  t_M^+ = N_{\rm max}-N^{(0)}+1.
\end{equation}

Beginning from the source location $z = 0$, which we assume   is
not a turning point,
$z_{-}^{(t)}$ is defined as the unique, negative arc-length satisfying
\begin{equation}
k = \frac{\pi N_{-}^{(t-1)}}{D\big(z_{-}^{(t)}\big)}, \quad t = 1,
\ldots , t_M^{-},
\label{eq:DC14}
\end{equation}
where the uniqueness is due to the monotonicity of $D(z)$. Similarly,
the jump location $z_{+}^{(t)}$ is defined as the unique, positive arc
length satisfying
\begin{equation}
k = \frac{\pi \big(N_{+}^{(t-1)}+1\big)}{D\big(z_{+}^{(t)}\big)},
\quad t = 1, \ldots , t_M^+.
\label{eq:DC15}
\end{equation}

The analysis of the modes is similar on the left and right of the source, 
so we focus attention in this section on a sector
$
z \in \big(z_{-}^{(t)},z_{-}^{(t-1)}\big)
$
of the waveguide, for some $1 \le t \le t_{M}^-$, and simplify the
notation for the number \eqref{eq:DC13pp} of propagating modes
\begin{equation}
\mathcal{N} = N_{-}^{(t-1)}.
\label{eq:FM1}
\end{equation}
These modes are a superposition of right and left going waves, with
random amplitudes that model cumulative
scattering in the waveguide, as we explain in the next section.

\subsubsection{The left and right going waves}
\label{sect:TP1}
We decompose the propagating modes in left and right going waves,
using a flow of smooth and invertible matrices ${\bf M}_j^\ep(\om,z)$,
\begin{equation}
\begin{pmatrix} a_j^\ep(\om,z) \\ b_j^\ep(\om,z) \end{pmatrix} =
{\bf M}_j^{\ep,-1}(\om,z) \begin{pmatrix} u_j^\ep(\om,z)\\
  v_j^\ep(\om,z)\end{pmatrix},
\label{eq:FM5}
\end{equation}
where ${\bf M}_j^{\ep,-1}$ denotes the inverse of ${\bf M}_j^\ep$ and
\begin{equation}
v_j^\ep(\om,z) = - i \ep \partial_zu_j^\ep(\om,z), \qquad j = 1,
\ldots, \mathcal{N}.
\label{eq:FM4}
\end{equation}
We obtain from \eqref{eq:DC3} that
\begin{align}
 \partial_z \begin{pmatrix} a_j^\ep(\om,z)
   \\ b_j^\ep(\om,z) \end{pmatrix} &= {\bf M}_j^{\ep,-1}(\om,z)
 \left\{ \frac{i}{\ep} \begin{pmatrix} 0 & 1 \\ k^2(\om) - \mu_j^2(z)
   & 0 \end{pmatrix}{\bf M}_j^\ep(\om,z)- \partial_z {\bf
   M}_j^\ep(\om,z) \right.  \nonumber \\ & \left. + \left[\frac{i
     \sigma}{\sqrt{\ep}} \mu_j^2(z) \nu \Big(\frac{z}{\ep}\Big) + i
   \sigma^2 g_j^\ep(\om , z) \right]
 \begin{pmatrix}0 & 0 \\ 1 &
  0 \end{pmatrix} {\bf M}_j^\ep(\om,z)\right\} 
\begin{pmatrix} a_j^\ep(\om,z)
   \\ b_j^\ep(\om,z) \end{pmatrix} \nonumber \\ & - i\,
\mathcal{C}_j^\ep(\om,z) {\bf M}_j^{\ep,-1}(\om,z) \begin{pmatrix} 0
  \\ 1
\end{pmatrix},
\label{eq:FM6}
\end{align}
and the decomposition is achieved by a flow ${\bf M}_j^\ep(\om,z)$
that removes to leading order the large deterministic term in
\eqref{eq:FM6}, the first line in the right hand side.

The matrix ${\bf M}_j^\ep(\om,z)$ must have the structure
\begin{equation}
{\bf M}_j^\ep(\om,z) = \begin{pmatrix} M_{j,11}^\ep(\om,z) & -
  \overline{M_{j,11}^\ep}(\om,z) \\ M_{j,21}^\ep(\om,z) &
  \overline{M_{j,21}^\ep}(\om,z) \end{pmatrix}, 
\label{eq:res6}
\end{equation}
where the bar denotes complex conjugate, so that the decomposition
\eqref{eq:FM5} conserves energy. The expression of the components in
\eqref{eq:res6} depends on the mode index, more precisely on the mode
wave number denoted by
\begin{equation}
\beta_j(\om,z) = \sqrt{k^2-\mu_j^2(z)}.
\label{eq:FM7}
\end{equation}
Note that $\beta_j$ is bounded away from zero for all $j = 1, \ldots,
\mathcal{N}-1$, and it approaches zero as $z \searrow z_{-}^{(t)}$,
for $j = \mathcal{N}$.  This last mode is a turning wave which
transitions from a propagating wave at $z \in
(z_{-}^{(t)},z_{-}^{(t-1)})$ to an evanescent wave at $z <
z_{-}^{(t)}$, as described in section \ref{sect:TP2}. Here we give the
decomposition of the modes indexed by $j \le \mathcal{N}-1$.

The entries of  \eqref{eq:res6} are defined by
\begin{align}
M_{j,11}^\ep(\om,z) &= \frac{1}{\sqrt{\beta_j(\om,z)}}
\exp\Big[\frac{i}{\ep} \int_{0}^z dz'
  \beta_j(\om,z')\Big], \nonumber \\ M_{j,21}^\ep(\om,z) &=
\beta_j(\om,z) M_{j,11}^\ep(\om,z), \label{eq:FM8}
\end{align}
for $j = 1, \ldots, \mathcal{N}-1$. This definition looks the same as
in perfect waveguides with straight and parallel boundary, except that
the mode wave number $\beta_j$ varies with $z$. We obtain from
\eqref{eq:res6}-\eqref{eq:FM8} that the determinant of ${\bf
  M}_{j}^\ep(\om,z)$ is constant
\begin{equation}
\mbox{det}\, {\bf M}_{j}^\ep(\om,z) = 2, \qquad  \forall \, z \in 
\big(z_{-}^{(t)},z_{-}^{(t-1)}\big),
\label{eq:FM9}
\end{equation}
so the matrix is invertible, and the decomposition \eqref{eq:FM5} 
can be rewritten as 
\begin{align}
u_j^\ep(\om,z) = \frac{1}{\sqrt{\beta_j(\om,z)}} \Big[a_j^\ep(\om,z)
e^{\frac{i}{\ep} \int_{0}^z dz' \beta_j(\om,z')}
- b_j^\ep(\om,z) e^{-\frac{i}{\ep} \int_{0}^z dz'
  \beta_j(\om,z')} \Big], \label{eq:FM10}
\end{align}
and  
\begin{align} \ep \partial_z u_j^\ep(\om,z)
= i\sqrt{\beta_j(\om,z)} \Big[a_j^\ep(\om,z)
e^{\frac{i}{\ep} \int_{0}^z dz' \beta_j(\om,z')}+ b_j^\ep(\om,z) e^{- \frac{i}{\ep} \int_{0}^z dz'
  \beta_j(\om,z')}\Big] . \label{eq:FM11}
\end{align}

Note that equations \eqref{eq:FM10}-\eqref{eq:FM11} are just the method of variation of
parameters for the perturbed wave equation satisfied by the $j$-th
mode. They decompose the mode in a right going wave with amplitude
$a_j^\ep$ and a left going wave with amplitude $b_j^\ep$. In perfect
waveguides these amplitudes would be constant, meaning physically that
the waves are independent.  In our case the amplitudes are random
fields, satisfying the system of stochastic differential equations
\begin{align} 
\partial_z \begin{pmatrix} a_j^\ep(\om,z)
  \\ b_j^\ep(\om,z) \end{pmatrix} = {\bf H}_j^\ep(\om,z) \begin{pmatrix}
  a_j^\ep(\om,z) \\ b_j^\ep(\om,z) \end{pmatrix} - \frac{i}{2} 
\mathcal{C}_j^\ep(\om,z) \begin{pmatrix} \overline{M_{j,11}^\ep}(\om,z) \\
M_{j,11}^\ep(\om,z) \end{pmatrix},
\label{eq:FM12}
\end{align}
obtained by substituting \eqref{eq:res6} and \eqref{eq:FM8} in
\eqref{eq:FM6}.  Here ${\bf H}_j^\ep(\om,z)$ is the matrix valued random
process
\begin{align}
{\bf H}_j^\ep(\om,z) = \begin{pmatrix} H_{j}^{\ep(aa)}(\om,z) &
  H_{j}^{\ep(ab)}(\om,z) \\ H_{j}^{\ep(ba)}(\om,z) &
  H_{j}^{\ep(bb)}(\om,z) \end{pmatrix},
\label{eq:res24}
\end{align}
with entries satisfying the relations
\begin{equation}
H_{j}^{\ep(ba)}(\om,z) = \overline{H_{j}^{\ep(ab)}}(\om,z), \quad 
H_{j}^{\ep(bb)}(\om,z) = \overline{H_{j}^{\ep(aa)}}(\om,z),
\label{eq:res24p}
\end{equation}
and taking the values 
\begin{align}
H_{j}^{\ep(aa)}(\om,z) \approx \frac{i}{2 \beta_j(\om,z)} \Big[
  \frac{\sigma}{\sqrt{\ep}} \mu_j^2(z) \nu \Big(\frac{z}{\ep} \Big) +
  \sigma^2 g_j^\ep(\om,z)\Big],
\label{eq:FM14}
\end{align}
and 
\begin{align}
H_{j}^{\ep(ab)}(\om,z) \approx \Big[ \overline{H_{j}^{\ep(aa)}}(\om,z) -
  \frac{\partial_z \beta_j(\om,z)}{2 \beta_j(\om,z)}\Big] \exp \Big[-
  \frac{2 i}{\ep} \int_{0}^z dz' \beta_j(\om,z') \Big].
\label{eq:FM14p}
\end{align}
As before, the approximation means up to negligible terms in the limit
$\ep \to 0$.

Equations \eqref{eq:FM12} show that the amplitudes of the $j$-th mode
are coupled to each other by the process ${\bf H}_j^\ep$, and to the other
modes by $\mathcal{C}_j^\ep$, defined by the series \eqref{eq:DC5}.
The first terms in this series involve the propagating waves
$u_{q}^\ep(\om,z)$, for $q \ne j$, decomposed as in
\eqref{eq:FM10}-\eqref{eq:FM10}. We describe in the next two sections
the turning and the evanescent waves.

\subsubsection{The turning waves}
\label{sect:TP2}
The mode indexed by $j = \mathcal{N}$ transitions at $z = z_{-}^{(t)}$
from propagating  to evanescent. This transition is
captured by the matrix ${\bf M}_{\mathcal{N}}^\ep(\om,z)$, which has
the same structure as in \eqref{eq:res6}, but its entries are
defined in terms of Airy functions \cite[chapter
  10]{abramowitz1972handbook}. This is because near the simple turning
point $z_{-}^{(t)}$, equation \eqref{eq:DC3} for $j = \mathcal{N}$ is a
perturbation of Airy's equation. 
We refer to \cite{anyanwu1978asymptotic, snyder1983optical} for
classic studies of turning waves in waveguides, and to
\cite{borcea2016turning} for an analysis of their interaction with the
random boundary. The setup in \cite{borcea2016turning} is the same as
here, with the exception that we consider a larger standard deviation
of the random fluctuations, to observe mode coupling in the waveguide.

We use the same ${\bf M}_{\mathcal{N}}^\ep(\om,z)$ as in
\cite{borcea2016turning}, with entries
\begin{align}
M_{\mathcal{N},11}^\ep(\om,z) = \ep^{-1/6}\sqrt{\pi}
Q_{\mathcal{N}}(\om,z) \exp \Big[ -i
  \frac{\phi_{\mathcal{N}}\big(\om,0\big)}{\ep} + \frac{i \pi}{4}
  \Big] \nonumber \\ \times\Big[ A_i\big(-\eta_{\mathcal{N}}^\ep(\om,z)\big) -i
  B_i\big(-\eta_{\mathcal{N}}^\ep(\om,z)\big)\Big],
\label{eq:TW1}
\end{align}
and 
\begin{align}
M_{\mathcal{N},21}^\ep(\om,z) = -i \ep \partial_z M_{\mathcal{N},11}^\ep(\om,z),
\label{eq:TW2}
\end{align}
for $z \in \big(z_{-}^{(t)}-\delta,z_{-}^{(t-1)}\big),$ where $\delta$
is a small, positive number, independent of $\ep$. We go slightly
beyond the turning point to capture the transition of the wave to an
evanescent one. The phase in definition \eqref{eq:TW1} is given by the
function
\begin{equation}
\phi_{\mathcal{N}}(\om,z) = \int_{z_{-}^{(t)}}^z dz' \, \sqrt{|k^2 - \mu_{\mathcal{N}}^2(z')|},
\label{eq:TW3}
\end{equation}
evaluated at the source location $z = 0$, and the amplitude factor
\begin{align}
 Q_{\mathcal{N}}(\om,z) = \frac{\big|3 \phi_{\mathcal{N}}(\om,z)/2\big|^{1/6}}{\big|k^2 -
    \mu_{\mathcal{N}}^2(z)\big|^{1/4}},
\label{eq:TW4}
\end{align} 
is shown in \cite[Section 3.1.1]{borcea2016turning} to be bounded, and
at least twice continuously differentiable. The Airy functions $A_i$
and $B_i$ are evaluated at
\begin{equation}
\eta_{\mathcal{N}}^\ep(\om,z) = \mbox{sign} \Big(z - z_{-}^{(t)}\Big) \left[\frac{3
    |\phi_{\mathcal{N}}(\om,z)|}{2 \ep}\right]^{2/3},
\label{eq:TW5}
\end{equation}
where $|\eta^\ep(\om,z)|$ is of order one in the vicinity $ \big| z -
z_{-}^{(t)} \big| \le O \big(\ep^{2/3}\big)$ of the turning point, and
it is much larger than one in the rest of the domain
$\big(z_{-}^{(t)}-\delta,z_{-}^{(t-1)}\big)$.

We recall from \cite[Lemma 3.1]{borcea2016turning} that the matrix
${\bf M}_{\mathcal{N}}^\ep(\om,z)$ is invertible, with constant
determinant
\begin{equation}
\mbox{det}\, {\bf M}_{\mathcal{N}}^\ep(\om,z) = 2, \qquad \forall \, z \in
\big(z_{-}^{(t)}-\delta,z_{-}^{(t-1)}\big),
\label{eq:TW7}
\end{equation}
so the decomposition \eqref{eq:FM5} is well defined. Moreover,
\cite[Lemma 3.2]{borcea2016turning} shows that at $z - z_{-}^{(t)} \gg
\ep^{2/3}$ the expressions \eqref{eq:TW1}-\eqref{eq:TW2} become like
\eqref{eq:FM8},
\begin{align}
M_{\mathcal{N},11}^\ep(\om,z) &=
\frac{1}{\sqrt{\beta_{\mathcal{N}}(\om,z)}} \exp\Big[\frac{i}{\ep}
  \int_{0}^z dz' \beta_{\mathcal{N}}(\om,z')\Big] +
O(\ep), \nonumber \\ M_{\mathcal{N},21}^\ep(\om,z) &=
\beta_{\mathcal{N}}(\om,z) M_{\mathcal{N},11}^\ep(\om,z) +
O(\ep), \label{eq:FM8E}
\end{align}
so the turning wave behaves just like any other propagating wave until
it reaches the vicinity of the turning point from the right. On the
left side of the turning point, at $ z_{-}^{(t)} - z \gg \ep^{2/3}$,
the entries of ${\bf M}_{\mathcal{N}}^\ep(\om,z)$ grow exponentially,
as given in \cite[Lemma 3.3]{borcea2016turning}. The wave is
evanescent in this region, and must be decaying in order to have
energy conservation. This is ensured by the radiation condition
\begin{equation}
a_{\mathcal{N}}^\ep\big(\om,z_{-}^{(t)}-\delta\big) = i \exp\Big[
  \frac{2i}{\ep} \phi_{\mathcal{N}}\big(\om,0\big) \Big]
b_{\mathcal{N}}^\ep\big(\om,z_{-}^{(t)}-\delta\big),
\label{eq:TW8}
\end{equation}
which sets to zero the coefficients of the growing Airy function $B_i$
and its derivative $B_i'$ in the expression of $u_{\mathcal{N}}^\ep$
and $\partial_z u_{\mathcal{N}}^\ep$ at the end $z_{-}^{(t)}-\delta$ of the domain. We refer to \cite[Section 3.1]{borcea2016turning} for
more details, and for the proof that the result does not depend on the
particular choice of $\delta$ which is small, but
bounded away from $0$ in the limit $\ep \to 0$.

The evolution equation of the turning mode amplitudes is of
the same form as in \eqref{eq:FM12}, with the  following entries of the matrix
\eqref{eq:res24}-\eqref{eq:res24p} indexed by $j = \mathcal{N}$,
\begin{equation}
H_{\mathcal{N}}^{\ep(aa)}(\om,z) \approx
\frac{i\big|M_{\mathcal{N},11}^\ep(\om,z)\big|^2}{2} \Big[
  \frac{\sigma}{\sqrt{\ep}} \mu_j^2(z) \nu \Big(\frac{z}{\ep} \Big) +
  \sigma^2 g_j^\ep(\om,z)\Big],
\label{eq:TW9}
\end{equation}
and 
\begin{equation}
H_{\mathcal{N}}^{\ep(ab)}(\om,z) \approx -
\frac{i\big[M_{\mathcal{N},11}^\ep(\om,z)\big]^2}{2} \Big[
  \frac{\sigma}{\sqrt{\ep}} \mu_j^2(z) \nu \Big(\frac{z}{\ep} \Big) +
  \sigma^2 g_j^\ep(\om,z)\Big].
\label{eq:TW10}
\end{equation}
 These expressions reduce to those in \eqref{eq:FM14}-\eqref{eq:FM14p}
 at $z - z_{-}^{(t)} \gg \ep^{2/3}$, with the extra term involving
 $\partial_z\beta_{\mathcal{N}}$ in \eqref{eq:FM14p} coming from an
 $O(\ep)$ correction of the amplitudes, induced by the residual in
 \eqref{eq:FM8E}.

\subsubsection{Coupling with the evanescent waves}
\label{sect:TP3}
The modes indexed by $j > \mathcal{N}$ in equations \eqref{eq:DC3} are
evanescent waves, with wavenumber 
\begin{equation}
\beta_j(\om,z) = \sqrt{\mu_j^2(z)-k^2}.
\label{eq:BetaEv}
\end{equation}
We analyze these waves in appendix \ref{ap:evanesc}, and show that
they can be expressed in terms of the propagating ones. Explicitly, at
arc length $z \in \big(z_{-}^{(t)},z_{-}^{(t-1)}\big)$, satisfying
$z_{-}^{(t-1)} - z \gg \ep$, we obtain that 
\begin{align}
u_j^\ep(\om,z) \approx &-\frac{\sigma \sqrt{\ep}}{2 \beta_j(\om,z)}
\sum_{q=1}^{\mathcal{N}}\int_{-\infty}^\infty \hspace{-0.05in} d \xi \hspace{-0.03in} \left[
\gamma_{jq}^{(e)}\Big(\om,\frac{z}{\ep}+ \xi\Big) \frac{a_q^\ep(\om,z)}{\sqrt{\beta_q(\om,z)}} e^{
  \frac{1}{\ep}\int_{0}^z d z' \beta_q(\om,z') + i \xi
  \beta_q (\om,z) }  \right. \nonumber \\
&\left.-\overline{\gamma_{jq}^{(e)}}\Big(\om,\frac{z}{\ep}+ \xi\Big)
\frac{b_q^\ep(\om,z)}{\sqrt{\beta_q(\om,z)}} e^{-\frac{1}{\ep}
  \int_{0}^z d z' \beta_q(\om,z') - i \xi
  \beta_q(\om,z)} \right]e^{-|\xi|
  \beta_j(\om,z)} ,\label{eq:CE2}
\end{align}
where the approximation means that we neglect the terms that are
negligible in the limit $\ep \to 0$.  Here we introduced the notation
\begin{equation}
\gamma_{jq}^{(e)}\Big(\om,\frac{z}{\ep}+ \xi\Big) = \Gamma_{jq} \nu''
\Big(\frac{z}{\ep} + \xi\Big) + i \beta_q(\om,z) \Theta_{jq}
\nu'\Big(\frac{z}{\ep} + \xi)\Big),
\label{eq:CE2p}
 \end{equation}
with coefficients $\Gamma_{jq}$ and $\Theta_{jq}$ defined in
\eqref{eq:DC6}, and recall that the bar denotes complex conjugate.
For the derivative we have
\begin{align}
\hspace{-0.05in}\ep \partial_zu_j^\ep(\om,z) \approx &-\frac{\sigma \sqrt{\ep}}{2
  \beta_j(\om,z)} \sum_{q=1}^{\mathcal{N}}\int_{-\infty}^\infty \hspace{-0.05in} d \xi
\hspace{-0.03in} \left[\theta_{jq}^{(e)}\Big(\om,\frac{z}{\ep}+
\xi\Big)\frac{a_q^\ep(\om,z)}{\sqrt{\beta_q(\om,z)}} e^{
  \frac{1}{\ep}\int_{0}^z d z' \beta_q(\om,z') + i \xi
  \beta_q (\om,z)}\right. \nonumber \\ & \hspace{-0.2in}\left.  -
\overline{\theta_{jq}^{(e)}}\Big(\om,\frac{z}{\ep}+ \xi\Big)
  \frac{b_q^\ep(\om,z)}{\sqrt{\beta_q(\om,z)}} e^{
    -\frac{1}{\ep} \int_{0}^z d z' \beta_q(\om,z') - i \xi
    \beta_q(\om,z)} \right]e^{-|\xi|
  \beta_j(\om,z)} ,\label{eq:CE4}
\end{align}
with notation
\begin{align}
\theta_{jq}^{(e)}\Big(\om,\frac{z}{\ep}+ \xi\Big) = \Gamma_{jq} \nu'''
\Big(\frac{z}{\ep} + \xi\Big) - \beta_q^2(\om,z) \Theta_{jq}
\nu'\Big(\frac{z}{\ep} + \xi\Big)\nonumber \\+ i \beta_q(\om,z)
(\Gamma_{jq} + \Theta_{jq}) \nu''\Big(\frac{z}{\ep} + \xi)\Big)
.\label{eq:CE5}
\end{align}

\subsection{Closed system for the propagating  modes}
\label{sect:dec3}
The propagating mode amplitudes satisfy the system of equations
\eqref{eq:FM12}, with coupling modeled by the series
\eqref{eq:DC5}. Substituting the expressions \eqref{eq:CE2} and
\eqref{eq:CE4} of the evanescent waves in \eqref{eq:DC5}, we obtain a
closed system of equations for the propagating modes, as we now
explain.

\subsubsection{Propagation between turning points}
We begin with  $z \in
\big(z_{-}^{(t)},z_{-}^{(t-1)}\big)$ satisfying $z - z_{-}^{(t)} \gg
\ep^{2/3}$ and $z_{-}^{(t-1)} - z \gg \ep$. In this region the turning
wave indexed by $j = \mathcal{N}$ behaves like all the other
propagating modes, and the evanescent modes have the expression
\eqref{eq:CE2} and \eqref{eq:CE4}. The system of equations for the
right and left going amplitudes is
\begin{equation}
\partial_z \begin{pmatrix} \bm{a}^\ep(\om,z)
  \\ \bm{b}^\ep(\om,z) \end{pmatrix} = \bm{\Upsilon}^\ep(\om,z) 
\begin{pmatrix} \bm{a}^\ep(\om,z)
  \\ \bm{b}^\ep(\om,z) \end{pmatrix},
\label{eq:CS1}
\end{equation}
where $\bm{a}^\ep$ and $\bm{b}^\ep$ are the complex column vectors in
$\mathbb{C}^{\mathcal{N}}$ with entries $a_j^\ep$ and $b_j^\ep$, for
$1 \le j \le \mathcal{N}$. The complex matrix
$\bm{\Upsilon}^\ep(\om,z)$ depends on the random fluctuations $\nu$
and the slow changes of the waveguide, and has the block structure
\begin{equation} 
\bm{\Upsilon}^\ep(\om,z) = 
\begin{pmatrix}
  \bm{\Upsilon}^{\ep(aa)}(\om,z) & \bm{\Upsilon}^{\ep (ab)} (\om,z)
  \\ \bm{\Upsilon}^{\ep(ba)}(\om,z) & \bm{\Upsilon}^{\ep (bb)}
    (\om,z) \end{pmatrix},
\label{eq:CS1p}
\end{equation}
with $\mathcal{N}\times \mathcal{N}$ 
blocks  satisfying the relations
\begin{equation}
\bm{\Upsilon}^{\ep(ba)}(\om,z) =
\overline{\bm{\Upsilon}^{\ep(ab)}}(\om,z), \quad
\bm{\Upsilon}^{\ep(bb)}(\om,z) =
\overline{\bm{\Upsilon}^{\ep(aa)}}(\om,z).
\label{eq:CS2}
\end{equation}
Their entries are defined as follows: Off the diagonal, we have
\begin{align}
&\Upsilon_{jq}^{\ep(aa)}(\om,z) = -\frac{ie^{\frac{i}{\ep}
  \int_{0}^z dz' \,
  \big(\beta_q(\om,z')-\beta_j(\om,z')\big)}}{2
  \sqrt{\beta_j(\om,z)\beta_q(\om,z)}} \left\{
\frac{\sigma}{\sqrt{\ep}}\Big[\Gamma_{jq}
  \nu''\Big(\frac{z}{\ep}\Big)+ i \beta_q(\om,z) \Theta_{jq}
  \nu'\Big(\frac{z}{\ep}\Big)\Big] \right. \nonumber \\ &\hspace{0.2in}\left. +
\sigma^2 \Big[\widetilde{\gamma}_{jq}\Big(\om,\frac{z}{\ep}\Big) + i
  \beta_q \widetilde{\theta}_{jq}\Big(\om,\frac{z}{\ep}\Big) \Big] +
\gamma_{jq}^{o}(z) + i \beta_q(\om,z)
\theta_{jq}^{o}(z)\right\}, \qquad j \ne q,
\label{eq:CS4}
\end{align}
and 
\begin{align}
\Upsilon_{jq}^{\ep(ab)}(\om,z) =
\overline{\Upsilon_{jq}^{\ep(aa)}}(\om,z) e^{ - \frac{2i}{\ep}
  \int_{0}^z dz' \, \beta_j(\om,z')}, \qquad j \ne q,
\label{eq:CS4p}
\end{align}
and on the diagonal we have 
\begin{equation}
\Upsilon_{jj}^{\ep(aa)}(\om,z) = H_{j}^{\ep(aa)} (\om,z) + \frac{i\sigma^2}{2
 \beta_j(\om,z)} \eta_j\Big(\om,\frac{z}{\ep}\Big),
\label{eq:CS5}
\end{equation}
and 
\begin{align}
\Upsilon_{jj}^{\ep(ab)}(\om,z) =
\Big[\overline{\Upsilon_{jj}^{\ep(aa)}}(\om,z)-\frac{\partial_z
    \beta_j(\om,z)}{2 \beta_j(\om,z)}\Big] e^{ - \frac{2i}{\ep}
  \int_{0}^z dz' \, \beta_j(\om,z')}.
\label{eq:CS5p}
\end{align}
The coefficients in these definitions are given in \eqref{eq:FM14},
and \eqref{eq:DC6}-\eqref{eq:DC10}, except for $\eta_j$,
$\widetilde{\gamma}_{jq}$ and $\widetilde{\theta}_{jq}$, which include
the interaction with the evanescent modes. These are defined by
\begin{align*}
\widetilde{\gamma}_{jq}\Big(\om,\frac{z}{\ep}\Big) =
          {\gamma}_{jq}\Big(\frac{z}{\ep}\Big) - \sum_{l >
            \mathcal{N}} \frac{\Gamma_{jl}}{2 \beta_l(\om,z)}
          \nu''\Big(\frac{z}{\ep}\Big) \int_{-\infty}^\infty \hspace{-0.05in}d \xi \,
          \gamma_{lq}^{(e)} \Big(\om,\frac{z}{\ep}+ \xi\Big) e^{-|\xi| \beta_l(\om,z)+i \xi
            \beta_q(\om,z)}, 
\end{align*}
and 
\begin{align*}
\widetilde{\theta}_{jq}\Big(\om,\frac{z}{\ep}\Big) =
          {\theta}_{jq}\Big(\frac{z}{\ep}\Big) - \sum_{l >
            \mathcal{N}} \frac{\Theta_{jl}}{2 \beta_l(\om,z)}
          \nu'\Big(\frac{z}{\ep}\Big) \int_{-\infty}^\infty \hspace{-0.05in}d \xi \,
          \theta_{lq}^{(e)} \Big(\om,\frac{z}{\ep}+ \xi\Big) e^{-|\xi| \beta_l(\om,z)+i \xi
            \beta_q(\om,z)}, 
\end{align*}
and 
\begin{align*}
\eta_j\Big(\om,\frac{z}{\ep}\Big) = &\sum_{l > \mathcal{N}} \frac{1}{2
  \beta_l(\om,z)} \int_{-\infty}^\infty \hspace{-0.05in} d \xi \, e^{-|\xi|
  \beta_l(\om,z)+i \xi \beta_j(\om,z)} \nonumber
\\ &\times \Big[\Gamma_{jl} \nu''\Big(\frac{z}{\ep}\Big)
  \gamma_{lj}^{(e)}\Big(\om,\frac{z}{\ep}+ \xi\Big) + \Theta_{jl}
  \nu'\Big(\frac{z}{\ep}\Big) \theta_{lj}^{(e)}\Big(\om,\frac{z}{\ep}+
  \xi\Big)\Big], 
\end{align*}
with $\gamma_{jq}$ and $\theta_{jq}$ given in
\eqref{eq:DC7}-\eqref{eq:DC8} and $\gamma_{lq}^{(e)}$,
$\theta_{lq}^{(e)}$ given in \eqref{eq:CE2p} and \eqref{eq:CE5}.  Note
that the coefficients $\Gamma_{jl}/\beta_l$ and $\Theta_{jl}/\beta_l$
decay as $1/l^2$ for $l \gg 1$, and the integrals in $\xi$ are
bounded, so the series defining $\widetilde{\gamma}_{jq}$, $\widetilde{\theta}_{jq}$ and $\eta_j$ are absolutely
convergent.

\subsubsection{Vicinity of turning points}
\label{sect:victurn}
Let us consider a vicinity $|z-z_{-}^{(t)}| = O(\ep^s)$ of the turning
point $z_{-}^{(t)}$, for some $s > 0$, and change for a moment
variables to $z = z_{-}^{(t)} + \ep^s \zeta$, so that $\zeta =
O(1)$. In the new variable, we observe that the  coupling terms
in the evolution equations \eqref{eq:FM12} for the turning wave
indexed by $j = \mathcal{N}$, modeled by the series \eqref{eq:DC5},
are proportional to
\begin{equation}
\frac{\ep^{s/2}}{\sqrt{\ep^{1-s}}}
\widetilde{\nu}\Big(\frac{\zeta}{\ep^{1-s}}\Big) + O(\ep^s), \quad \widetilde \nu
= \nu'' ~ \mbox{or} ~ \nu'.
\end{equation}
In the limit $\ep \to 0$, described in detail in section
\ref{sect:results}, all these terms tend to zero. Thus, the turning
wave does not interact with the other modes near the turning point.

We also obtain that the right hand side of equation \eqref{eq:FM12}
for $1 \le j \le \mathcal{N}-1$ tends to zero as $\ep \to 0$, so the
propagating mode amplitudes remain constant as they traverse the
vicinity of the turning point $z_{-}^{(t)}$.  A similar argument shows
that the propagating mode amplitudes remain constant as they traverse
the vicinity of the turning point $z_{-}^{(t-1)}$ at the other end of
the interval.

It remains to describe the turning mode that undergoes a transition
near $z_{-}^{(t)}$. To do so, we return to the original coordinate
$z$, but stay in the vicinity of $z_{-}^{(t)}$. We obtain from
\eqref{eq:FM12} with $j = \mathcal{N}$, after neglecting the coupling
terms, that
\begin{align} 
\partial_z \begin{pmatrix} a_{\mathcal{N}}^\ep(\om,z)
  \\ b_{\mathcal{N}}^\ep(\om,z) \end{pmatrix} \approx
{\bf H}_{\mathcal{N}}^\ep(\om,z) \begin{pmatrix}
  a_{\mathcal{N}}^\ep(\om,z)
  \\ b_{\mathcal{N}}^\ep(\om,z) \end{pmatrix},
\label{eq:CS9}
\end{align}
where the matrix ${\bf H}_{\mathcal{N}}^\ep$ is defined by
\eqref{eq:res24} and \eqref{eq:TW9}-\eqref{eq:TW10}. These equations
give 
\begin{equation}
\partial_z \Big[ \big|a_{\mathcal{N}}^\ep(\om,z)\big|^2 -
  \big|b_{\mathcal{N}}^\ep(\om,z)\big|^2\Big] \approx 0,
\end{equation}
and using the radiation condition \eqref{eq:TW8}, we conclude that
near the turning point we have energy conservation\footnote{All the energy
  conservation relations are approximate at a finite $\ep$, due to the
  interaction with the evanescent modes. We will see in section
  \ref{sect:results} that there is no energy loss in the limit $\ep  \to 0$.}
\begin{equation}
\big|a_{\mathcal{N}}^\ep(\om,z)\big|^2 \approx
\big|b_{\mathcal{N}}^\ep(\om,z)\big|^2.
\label{eq:CS10}
\end{equation}
The impinging left going wave with amplitude $b^\ep$ is reflected back
at the turning point to give the right going wave with amplitude
$a^\ep$, determined by the reflection coefficient
\begin{equation}
{\rm R}^{\ep}_{\mathcal{N}}(\om,z) =
\frac{a_{\mathcal{N}}^\ep(\om,z)}{b_{\mathcal{N}}^\ep(\om,z)} \approx
i \exp\big[ \frac{2 i}{\ep} \phi_{\mathcal{N}}(\om,0) + i
  \vartheta_{\mathcal{N}}^\ep(\om,z)\big].
\label{eq:CS11}
\end{equation}
This is a complex number with modulus $\big|{\rm
  R}^{\ep}_{\mathcal{N}}(\om,z)\big| \approx 1$, because there is no
loss of power in the limit $\ep \to 0$, and with random phase $
\vartheta_{\mathcal{N}}^\ep(\om,z).  $

The phase $\vartheta_{\mathcal{N}}^\ep$ is described in detail in
\cite[Lemmas 4.1 and 4.2]{borcea2016turning}, for the purpose of
characterizing the reflection of a broad-band pulse at the turning
point. The standard deviation of the random boundary fluctuations
considered in \cite{borcea2016turning} is smaller than what we have in
\eqref{eq:sc4}, by a factor of $|\ln \ep|^{1/2}$, so that as $\ep \to 0$
there is no mode coupling at any $z$, small or order one. Here we have
mode coupling away from the turning points, due to the stronger random
boundary fluctuations, and we are interested in the transport of
energy by single frequency modes in the waveguide. The mode powers are
not affected by the phase, so the details of
$\vartheta_{\mathcal{N}}^\ep(\om,z)$ are not important in the context of
this paper.

\subsubsection{Source excitation and matching conditions} 
The evolution equations of the left and right going mode amplitudes,
described above, are complemented by matching conditions at the
turning points, by radiation conditions at $|z| > Z_M$, and by 
initial conditions at $z=0$, where the source lies.

Starting from the source location $z=0$, which  is not a
turning point, we have the jump conditions,
\begin{align}
a_{j}^\ep(\om,0+)-a_{j}^\ep(\om,0-) &= \frac{\hat f(\om)
  y_j(\rho_\star,0)}{2 i \sqrt{\beta_j(\om,0)}}, \nonumber \\
   b_{j}^\ep(\om,0+)-b_{j}^\ep(\om,0-) &= \frac{\hat f(\om)
  y_j(\rho_\star,0)}{2 i \sqrt{\beta_j(\om,0)}},  \quad 1 \le j \le  N^{(0)},
\label{eq:CS14}
\end{align}
where we recall that $N^{(0)} $ is  the number of propagating modes
at $z = 0$ and we denote $a(0+) = \lim_{z \searrow 0} a(z)$ and $a(0-) = \lim_{z \nearrow 0} a(z)$.

On the left of the source, at turning points $z_{-}^{(t)}$, for $1 \le t \le 
t_M^{-}$, we have the continuity relations
\begin{align}
\hspace{-0.1in} a_{j}^\ep(\om,z_{-}^{(t)}+)= a_{j}^\ep(\om,z_{-}^{(t)}-),\quad 
b_{j}^\ep(\om,z_{-}^{(t)}+)&= b_{j}^\ep(\om,z_{-}^{(t)}-),  \label{eq:CS15}
\end{align}
for $1 \le j  \le N_{-}^{(t-1)}-1$, where we recall definition
\eqref{eq:DC13pp} of $N_{-}^{(t-1)}$. We also have the reflection of
the turning mode, modeled by equation 
\begin{align}
a_{N_{-}^{(t-1)}}^\ep(\om,z_{-}^{(t)}+) =
{\rm R}^{\ep}_{N_{-}^{(t-1)}}(\om,z_{-}^{(t)})
b_{N_{-}^{(t-1)}}^\ep(\om,z_{-}^{(t)}+),
\label{eq:CS16}
\end{align}
where ${\rm R}^{\ep}_{N_{-}^{(t-1)}}$ is the complex reflection
coefficient defined as in \eqref{eq:CS11}.

At $z < -Z_M$, where the waveguide has straight and parallel
boundaries and supports $N_{\rm min}$ propagating modes, the wave is
outgoing (propagating to the left), so we have the conditions
\begin{align}
\hspace{-0.1in}a_{j}(z) &= a_{j}(-Z_M+) = 0, \quad b_{j}(z) =
b_{j}(-Z_M+), \qquad z < -Z_M,
\label{eq:CS17}
\end{align} 
for $j = 1, \ldots, N_{\rm min}.$

Similarly, on the right of the source, at turning points
$z_{+}^{(t)}$, for $1\le t \le t_M^+$, we have the continuity
relations
\begin{align}
a_{j}^\ep(\om,z_{+}^{(t)}+)= a_{j}^\ep(\om,z_{+}^{(t)}-), \quad 
b_{j}^\ep(\om,z_{+}^{(t)}+)&= b_{j}^\ep(\om,z_{+}^{(t)}-), 
\label{eq:CS15P}
\end{align}
for $1\le j \le N_{+}^{(t-1)}$, where we recall definition
\eqref{eq:DC13pp} of $N_{+}^{(t-1)}$. The number of propagating modes
increases by one at $z_{+}^{(t)}$, to equal $ N_{+}^{(t)}$, and the
amplitude of the turning wave, indexed by $j =N_{+}^{(t)}$, starts
from zero there
\begin{align}
a_{N_{+}^{(t)}}^\ep(\om,z_{+}^{(t)}) =
b_{N_{+}^{(t)}}^\ep(\om,z_{+}^{(t)}+) = 0.
\label{eq:CS16P}
\end{align}

At $z > Z_M$, where the waveguide has straight and parallel boundaries
and supports $N_{\rm max}$ propagating modes, the wave is outgoing
(propagating to the right), so we have the conditions
\begin{align}
a_{j}(z) = a_{j}(Z_M-), \quad  b_{j}(z) &= b_{j}(Z_M-) = 0, \qquad z > Z_M,
\label{eq:CS17P}
\end{align} 
for $j = 1, \ldots, N_{\rm max}. $

\section{The asymptotic limit}
\label{sect:results}
To quantify the net effect of the waveguide variations on the
propagating waves, we take the asymptotic limit $\ep \to 0$ of the
random  mode amplitudes. The limit is taken in
each sector of the waveguide, bounded by two consecutive turning
points, as explained in section \ref{sect:CSC1}.  We introduce in
section \ref{sect:forwsc} a simplification, known as the forward
scattering approximation, which applies to smooth enough random
fluctuations $\nu$. The $\ep \to 0$ limit of the mode amplitudes under
this approximation is described in section \ref{sect:coupledpower}.

\subsection{The propagator matrix}
\label{sect:CSC1}
The discussion below applies to any sector of the waveguide, so let us
consider as in section \ref{sect:TP1} the sector $z \in
\big(z_{-}^{(t)},z_{-}^{(t-1)}\big)$, supporting ${\cal
  N}=N_-^{(t-1)}$ propagating modes. 

The mode amplitudes satisfy the system of equations \eqref{eq:CS1},
with $2 \mathcal{N} \times 2 \mathcal{N}$ random propagator matrix
$\calP^\ep(\om,z;z_{-}^{(t-1)})$. This solves the equation
\begin{equation}
\partial_z \calP^\ep(\om,z;z_{-}^{(t-1)})=
\bm{\Upsilon}^{\ep}(\om,z) \calP^\ep(\om,z;z_{-}^{(t-1)}),
    \label{eq:propagator1}
\end{equation}
backward in $z$, starting from
\begin{equation}
\calP^\ep(\om,z_{-}^{(t-1)};z_{-}^{(t-1)} )={\bf I},
  \label{eq:propagator1p}
\end{equation}
where ${\bf I}$ is the $2 \mathcal{N} \times 2 \mathcal{N}$ identity
matrix and $\bm{\Upsilon}^{\ep}(\om,z)$ is defined in
\eqref{eq:CS1p}-\eqref{eq:CS5p}. 

The propagator defines the solution of \eqref{eq:CS1},
\begin{equation}
\begin{pmatrix} {\itbf a}^\ep(\om,z) \\  {\itbf b}^\ep(\om,z)
\end{pmatrix} = \calP^\ep(\om,z;z_{-}^{(t-1)} ) 
\begin{pmatrix} {\itbf a}^\ep(\om,z_{-}^{(t-1)}) 
\\ {\itbf b}^\ep(\om,z_{-}^{(t-1)})
\end{pmatrix}, 
\label{eq:propagator1pp}
\end{equation}
and due to the symmetry relations \eqref{eq:CS2} of the blocks of
$\bm{\Upsilon}^\ep$, we note that 
\begin{equation}
\begin{pmatrix} \overline{{\itbf b}^\ep}(\om,z) \\  \overline{{\itbf a}^\ep}(\om,z)
\end{pmatrix} = \calP^\ep(\om,z;z_{-}^{(t-1)} )  \begin{pmatrix}
 \overline{{\itbf b}^\ep}(\om,z_{-}^{(t-1)}) \\ \overline{{\itbf a}^\ep}(\om,z_{-}^{(t-1)})
\end{pmatrix}
\label{eq:propagator1ppp}
\end{equation}
is also a solution.  Writing explicitly these equations, and using the
uniqueness of solution of \eqref{eq:CS1}, we conclude that the
propagator has the block form
\begin{equation}
\label{formpropagator}
\calP^\ep(\om,z;z_{-}^{(t-1)} ) =\begin{pmatrix}
\overline{\calP^{\ep(bb)}}(\om,z;z_{-}^{(t-1)}) &
\overline{\calP^{\ep(ba)}}(\om,z;z_{-}^{(t-1)})\\ 
\calP^{\ep(ba)}(\om,z;z_{-}^{(t-1)})
& \calP^{\ep(bb)}(\om,z;z_{-}^{(t-1)})\\
\end{pmatrix} \, .
\end{equation}
The blocks are $\mathcal{N} \times \mathcal{N}$ complex matrices,
where $\calP^{\ep(bb)}$ describes the coupling between the
components of ${\itbf b}^\ep$, the vector of left-going mode amplitudes,
and $\calP^{\ep(ba)}$ describes the coupling between the
components of ${\itbf b}^{\ep}$ and of ${\itbf a}^\ep$, the vector of
right-going mode amplitudes.

The limit of $\calP^{\ep}(\om,z;z_{-}^{(t-1)})$ as $\ep
\rightarrow 0$ can be obtained and identified as a multi-dimensional
diffusion process, meaning that the entries of the limit matrix
satisfy a system of linear stochastic equations.  This follows from
the application of an extension of the diffusion approximation theorem
proved in \cite{papanicolaou1974asymptotic} and presented in
\cite[Chapter 6]{fouque07}. This extension is stated in Theorem
\ref{prop.diflim} and is proved in section \ref{app:adif} for a
general system of  equations with real valued unknown
vector ${\bm{X}}^\ep$. In our case ${\bm{X}}^\ep$ is obtained by
concatenating the moduli and arguments of the entries in
$\calP^{\ep(bb)}$ and $\calP^{\ep(ba)}$.

\subsection{The forward scattering approximation}
\label{sect:forwsc}
When we use Theorem \ref{prop.diflim}, we obtain that the limit
entries of $\calP^{\ep(bb)}$ are coupled to the limit
entries of $\calP^{\ep(ba)}$ through coefficients that are
proportional to the power spectral density\footnote{The power spectral
  density is the Fourier transform of the auto-correlation function $\cR$
  defined in \eqref{eq:rp2}. It is a non-negative and even function
  that decays rapidly when $\cR$ and therefore $\nu$ are smooth in
  $z$.} $\hat \cR$ of the random fluctuations $\nu$, evaluated at the
sum of the mode wavenumbers,
\begin{equation}\hat \cR\big(\beta_j (\om,z)+ \beta_l(\om,z)\big) = 
2 \int_{0}^\infty d\zeta \, \cR (\zeta) \cos[( \beta_j(\om,z)+\beta_l
  (\om,z))\zeta]\, ,
\label{eq:pspsum}
\end{equation}
for $j, l = 1, \ldots, \mathcal{N}$. This can be traced back to the
phase factors 
\[
\frac{1}{\ep} \int_{0}^z d z'
  \big[\beta_j(\om,z')+\beta_l(\om,z')\big]
\] 
in matrix $\bm{\Upsilon}^{\ep (ba)} (\om,z)$ defined in
\eqref{eq:CS4p}. The limit entries of $\calP^{\ep(bb)}(z)$ are
coupled to each other through $\hat \cR\big(\beta_j (\om,z)-
\beta_l(\om,z)\big)$, because the
phase factors in $\bm{\Upsilon}^{\ep (bb)} (\om,z)$ defined in
\eqref{eq:CS2}-\eqref{eq:CS4} are
\[
\frac{1}{\ep} \int_{0}^z d z'
\big[\beta_j(\om,z')-\beta_l(\om,z')\big],
\] 
for $j, l = 1, \ldots, \mathcal{N}$.

To simplify the analysis of the cumulative scattering effects in the
limit $\ep \to 0$, we assume that the power spectral density $\hat
\cR$ peaks at zero and decays rapidly away from it\footnote{An example
  is the Fourier transform of the Gaussian auto-correlation function used in the
  numerical simulations in section \ref{sect:netscat}.}, so that 
\begin{equation}
\hat \cR\big(\beta_j (\om,z)+ \beta_l(\om,z)\big)  \approx 0\, , \quad 
\forall \, j,l = 1, \ldots, \mathcal{N}.
\label{eq:FSCCOND}
\end{equation}
With this assumption we can neglect the coupling between the blocks
$\calP^{\ep(bb)}(\om,z)$ and $\calP^{\ep(ba)}(\om,z)$ of the
propagator. Since $ \calP^{\ep(ba)}$ starts from zero at $z =
z_{-}^{(t-1)}$, we obtain
\begin{equation}
\label{eq:fwdprop}
\calP^\ep(\om,z;z_{-}^{(t-1)}) \approx \begin{pmatrix}
\overline{\calP^{\ep(bb)}}(\om,z;z_{-}^{(t-1)}) & {\bf 0}\\
{\bf 0} &  {\calP^{\ep(bb)}(\om,z;z_{-}^{(t-1)})}\\
\end{pmatrix} ,
\end{equation}
and equation \eqref{eq:propagator1pp} gives
\begin{equation}
{\itbf b}^\ep(\om,z) \approx
{\calP^{\ep(bb)}(\om,z;z_{-}^{(t-1)})}{\itbf b}^\ep(\om,z_{-}^{(t-1)}) \, , \quad z < z_{-}^{(t-1)}.
\label{evola}
\end{equation}
This is the forward scattering approximation. It describes the left
going amplitudes ${\itbf b}^\ep(\om,z)$ of the waves, propagating forward
from the source, independent of the right-going amplitudes ${\itbf a}^\ep(\om,z)$ of the waves, propagating backward, toward the source.

Note that since $\beta_j$ decrease monotonically with $j$, the
smallest argument of the power spectral density in \eqref{eq:FSCCOND}
is at $j = l = \mathcal{N}$.  The wave number $\beta_{\mathcal{N}}(z)$
is of order $k/\sqrt{\mathcal{N}}$ away from the turning point
$z_{-}^{(t)}$, but tends to zero as $z \searrow z_{-}^{(t)}.$ The
left and right going amplitudes of the turning mode are coupled
near $z_{-}^{(t)}$, as described by the reflection coefficient
\eqref{eq:CS11}. We assume that this coupling holds only at $z - z_{-}^{(t)}
< \delta$, where $\delta$ is a small and positive number,
independent of $\ep$. Over the small distance $\delta$ there is
negligible interaction between the turning mode and the others, as
explained in section \ref{sect:victurn}. In the remaining interval $
z \in (z_{-}^{(t)}+\delta, z_{-}^{(t-1)})$ we have
\begin{equation}
\cR\big(2\beta_{\mathcal{N}} (\om,z)\big) \lesssim \hat
\cR\big(2\beta_{\mathcal{N}} (\om,z_{-}^{(t)} + \delta)\big) \approx
0,
\end{equation}
so we can use the forward scattering approximation.

\begin{figure}[t]
\begin{center}
\includegraphics[width=0.5\textwidth]{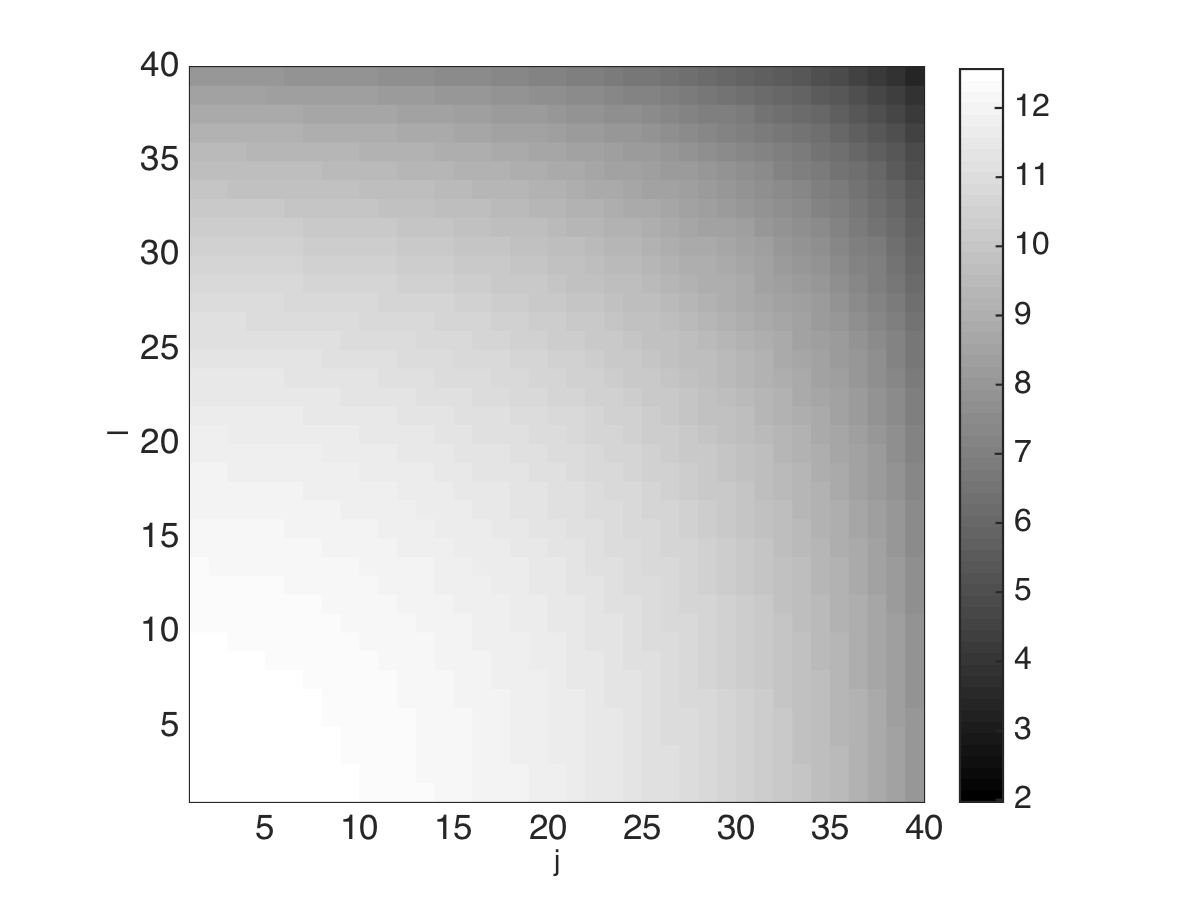}
\hspace{-0.3in}\includegraphics[width=0.5\textwidth]{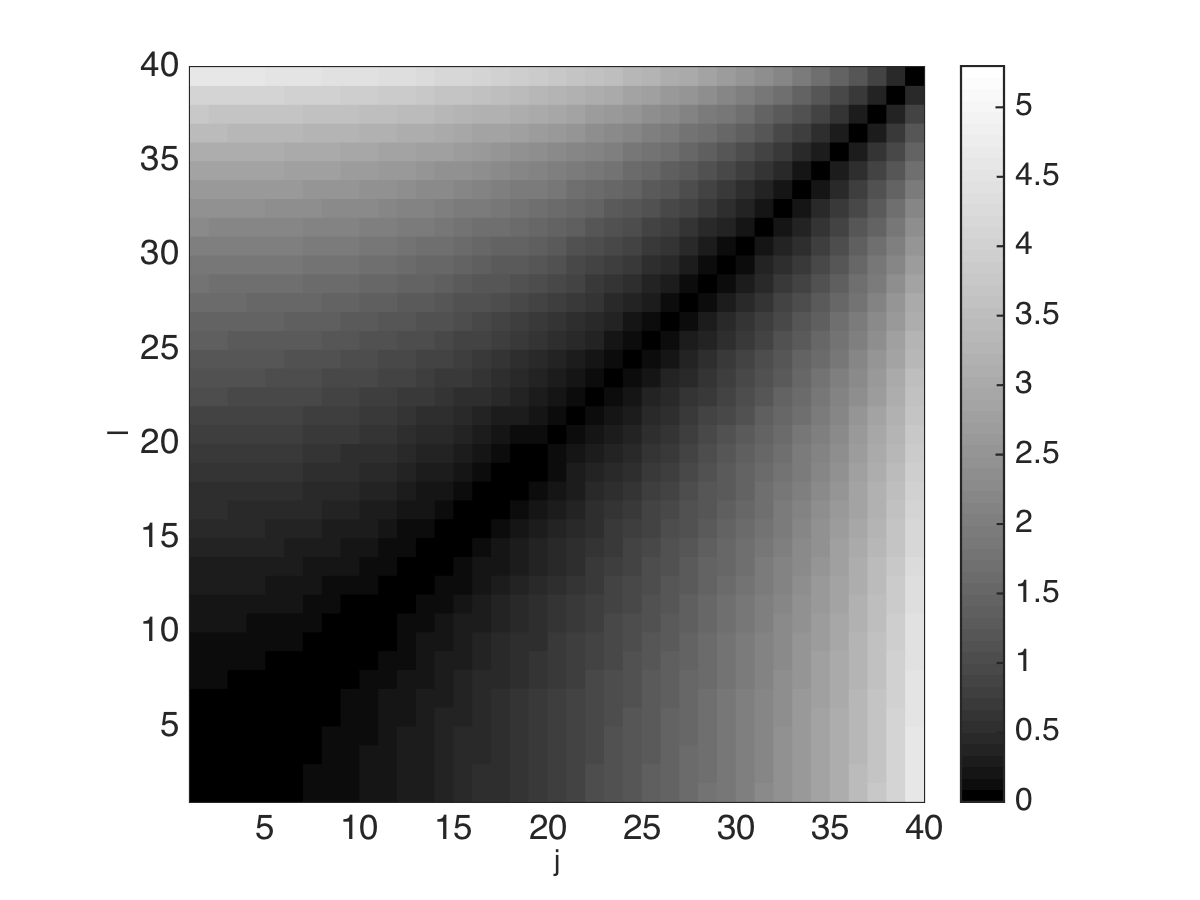}
\end{center}
\vspace{-0.1in}
\caption{Plot of the matrix with entries $\beta_j+\beta_l$ on the left
  and $|\beta_j-\beta_l|$ on the right, v.s. $j,l = 1, \ldots,
  \mathcal{N}$, for the case of $\mathcal{N} = 40$ propagating
  modes. The scaled wavenumber is $k = 2 \pi$ and the waveguide width
  is $D = 20.25$. Note that the entries in the left plot are larger
  than $2 \beta_\mathcal{N} = 1.97$, whereas the entries near the
  diagonal in the right plot are small.}
\label{fig:compbeta}
\end{figure}

Note that there is mode coupling in this approximation, but only
between the forward going mode amplitudes. This is due to the fact
that $|\beta_j(\om,z)-\beta_l(\om,z)|$ is small at least for nearby
indexes $j,l$, as illustrated in Figure \ref{fig:compbeta}. The power
spectral density evaluated at such differences is not negligible, and
the net coupling effect is described in the next section.
\subsection{The coupled mode diffusion process}
\label{sect:coupledpower}
The $\ep \to 0$ limit of the forward going mode amplitudes is stated
in the next theorem. We derive it using Theorem \ref{prop.diflim}
for the vector ${\itbf X}^\ep \in \mathbb{R}^{2 \mathcal{N}}$ obtained
by concatenating the moduli and arguments of $b_j^\ep$, with $j = 1,
\ldots, \mathcal{N}$.  The differential equations for ${\itbf X}^\ep$ follow from the system
\begin{equation}
\partial_z {\itbf b}^\ep(\om,z) \approx \bm{\Upsilon}^{\ep(bb)}(\om,z) {\itbf b}^\ep(\om,z)\, , \quad z < z_{-}^{(t-1)},
\label{eq:FORWEQ}
\end{equation}
with given ${\itbf b}^\ep(\om,z_{-}^{(t-1)})$. As explained in
the previous section, the approximation in \eqref{eq:FORWEQ} means
that there is an error that vanishes in the limit $\ep \to 0$.

\vspace{0.05in}
\begin{theorem}
\label{propdiff}%
The complex mode amplitudes $\{{b}_j^\ep(\omega,z) \}_{ j=1}^{ \mathcal{N}}$ converge in distribution as $\ep \rightarrow 0$ to an
inhomogeneous diffusion Markov process $\{{b}_j(\omega,z)
\}_{j=1}^{{\cal N}}$ with generator $- {\cal L}_z^{{\cal N}}$ given
below.\footnote{The minus sign in front of the generator is because we
  solve the Kolmogorov equation for the moments of the limit process
  backward in $z$, starting from $z_{-}^{(t-1)}$.}
\end{theorem}

\vspace{0.05in}
Let us write the limit process as
$$ 
{b}_j(\omega,z) = P_j^{1/2}(\omega,z) e^{i \psi_j(\omega,z)},
\quad j=1,\ldots,{\cal N},
$$ in terms of the power $P_j = | b_j|^2$ and the phase $\psi_j =
\mbox{arg} \, b_j$. Then, we can express the infinitesimal generator
of the limit diffusion as the sum of two operators
\begin{equation}
\label{gendiffa}
{\cal L}_z^{{\cal N}} = {\cal L}_{P,z}^{{\cal N}} + {\cal L}_{\psi,z}^{{\cal N}} . 
\end{equation}
The first is a partial differential operator in the powers
\begin{equation}
\label{gendiff2P}
{\cal L}_{P ,z}^{{\cal N}}= \hspace{-0.07in}\sum_{{\scriptsize \begin{array}{c}j, l = 1
      \\ j \ne l \end{array}} }^{{\cal N}} \hspace{-0.07in} G_{jl}^{(c)}(\omega,z)
\left[ P_l P_j \left( \frac{\partial}{\partial P_j}
  -\frac{\partial}{\partial P_l} \right) \frac{\partial}{\partial P_j}
  + (P_l-P_j) \frac{\partial}{\partial P_j} \right] \, ,
\end{equation}
with symmetric matrix ${\bf G}^{(c)}(\om,z) =
\big(G^{(c)}_{jl}(\om,z)\big)_{j,l=1}^{{\cal N}}$ of coefficients that
are non-negative off the diagonal
\begin{equation}
 G_{jl}^{(c)}(\omega,z) = \frac{\sigma^2
   \mu_j^2(z)\mu_l^2(z)}{4\beta_j(\om,z) \beta_l(\om,z)} \hat \cR
       [\beta_j(\om,z)-\beta_l(\om,z)] \, , \quad j\neq
       l\, \label{defgamma1},
\end{equation} 
and sum to zero in the rows
\begin{equation}
\label{defgamma1b}
G_{jj}^{(c)}(\omega,z) = - \sum_{l = 1, l \ne j}^{{\cal N}} G_{jl}^{(c)}(\omega,z)\, .
\end{equation}
The second partial differential operator in \eqref{gendiffa} is with
  respect to the phases
\begin{align}
\nonumber {\cal L}_{\psi,z}^{{\cal N}} &= \frac{1}{8}\hspace{-0.08in}
\sum_{{\scriptsize \begin{array}{c}j, l = 1 \\ j \ne
      l \end{array}}}^{{\cal N}}\hspace{-0.07in} G_{jl}^{(c)}(\omega,z) \left[
  \frac{P_j}{P_l} \frac{\partial^2}{ \partial \psi_l^2} +
  \frac{P_l}{P_j} \frac{\partial^2}{ \partial \psi_j^2} + 2
  \frac{\partial^2}{\partial \psi_j \partial \psi_l} \right] +
\frac{1}{2} \sum_{j , l=1}^{{\cal N}} G_{jl}^{(0)}(\omega,z)
\frac{\partial^2}{\partial \psi_j \partial \psi_l} \\ & +
\frac{1}{2} \hspace{-0.08in}\sum_{{\scriptsize \begin{array}{c}j, l = 1 \\ j \ne
      l \end{array}} }^{{\cal N}} \hspace{-0.07in} G_{jl}^{(s)}(\omega,z)
\frac{\partial}{\partial \psi_j} + \sum_{j=1}^{{\cal N}}
\kappa_j^{{\cal N}}(\om,z) \frac{\partial}{\partial \psi_j}
\,, \label{gendiff2T}
\end{align}
with coefficients
\begin{align}
G_{jl}^{(0)}(\omega,z) = \frac{\sigma^2
  \mu_j^2(z)\mu_l^2(z)}{4 \beta_j(\om,z) \beta_l(\om,z)} \hat \cR (0)  
\, \label{defgamma10} ,  \quad j,l = 1, \ldots, \mathcal{N},
\end{align} 
and
\begin{align}
 G_{jl}^{(s)}(\omega,z) = \frac{\sigma^2
   \mu_j^2(z)\mu_l^2(z)}{2\beta_j(\om,z) \beta_l(\om,z)}
 \int_{0}^\infty \hspace{-0.05in} d\zeta \, \cR (\zeta )\sin \left[
   (\beta_j(\omega,z)-\beta_l(\omega,z))\zeta \right]
 \, \label{defgamma1s}, \quad
\end{align} 
for $j,l = 1, \ldots, \mathcal{N}$ and $j \ne l$.  The coefficient
$\kappa_j^{{\cal N}}$ in the last term of (\ref{gendiff2T}) is
\begin{align}
\nonumber &\hspace{-0.2in} \kappa_j^{{\cal N}}(\om,z) =
\frac{\sigma^2}{2\beta_j(\om,z)} \Big\{ \Big(\frac{\pi^2
  j^2}{12}+\frac{1}{16}\Big) \cR''(0)-\frac{3 \mu_j^2(z)}{4}
\cR(0)\Big\} \\ \nonumber & - \sum_{{\scriptsize \begin{array}{c}l =
      1 \\ j \ne l \end{array}} }^{{\cal N}} \frac{\mu_j^2(z)
  \mu_l^2(z)}{4 \beta_j(\om,z)
  \beta_l(\om,z)[\beta_j(\om,z)-\beta_l(\om,z)]} \Big[ \cR(0)
  +\frac{\cR''(0)}{[\beta_j(\om,z)+\beta_l(\om,z)]^2}\Big] \\ \nonumber
&+\sum_{l > {\cal N}} \, \frac{\sigma^2 \mu_j^2(z)\mu_l^2(z)}{2\beta_j
  \beta_l [\beta_j^2(\om,z)+\beta_l^2(\om,z)]^2} \Big\{ -
\beta_l(\om,z) \cR''(0) + \int_0^\infty \hspace{-0.05in} d \zeta \,
\cR''(\zeta) e^{-\beta_l (\om,z) \zeta} \\ & \times \Big[
  [\beta_l^2(\om,z) - \beta_j^2(\om,z)] \cos(\beta_j (\om,z) \zeta) -
  2 \beta_j(\om,z) \beta_l (\om,z) \sin(\beta_j (\om,z)\zeta) \Big]
\Big\}.
  \label{eq:kappaj}
\end{align}

Note that the coefficients of the partial derivatives with respect to
the mode powers $P_j$ are independent of the phases $\psi_j$. This
means that 
$\{|b_j^\ep(\om,z)|^2\}_{j=1}^{\mathcal{N}}$ converge in distribution
in the limit $\ep \to 0$ to the inhomogeneous diffusion Markov process
$\{P_j(\om,z)\}_{j=1}^{\mathcal{N}}$ with infinitesimal generator
$-\mathcal{L}_{P,z}^\mathcal{N}$ defined in \eqref{gendiff2P}. The
total power of the propagating modes satisfies
\begin{equation}
\mathcal{L}_{P,z}^\mathcal{N} \Big[\sum_{j=1}^\mathcal{N}
  P_j(\om,z)\Big] = \sum_{{\scriptsize \begin{array}{c}j,l = 1 \\ j
      \ne l \end{array}} }^{{\cal N}} G_{jl}^{(c)}(\om,z) \Big[
  P_l(\om,z) - P_j(\om,z) \Big] = 0,
\end{equation}
where we used \eqref{defgamma1b} and the symmetry of matrix
${\bf G}^{(c)}(\om,z)$. This implies that the total power is 
conserved
\begin{equation}
\sum_{j=1}^\mathcal{N} P_j(\om,z) = \mbox{constant}, \qquad z \in
\big(z_{-}^{(t)}, z_{-}^{(t-1)}\big).
\label{eq:EXACTCons}
\end{equation}

The evanescent waves do not contribute to the expression of the
infinitesimal generator $\mathcal{L}_{P,z}^\mathcal{N}$, so they do
not exchange energy with the propagating modes in the limit $\ep \to
0$. However, they appear in the last coefficient \eqref{eq:kappaj} of
the operator $\mathcal{L}_{\psi,z}^\mathcal{N}$, so they affect the
phases of the mode amplitudes.

The limit Markov process $\{b_j(\om,z)\}_{j=1}^{\mathcal N}$ is
inhomogeneous due to the slow variations of the waveguide which make
the coefficients of the operators \eqref{gendiff2P} and
\eqref{gendiff2T} $z$ dependent. The slow variations also
change the number of propagating modes at the turning points, and
this leads to partial reflection of power, as described in the next
section.

\section{Transport and reflection of power in the waveguide}
\label{sect:netscat}
We now use the infinitesimal generator \eqref{gendiffa} to
quantify the cumulative scattering effects in the waveguide.  We begin
in section \ref{sect:FWDTransp} with the modes transmitted through the
left part of the waveguide. The right going modes are discussed in
\ref{sect:REFTransp}. They are defined by the direct excitation from
the source and the reflection at the turning points. We end with some numerical illustrations in
section \ref{sect:Numerics}.

\subsection{The left going waves}
\label{sect:FWDTransp}
The wave propagation  from the source at $z=0$ to the end $z=-Z_M$
of the support of  variations of the waveguide can be described in
the limit $\ep \to 0$ as follows:

The left going mode amplitudes start  with the values
\begin{equation}
 b_{j}(\om,0-) = b_{j,o}(\om) = - \frac{f(\om) y_j(\rho_\star,0)}{2 i
   \sqrt{\beta_j(\om,0)}}, \qquad j = 1, \ldots, N^{(0)},
\label{eq:INIb0}
\end{equation}
obtained from equation \eqref{eq:CS14} and the observation that at $z
> 0$, where the opening $D(z)$ increases, the waves
are right going.

In the sector $\big(z_-^{(1)},0\big)$ the amplitudes $\{{b}_j(\omega,z)
\}_{j=1}^{N^{(0)}}$ evolve according to the diffusion Markovian
dynamics with generator $-{\cal L}_z^{N^{(0)}}$, starting from
$\{{b}_{j,o}(\om)\}_{j=1}^{N^{(0)}}$. The first $N^{(0)}-1$ left
going modes pass through the turning point
\begin{equation}
b_{j}(\om,z_-^{(1)}-) = b_{j}(\om,z_-^{(1)}+), \qquad j = 1, \ldots, N^{(0)}-1,
\label{eq:WT1}
\end{equation} 
but the last mode is reflected back. 

In the sector $\big(z_-^{(2)},z_-^{(1)}\big)$ there are $N_{-}^{(1)} =
N^{(0)}-1$ left going modes, with amplitudes evolving according to the
diffusion Markovian dynamics with generator $-{\cal
  L}_z^{N_{-}^{(1)}}$, starting from the values \eqref{eq:WT1} at $z =
z_-^{(1)}-$. At the next turning point $z_{-}^{(2)}$, only the first
$N_{-}^{(1)}-1$ modes pass through
\begin{equation}
b_{j}(\om,z_-^{(2)}-) = b_{j}(\om,z_-^{(2)}+), \qquad j = 1, \ldots, N_{-}^{(1)}-1,
\label{eq:WT2}
\end{equation} 
and the last mode is reflected back.

We continue this way until we reach $z = -Z_M$, with amplitudes
$\{b_j(\om,-Z_M)\}_{j=1}^{N_{\rm min}}$ obtained from the
diffusion Markovian dynamics with generator $-{\cal L}_z^{N_{\rm min}}$
over the interval $(-Z_M,z_{-}^{(t_M^{-})})$, starting with the values
$\{b_j(\om,z_{-}^{(t_M^{-})}-)\}_{j=1}^{N_{\rm min}}$ determined as
explained above, from the previous waveguide sectors.  

The waveguide has no variations at $z < Z_{-M}$, so the left going
mode amplitudes remain equal to their values at $-Z_M$, as stated in
equation \eqref{eq:CS17}. The emerging wave is obtained
from \eqref{eq:DC1} and \eqref{eq:FM10},
\begin{align}
p^\ep(\om,\rho,z) \approx -\sum_{j=1}^{N_{\rm min}} 
\frac{y_j(\rho,-Z_M)b_j(\om,-Z_M)}{\sqrt{\beta_j(\om,-Z_M)}} \exp \Big[-
  \frac{i}{\ep}\int_{0}^{-Z_M} dz' \beta_j(\om,z')\nonumber \\
  - \frac{i}{\ep}
  \beta_j(\om,-Z_M)(z+Z_M)\Big], \quad \mbox{for} ~ z < -Z_{M}.\label{eq:WT3}
\end{align}

\subsection{The mean transmitted wave field}
\label{eq:MeanTr}
With the infinitesimal generator \eqref{gendiffa} and Kolmogorov's
equation, we now calculate  the mean mode amplitudes 
\begin{equation}
\big<b_j(\om,z)\big> = \EE \big[b_j(\om,z)\big].
\label{eq:MEANNot}
\end{equation} 
In the first sector $(z_{-}^{(1)},0)$, these
satisfy the evolution equations
\begin{equation}
\partial_z \big<b_j(\om,z)\big> = -\Big[
G_{jj}^{(c)}(\om,z) - G_{jj}^{(0)}(\om,z) + i G_{jj}^{(s)}(\om,z) +2 i
\kappa_{j}^{N^{(0)}}(\om,z) \Big]\frac{\big<b_j(\om,z)\big>}{2}, 
\label{eq:WT4}
\end{equation}
solved backward in $z$, for $z \in \big(z_{-}^{(1)},0\big)$, starting
from the values 
\begin{equation}
\big<b_j(\om,0-)\big> = b_{j,o}(\om), \quad j = 1, \ldots, N^{(0)}.
\end{equation}
The coefficients in \eqref{eq:WT4} are defined by \eqref{defgamma1b},
\eqref{defgamma10}, \eqref{eq:kappaj} and 
\begin{equation}
G_{jj}^{(s)}(\om,z) = -\sum_{l = 1,l \ne j}^{N^{(0)}}
G_{lj}^{(s)}(\om,z),
\label{eq:WT6}
\end{equation}
with $G_{lj}^{(s)}(\om,z)$ given in \eqref{defgamma1s}.  
Because $
-G_{jj}^{(c)}(\om,z) + G_{jj}^{(0)}(\om,z) > 0$ (by Wiener-Khintchine theorem), 
we conclude from \eqref{eq:WT14} 
that the mean mode amplitudes decay with $|z|$, and therefore
\begin{equation}
\left| \big<b_j(\om,z_{-}^{(1)})\big> \right| < \big|
b_{j,o}(\om)\big|, \qquad 1 \le j \le N^{(0)}.
\label{eq:WT8}
\end{equation}
This decay models the randomization of the left going modes, and
occurs on a $j$ dependent length scale, as illustrated in section
\ref{sect:Numerics}. Similar to the case of  waveguides with random
perturbations of straight boundaries \cite[Section 5]{alonso2011wave},
the modes with larger index $j$ randomize faster. Intuitively, this is
because these modes propagate slowly along $z$, at group velocity
$1/\partial_\om \beta_j(\om,z)$ that is small with respect to the wave
speed, and bounce more often at the random boundary.

A similar calculation applies to the other sectors
$\big(z_{-}^{(t)},z_{-}^{(t-1)}\big)$ of the waveguide, indexed by $t
= 1, \ldots t_M^{-}$. The only difference is that the starting values
of the mode amplitudes are random, so we use conditional expectations
\begin{equation}
\big< b_j(\om,z)\big> = \EE \left[ \EE \left[ b_j(\om,z) \big|
    {\cal F}_{z_{-}^{(t-1)}} \right] \right], \qquad z < z_{-}^{(t-1)},
\label{eq:CONDEXP}
\end{equation}
where $\mathcal{F}_{z_{-}^{(t-1)}}$ denotes the $\sigma$-algebra (
information) generated by the Markov limit process
$\left\{b_q(\om,z)\right\}_{q=1}^{N_{-}^{(t-1)}}$ at $z =
z_{-}^{(t-1)}$.  We obtain that $\big<b_j(\om,z)\big>$ satisfies an
equation like \eqref{eq:WT4}, with redefined coefficients for the
$N_{-}^{(t-1)}$ number of propagating modes, and starting value 
$\big<b_j(\om,z_{-}^{(t-1)})\big>$ calculated in the previous
waveguide sector.

Proceeding this way we reach $z = -Z_{M}$. The mean transmitted wave is the
expectation of \eqref{eq:WT3}, with $\big<b_j(\om,-Z_M)\big>$
obtained by solving equations \eqref{eq:WT4} for all the sectors of
the waveguide.  The scattering effects at the random boundary add up
in each sector, and the mean mode amplitudes decay, as explained
above,
\begin{equation}
\left| \big<b_j(\om,-Z_{M})\big> \right| < \left|
\big<b_j(\om,-z_{-}^{(t_M^{-})}\big> \right| < \ldots < \big|
b_{j,o}(\om)\big|, \qquad 1 \le j \le N_{\rm min}.
\end{equation}

\subsection{The transmitted power}
\label{eq:MeanPow}
Using the infinitesimal generator \eqref{gendiff2P} of the Markov
process $\{P_j(\om,z)\}$, the $\ep \to 0$ limit of the left going mode
powers, we now calculate the mean and standard deviation of the
transmitted power at $z < 0$.

We proceed as in the previous section, one sector of the waveguide at
a time, starting from the source. In the first sector $z \in
\big(z_{-}^{(1)},0\big)$, the mean powers
\begin{equation}
\big< P_j(\om,z) \big> = \EE \left[ P_j(\om,z) \right], \qquad j =
1, \ldots, N^{(0)},
\label{eq:WT10}
\end{equation}
evolve from the initial values $\big<P_j(\om,0-) \big> =
|b_{j,o}(\om)|^2$ according to equation
\begin{equation}
\partial_z \begin{pmatrix} \big< P_1(\om,z) \big> \\ \vdots
  \\\big< P_{N^{(0)}}(\om,z) \big> \end{pmatrix} = 
-{\bf G}^{(c)}(\om,z) \begin{pmatrix} \big< P_1(\om,z) \big> \\ \vdots
  \\\big< P_{N^{(0)}}(\om,z) \big> \end{pmatrix}, 
\label{eq:WT11}
\end{equation}
with matrix ${\bf G}^{(c)}(\om,z)$ defined in
\eqref{defgamma1}-\eqref{defgamma1b}, for $\mathcal{N} = N^{(0)}$.

In the next sectors $(z_{-}^{(t)},z_{-}^{(t-1)})$  we use conditional
expectations
\begin{equation}
\big< P_j(\om,z) \big> = \EE \left[ \EE \left[P_j(\om,z)\big|
    {\cal F}_{z_{-}^{(t-1)}} \right]\right], \qquad z < z_{-}^{(t-1)},
\label{eq:WT12}
\end{equation}
and obtain that the mean powers satisfy an equation like
\eqref{eq:WT11}, with $N_{-}^{(t-1)}$ unknowns and $N_{-}^{(t-1)}
\times N_{-}^{(t-1)}$ matrix ${\bf G}^{(c)}(\om,z)$. These equations
are solved backward in $z$, starting from the values $\big<
P_j(\om,z_{-}^{(t-1)}) \big>$ computed in the previous sectors. Proceeding this way, we reach 
$z = -Z_M$, and
obtain $\big< P_j(\om,-Z_M) \big>$, for $j = 1, \ldots, N_{\rm min}.$

Note that unlike the expectations \eqref{eq:MEANNot}, the mean powers
are coupled by the matrix ${\bf G}^{(c)}(\om,z)$. This coupling models
the exchange of power between the left going modes, induced by
cumulative scattering at the random boundary of the waveguide. The
exchange depends on the mode index, as illustrated in section
\ref{sect:Numerics}. Specifically,  the higher indexed
modes transfer  power more quickly than the others.

How much power is exchanged depends on the length of the sectors
$\big(z_{-}^{(t)},z_{-}^{(t-1)}\big)$ of the waveguide. In short
sectors, the exchange is mostly among the higher indexed modes. The
longer the sectors, the more modes participate in the exchange and the
power may become evenly distributed among the modes, independent of
the starting value at $z_{-}^{(t-1)}$. This equipartition of energy
has been explained in waveguides with straight walls in \cite[Section
  20.3]{fouque07}, for a matrix ${\bf G}^{(c)}$ with non-zero off
diagonal entries. By the Perron-Frobenius theorem, and due to energy conservation, such a matrix has a
simple eigenvalue equal to  zero, and the
other eigenvalues are negative.  It is straightforward to see from
equation \eqref{eq:WT11} that the solution converges at large $|z|$ to
a vector in the null space of ${\bf G}^{(c)}$. Equation
\eqref{defgamma1b} gives that this space is spanned by the vector of
all ones, so the power becomes evenly distributed at distances that
exceed the equipartition distance. This length scale is defined by the
inverse of the absolute value of the largest, non-zero eigenvalue of
${\bf G}^{(c)}$.

By the energy conservation \eqref{eq:EXACTCons}, the transmitted power in the first sector of
the waveguide is 
\begin{equation}
\mathcal{P}_{\rm trans}(\om,z) = \sum_{j=1}^{N^{(0)}} P_j(\om,z) =
\sum_{j=1}^{N^{(0)}} |b_{j,o}(\om)|^2, \qquad z \in (z_{-}^{(1)},0),
\label{eq:WT13}
\end{equation}
where the right hand side is the deterministic, total left going power emitted by the
source. At the turning point $z_{-}^{(1)}$ the $N^{(0)}$-th mode is
reflected back. The transmitted power to the next sector of the
waveguide, carried by the remaining $N_{-}^{(1)} = N^{(0)} - 1$  modes,
is random and given by 
\begin{equation}
\mathcal{P}_{\rm trans}(\om,z) = \sum_{j=1}^{N_{-}^{(1)}} P_j(\om,z) =
\sum_{j=1}^{N_{-}^{(1)}} P_j(\om,z_{-}^{(1)}), \qquad z \in
(z_{-}^{(2)},z_{-}^{(1)}).
\label{eq:WT14}
\end{equation}
This repeats for the other sectors, and beyond $z=-Z_M$ we have 
\begin{equation}
\mathcal{P}_{\rm trans}(\om,z) = \sum_{j=1}^{N_{\rm min}} P_j(\om,-Z_M),
\qquad z \le -Z_M.
\label{eq:WT15}
\end{equation}

In summary, the transmitted power is a piecewise constant
function with jumps at the turning points, and random values
determined by the sum of the mode powers entering each sector of the
waveguide. Its mean is obtained by taking expectations in
\eqref{eq:WT13}-\eqref{eq:WT15}, and using the mean mode powers
calculated as explained above.
 
The random fluctuations of $\mathcal{P}_{\rm trans}(\om,z)$
about the mean are quantified by its standard deviation
\begin{equation}
\mbox{StD} \left[ \mathcal{P}_{\rm trans}(\om,z)\right] = \Big\{
\sum_{j,l=1}^{N_{-}^{(t-1)}} \left[ \big< \mathcal{P}_{jl}(\om,z)
  \big> - \big< P_{j}(\om,z)\big> \big< P_{l}(\om,z) \big>
  \right]\Big\}^{1/2}
\label{eq:WT20}
\end{equation}
for $z \in (z_{-}^{(t)}, z_{-}^{(t-1)})$ and $1\le t \le t_M^{-}$.
 To calculate it we need the
second moments
\begin{equation}
\big< \mathcal{P}_{jl}(\om,z)\big> = \EE \left[ {P}_{j}(\om,z)
  {P}_{l}(\om,z)\right].
\label{eq:WT16}
\end{equation}
Again, these  are obtained in one sector of the waveguide at a
time, starting from the source, where
\begin{equation}
\big< \mathcal{P}_{jl}(\om,0)\big> = |b_{j,o}(\om)|^2
|b_{l,o}(\om)|^2, \qquad j,l = 1, \ldots, N^{(0)}.
\label{eq:WT17}
\end{equation}
The evolution equations of the moments \eqref{eq:WT16} at $z \in
\big(z_{-}^{(t)},z_{-}^{(t-1)}\big)$ are
\begin{equation}
\partial_z \big< \mathcal{P}_{jj}(\om,z)\big> = 2
G_{jj}^{(c)}(\om,z) \big< \mathcal{P}_{jj}(\om,z)\big> - 4
\sum_{l=1}^{N_{-}^{(t-1)}} G_{jl}^{(c)}(\om,z) \big<
\mathcal{P}_{lj}(\om,z)\big>, 
\label{eq:WT18}
\end{equation}
and
\begin{align}
\partial_z \big< \mathcal{P}_{jq}(\om,z)\big> = 2
G_{jq}^{(c)}(\om,z) \big< \mathcal{P}_{jq}(\om,z)\big> - 
\sum_{l=1}^{N_{-}^{(t-1)}} \Big[G_{jl}^{(c)}(\om,z) \big<
  \mathcal{P}_{lq}(\om,z)\big> \nonumber \\ + G_{lq}^{(c)}(\om,z)
  \big< \mathcal{P}_{jl}(\om,z)\big> \Big],
\label{eq:WT19}
\end{align}
for $j, q = 1, \ldots, N_{-}^{(t-1)}$ and $j \ne q$. These equations
are solved backward in $z$, with the starting values $\big<
\mathcal{P}_{jq}(\om,z_{-}^{(t-1)})\big>$ calculated from the
previous sector.  
\subsection{The right going waves}
\label{sect:REFTransp}
Even though we consider the forward scattering approximation in each
sector of the waveguide, there are both left and right going modes at
$z < 0$, due to reflection at the turning points. At $z > 0$ we also
have the right going waves emitted from the source.  The analysis of
the reflected mode amplitudes is more complicated, because they
quantify cumulative scattering in the waveguide sectors traversed both
ways: to the left by the incoming wave and to the right by the
reflected wave.

In each sector $\big(z_{-}^{(t)},z_{-}^{(t-1)}\big)$ we obtain from
\eqref{eq:fwdprop} that the right going mode amplitudes satisfy
\begin{equation}
{\itbf a}^\ep(\om,z) \approx
\overline{\calP^{\ep(bb)}}(\om,z;z_{-}^{(t-1)}) {\itbf
  a}^\ep(\om,z_{-}^{(t-1)}), \qquad t = 1, \ldots, t_M^-.
\label{eq:REFL1}
\end{equation}
This looks similar to equation \eqref{evola} that describes the
evolution of the left going waves, but we have different boundary
conditions, as we now explain. 

Starting from the leftmost turning point
$z_{-}^{(t_M^{-})}$, and denoting $\mathcal{N} = N_{-}^{(t_M^{-}-1)}$,
we obtain from \eqref{eq:CS16} the initial condition
\begin{equation}
a_j^{\ep}\big(\om,z_{-}^{(t_M^{-})}\big) = {\rm R}_{\mathcal{N}}^\ep
\big(\om,z_{-}^{(t_M^{-})}\big) b_{\mathcal
  N}^\ep\big(\om,z_{-}^{(t_M^{-})}\big) \delta_{j \mathcal{N}}, \qquad
j = 1, \ldots, \mathcal{N},
\label{eq:REFL2}
\end{equation}
for the vector ${\itbf a}^\ep(\om,z) \in \mathbb{C}^{\mathcal{N}}$, where
$\delta_{j \mathcal{N}}$ is the Kronecker delta symbol and
${\rm R}_{\mathcal{N}}^\ep$ is reflection coefficient defined in
\eqref{eq:CS11}. The amplitudes of the right going modes
impinging on the next turning point are obtained from \eqref{eq:REFL1}
\begin{equation}
{\itbf a}^\ep(\om,z_{-}^{(t_M^{-}-1)}-) \approx \Big[\calP^{\ep(bb)}(\om,z_{-}^{(t_M^{-})};z_{-}^{(t_M^{-}-1)})\Big]^T {\itbf
    a}^\ep(\om,z_{-}^{(t_M^{-})}),
\label{eq:REFL3}
\end{equation}
using that the propagator $\calP^{\ep(bb)}$ is approximately
unitary. This follows from the energy conservation relation
\eqref{eq:EXACTCons}, which holds in the limit $\ep \to 0$,
independent of the initial conditions.

On the right of the turning point $z_{-}^{(t_M^{-}-1)}$ there is an
extra right going mode. Renaming $\mathcal{N} = N_{-}^{(t_M^--2)}$, we
obtain the following initial condition for the vector ${\itbf a}^\ep(\om,z)$: Its first $\mathcal{N}-1$ components are
given in \eqref{eq:REFL3}, and the last component is
\begin{equation} 
a_{\mathcal{N}}^\ep\big(\om,z_{-}^{(t_M^{-}-1)}\big) =
{\rm R}_{\mathcal{N}}^\ep \big(\om,z_{-}^{(t_M^{-}-1)}\big) b_{\mathcal N}^\ep\big(\om,z_{-}^{(t_M^{-}-1)}\big).
\label{eq:REFL4}
\end{equation}
These amplitudes and the $\mathcal{N} \times \mathcal{N}$ propagator
$\calP^{\ep(bb)}(\om,z;z_{-}^{(t_M^{-}-2)})$ determine the
amplitudes of the right going modes impinging on the turning point
$z_{-}^{(t_M^--2)}$ and so on.

Proceeding this way we obtain the amplitudes
$\{a_j^\ep(\om,0-)\}_{j=1}^{N^{(0)}}$ on the left of the source.
The amplitudes at $z = 0+$ are given by these and the source conditions
\eqref{eq:CS14}. The analysis of forward propagation at $z>0$ is
similar to that in section \ref{sect:FWDTransp}, with the exception
that at the turning points $z_{+}^{(t)}$, for $1\le t \le 
t_{M}^+$, there is no reflection. We add instead a new mode with zero
initial condition, as stated in \eqref{eq:CS16P}.

\subsection{The net reflected power}
The calculation of the statistical moments of the right going mode
amplitudes in the limit $\ep \to 0$ requires the infinitesimal
generator of the limit propagator $\calP^{\ep(bb)}$, in
each sector of the waveguide. This operator can be obtained using
Theorem \ref{prop.diflim}, but the calculation is complex.  Here
we quantify only the net reflected power at each turning point,
without asking how this power gets distributed among the modes as they
propagate toward the right.

The net reflected power is determined by the transmitted power in the
left part of the waveguide, using energy conservation. Specifically,
starting from the leftmost turning point, the net reflected power is
\begin{equation}
\mathcal{P}_{\rm refl}(\omega,z) = P_{N_{-}^{(t_M^- -1)}}(\om,z_{-}^{(t_M^-)}+),
\qquad z \in \big(z_{-}^{(t_M^-)},z_{-}^{(t_M^- -1)}\big),
\label{eq:REFL5}
\end{equation}
where the right hand side is the power of the left going turning mode,
analyzed in section \ref{sect:FWDTransp}. Here we used the conservation relation
\[
\lim_{\ep \to 0} \sum_{j=1}^{N_{-}^{(t_M^- -1)}} |a_j^\ep(\om,z)|^2 = \mbox{constant}, 
\quad \mbox{for}~ z  \in  \big(z_{-}^{(t_M^-)},z_{-}^{(t_M^- -1)}\big),
\]
derived the same way as  \eqref{eq:EXACTCons}, equation \eqref{eq:REFL2} and 
$\lim_{\ep \to 0}|{\rm R}_{\mathcal{N}}^\ep| = 1$.

At the next turning point $z_{-}^{(t_M^- -1)}$ we add a new mode amplitude, and the net reflected power increases to
\begin{equation}
\mathcal{P}_{\rm refl}(\omega,z) = P_{N_{-}^{(t_M^- - 1)}}(\om,z_{-}^{(t_M^-)}+) +
P_{N_{-}^{(t_M^- -2)}}(\om,z_{-}^{(t_M^- - 1)}+), 
\label{eq:REFL6}
\end{equation}
for $z \in \big(z_{-}^{(t_M^- - 1)},z_{-}^{(t_M^- -2)}\big)$, and so
on.  Proceeding this way we obtain that the net reflected power is a
piecewise constant function at $z < 0$, with jumps at the turning
points $z_{-}^{(t)}$ indexed by  $1 \le t \le t_{M}^-$. At the source location
this equals
\begin{equation}
\mathcal{P}_{\rm refl}(\omega,0) = \sum_{t=1}^{t_{M}^-}
P_{N_{-}^{(t-1)}}(\om,z_{-}^{(t)}+),
\label{eq:REFL7}
\end{equation}
and its mean and standard deviation are determined by those of the
turning wave powers, calculated in section \ref{sect:FWDTransp}.
By comparing with (\ref{eq:WT13}-\ref{eq:WT15}) we obtain the global conservation of energy relation
\begin{equation}
\mathcal{P}_{\rm refl}(\omega,0) + \mathcal{P}_{\rm trans}(\omega,-Z_M)
=
\sum_{j=1}^{N^{(0)}} |b_{j,o}(\om)|^2 .
\end{equation}
Therefore the first two moments of the net transmitted and reflected powers are related through:
\begin{eqnarray}
\left< \mathcal{P}_{\rm refl}(\om,0)\right> 
&=&
\sum_{j=1}^{N^{(0)}} |b_{j,o}(\om)|^2 - \left< \mathcal{P}_{\rm trans}(\om,-Z_M)\right>,
\\
\mbox{StD} \left[ \mathcal{P}_{\rm refl}(\om,0)\right] 
&=&
\mbox{StD} \left[ \mathcal{P}_{\rm trans}(\om,-Z_M)\right] .
\end{eqnarray}

\subsection{The net power transmitted to the right}
There is no mode reflection at $z>0$, and the net transmitted power to the right is
\begin{equation}
\mathcal{P}_{\rm trans,right}(\om,z) = \lim_{\ep \to 0} \sum_{j=1}^{N^{(0)}}
|a_j^\ep(\om,0+)|^2, \qquad z > 0,
\label{eq:REFL8}
\end{equation}
where the equality means having the same statistical distribution, and 
\begin{equation}
a_j^\ep(\om,0+) = a_j^\ep(\om,0-) + a_{j,o}(\om), \qquad a_{j,o}(\om)
= \frac{\hat f(\om)y_j(\rho_\star,0)}{2 i \sqrt{\beta_j(\om,0)}}.
\end{equation}
The calculation of the statistical moments of \eqref{eq:REFL8} is
 as complicated as the calculation of the moments
of the limit right going mode amplitudes. Specifically, it requires the infinitesimal
generator of the $\ep \to 0$ limit of the propagator $\calP^{\ep(bb)}$,
in particular, we need to characterize the phases of the reflection coefficients 
${\rm R}_{N_-^{(t-1)}}^\ep \big(\om,z_{-}^{(t)}\big)$,
$t=1,\ldots,t_M^{-}$.
By extrapolating the results given in \cite{borcea2016turning} (in which the standard deviation of the fluctuations of the boundary was smaller),
we could anticipate that these phases are independent and uniformly distributed over $[0,2\pi]$.
We could then anticipate that the mean power transmitted to the right is
\begin{eqnarray}
\nonumber
\left< \mathcal{P}_{\rm trans,right}(\om,z) \right> &=& \left< \mathcal{P}_{\rm refl}(\om,0) \right>
+\sum_{j=1}^{N^{(0)}}
| a_{j,o}(\om)|^2
\\ &=&
\sum_{j=1}^{N^{(0)}}
| a_{j,o}(\om)|^2 +
\sum_{j=1}^{N^{(0)}} |b_{j,o}(\om)|^2 - \left< \mathcal{P}_{\rm trans}(\om,-Z_M)\right>
 ,
\end{eqnarray}
for any $z>0$.

\subsection{Numerical illustration}
\label{sect:Numerics}
In this section we illustrate with some plots the exchange of power
among the propagating modes in the left part $z<0$ of the waveguide,
due to a point source at $\bx_\star = (D(0)/7,0)$. For comparison, we
also consider other initial conditions, where the excitation at $z =
0$ is for a single mode at a time.

We take a waveguide with a straight axis that has a single turning
point, at arc length $z_{-}^{(1)} = -L = -1000 \la$, where $\la$ is the wavelength. The waveguide
opening $D(z/L)$ increases linearly in $z$ 
in the interval $[-L, 0]$, from the value $20 \la$ to $20.49\la$, and 
transitions as a cubic polynomial to the constant $19.999\la$ at $z < -L - 0.2 \la$ and 
$20.491\la$ at $z > 0.2 \la$. 
Thus, there are $N^{(0)} = 40$ propagating modes at $z > -L$ and $N_{-}^{(1)} = 39$ modes at $z < -L$. The  top and bottom  boundaries of
the waveguide are straight and parallel at $z\in(-\infty, -L-0.2\la)\cup (0.2 \la,\infty)$.  

The auto-correlation function $\cR$ of the process $\nu(\zeta)$ is a
Gaussian with standard deviation $1$. The correlation length of the fluctuations  is
$\ell = 3 \la$, so  $\ep = \ell/L = 0.003$, and
the standard deviation $\sigma$ of the fluctuations equals
$\sqrt{\ep}$.

We can describe approximately what to expect in terms of the
randomization of the mode amplitudes and the exchange of power among
the modes by looking at the following length scales calculated in a
waveguide with constant opening equal to $D(0)$: 

\vspace{0.06in} 1.The mode dependent scattering mean free path
\begin{equation}
L_{j,{\rm smf}} = \frac{2}{G_{jj}^{(0)}(\om,0)-G_{jj}^{(c)}(\om,0)}, \qquad j = 1, \ldots, 40,
\label{eq:Num1}
\end{equation}
which is the scale of decay of the mean
mode amplitudes, as seen from \eqref{eq:WT4}.

\vspace{0.06in}  2. The  mode dependent transport mean free paths,
\begin{equation}
L_{j,{\rm tmf}} = -\frac{2}{G_{jj}^{(c)}(\om,0)}, \qquad j = 1, \ldots, 40,
\label{eq:Num2}
\end{equation}
defined in terms of the diffusion coefficient $-G_{jj}^{(c)}$ of the mode power 
infinitesimal generator \eqref{gendiff2P}. The modes exchange power with their neighbors as 
they propagate at 
distances of  order \eqref{eq:Num2}.  

\vspace{0.06in}
3. The
equipartition distance $L_{\rm eq}$, which is defined as the inverse of
the absolute value of the largest, non-zero eigenvalue of matrix
${\bf G}^{(c)}(\om,0)$. At distances of order $L_{\rm eq}$, we expect that
the power gets evenly distributed among the modes, independent of the
excitation at $z = 0$.

\begin{figure}[t]
\begin{center}
\includegraphics[width=0.5\textwidth]{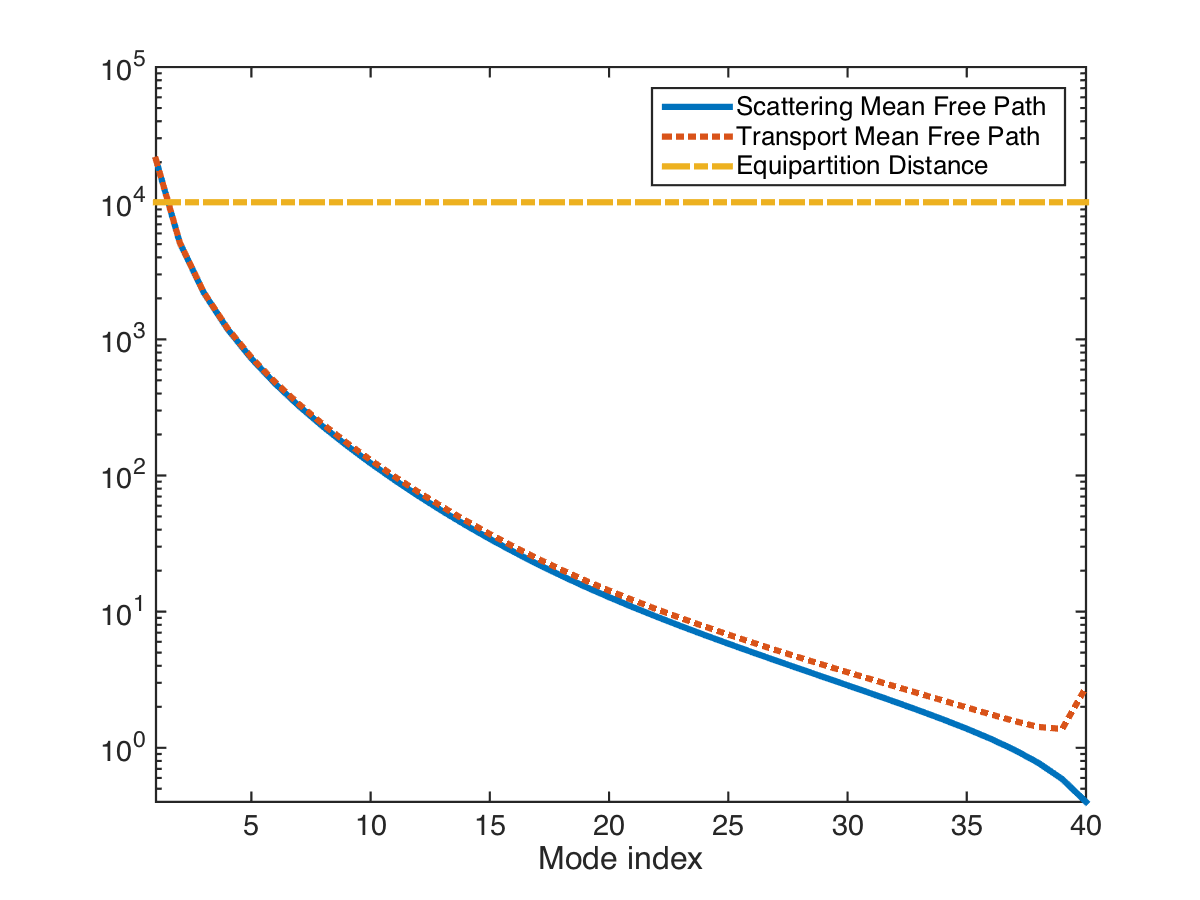}
\end{center}
\vspace{-0.1in}
\caption{The three length scales that quantify net scattering in a
  waveguide with constant opening $D(0)$. The solid blue line is for
  the scattering mean free path \eqref{eq:Num1}. The dashed red line
  is for the transport mean free path \eqref{eq:Num2}. The yellow
  dashed line is for the equipartition distance. The abscissa is the
  mode index $j = 1, \ldots, 40$ and the ordinate is in units of
  $\la$. }
\label{fig:scales}
\end{figure}

\vspace{0.06in} We display these scales in Figure
\ref{fig:scales}
and observe that at the distance $L = 1000\la$ between the source and
the turning point, we have 
\[
L \ge L_{j,{\rm smf}}, L_{j,{\rm tmf}}, \qquad j = 5, \ldots, 40.
\]
Thus, these modes should be randomized and
moreover, they should share their power with the other modes.  Because
$L < L_{\rm eq}$, we expect that at least the first five modes have not
shared all their power with the other modes.

\begin{figure}[t]
\begin{center}
\includegraphics[width=0.5\textwidth]{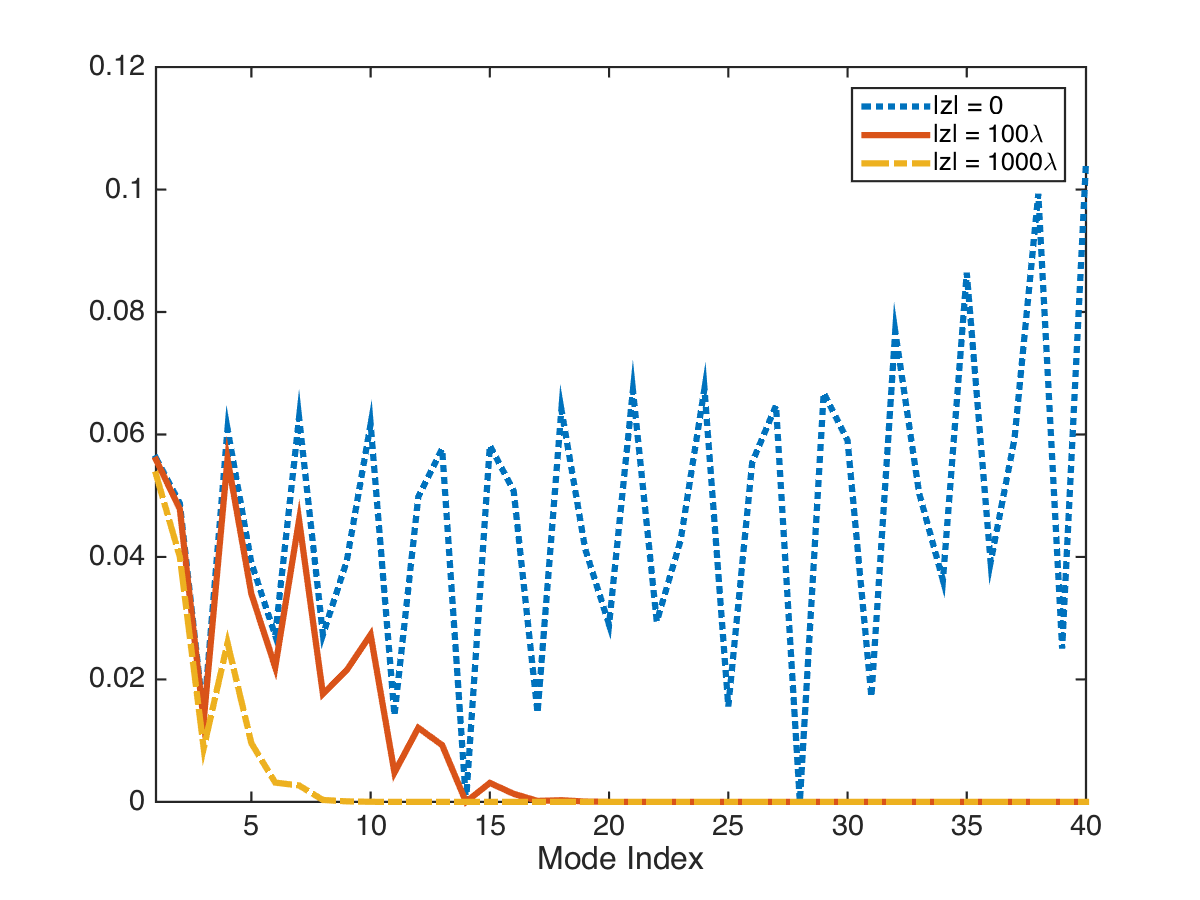}
\hspace{-0.1in}\includegraphics[width=0.5\textwidth]{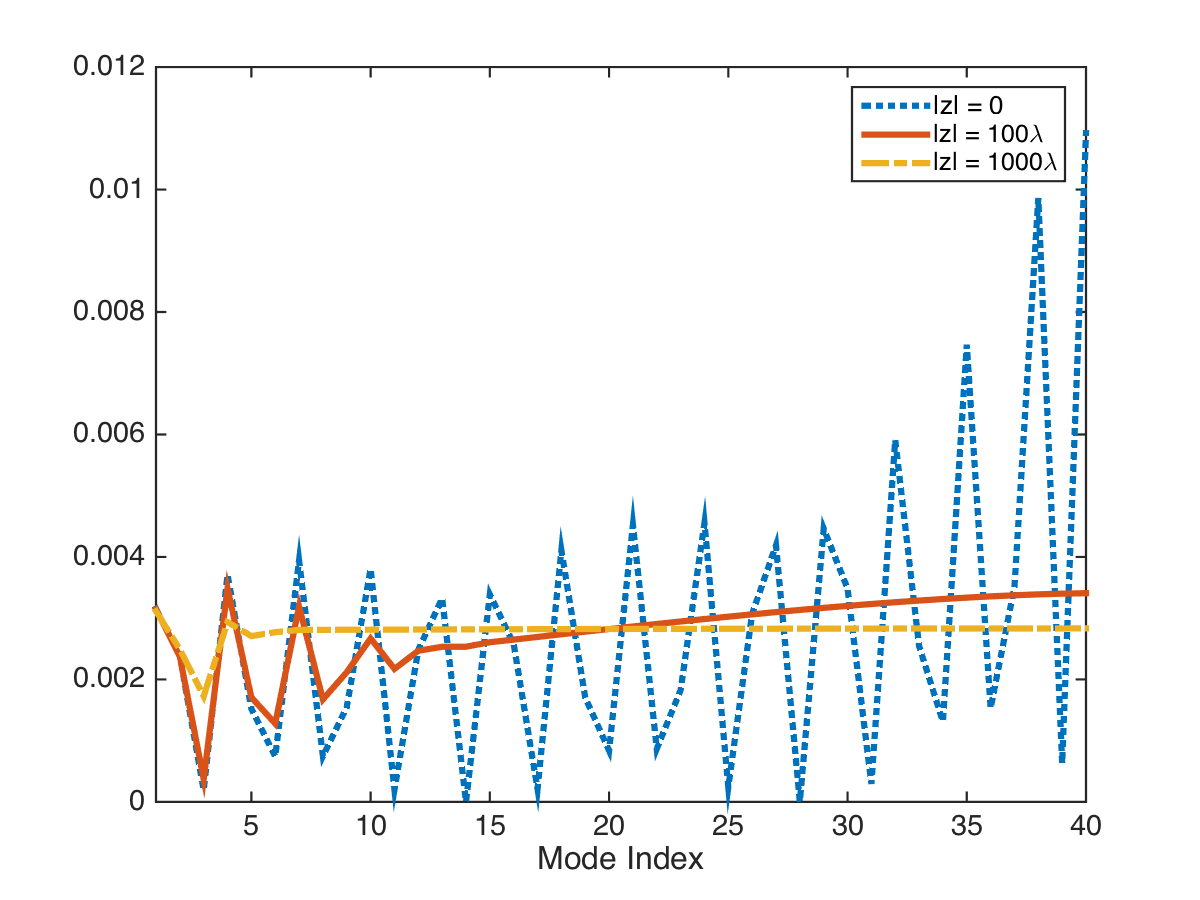}
\end{center}
\vspace{-0.1in}
\caption{Display of the absolute value of the mean mode amplitudes
  $|\big<b_j(\om,z)\big>|$ (left) and the mean mode powers
  $\big<P_j(\om,z)\big>$ v.s. the mode index $j$ at three different
  distances from the source: The blue dashed line corresponds to the
  initial values at $z=0$, due to a point source at location
  $(D(0)/7,0)$.  The full red line is for $|z| = 100\la$ and the
  yellow line is for $|z| = L = 1000\la$.  The abscissa is the mode
  index $j = 1, \ldots, 40$. }
\label{fig:Means}
\end{figure}

These expectations are confirmed by the results displayed in Figure
\ref{fig:Means}, where we show the absolute values
$|\big<b_j(\om,z)\big>|$ of the mean mode amplitudes (left plot)
and the mean mode powers $\big<P_j(\om,z)\big>$ (right plot) 
at three distances from the point source.
The dashed blue line is for $z = 0$, so it corresponds to the initial
values \eqref{eq:INIb0} of the mode amplitudes, which oscillate in $j$
due to the factor
\[
y_j(\rho_\star,0) = \sqrt{\frac{2}{D(0)}} \sin \Big[
\Big(\frac{\rho_\star}{D(0)} + \frac{1}{2} \Big) \pi j \Big], \qquad j = 1, \ldots, N^{(0)}, ~ ~ \rho_\star = \frac{D(0)}{7}.
\]
As we increase the
distance $|z|$ from the source, the left plot in Figure
\ref{fig:Means} illustrates the decay of the mean mode amplitudes. We note
that at $|z| = 100 \la$, the modes indexed by $j > 15$ have negligible
mean, and at the turning point $|z| = L = 1000\la$, the modes indexed
by $j > 5$ have negligible mean. This is as expected from Figure
\ref{fig:scales}, because because $L_{j,{\rm smf}} < 100 \la$ for $j > 15$
and $L_{j,{\rm smf}} < 1000 \la$ for $j > 5$.  The right plot in Figure
\ref{fig:Means} illustrates the effect of exchange of power among the
modes. The scattering mean free path and the transport mean free path
are almost the same in this simulation, as shown in Figure
\ref{fig:scales}, and we note that at the turning point $|z| = L =
1000 \la$ the modes indexed by $j > 5$ have almost the same power.

\begin{figure}[t]
\begin{center}
\includegraphics[width=0.5\textwidth]{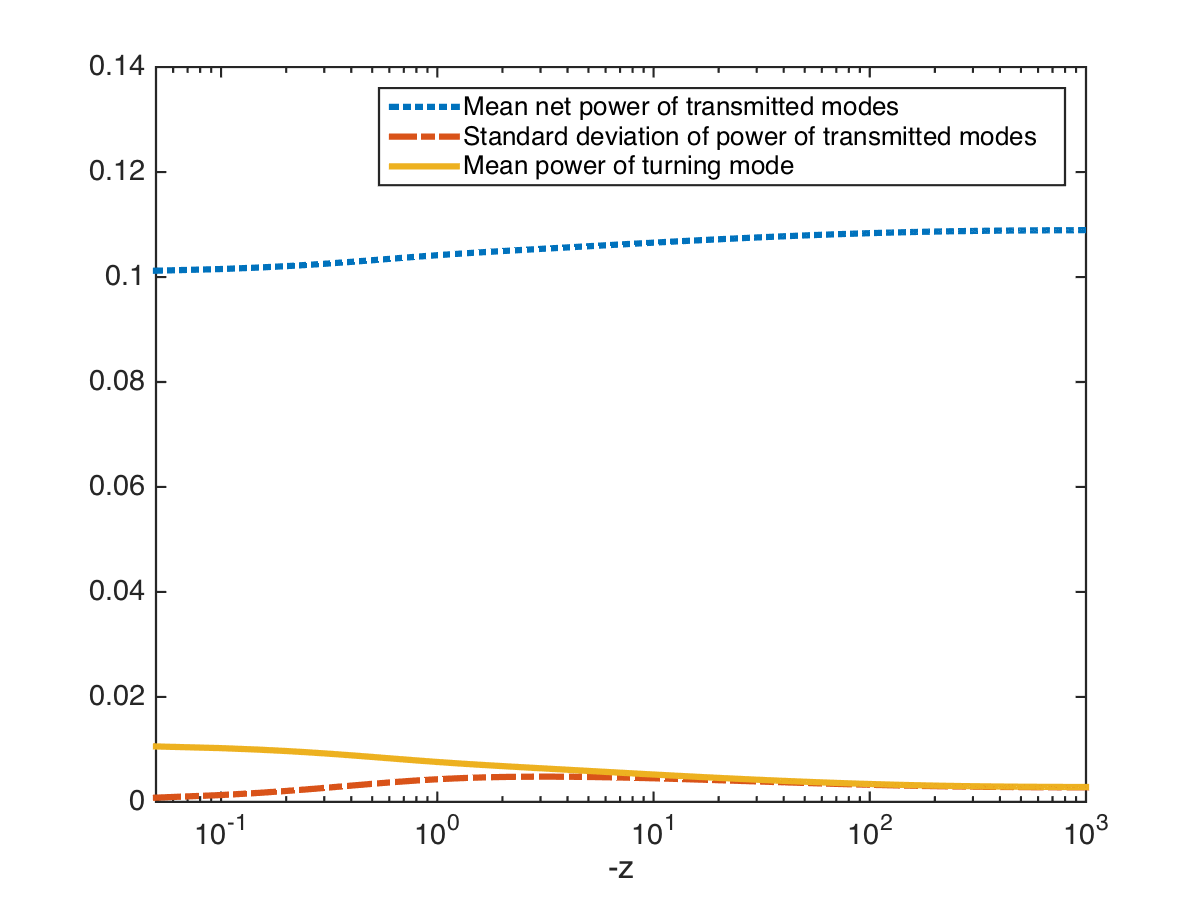}
\end{center}
\vspace{-0.1in}
\caption{Display of the mean net power of the transmitted modes in
  dashed blue line, of the standard deviation of this power in dashed
  red line, and the mean power of the turning mode, indexed by $j =
  40$.  The abscissa is the arc length in units of $\la$ (in
  logarithmic scale). }
\label{fig:TransPT}
\end{figure}

In Figure \ref{fig:TransPT} we display the mean and standard deviation
of the net power $\sum_{j=1}^{39} P_j(\om,z)$ of the modes that are
transmitted through the turning point, and the mean power of the
turning mode, as functions of $z$. At $|z| = L = 1000\la$, these
determine the transmitter power \eqref{eq:WT14} beyond the turning
point, and the reflected power \eqref{eq:REFL7}. Note that in this
case cumulative scattering at the random boundary is beneficial for
power transmission through the waveguide. In the absence of the random
fluctuations  there would be no power
exchange between the modes, and the transmitted power would equal
$\sum_{j=1}^{39} P_j(\om,0)$. As seen in Figure \ref{fig:Means}, the
turning mode has the largest mode amplitude initially, and all its
power would be reflected back. The cumulative scattering at the random
boundary leads to rapid exchange of the power of the turning mode, as
shown in the right plot of Figure \ref{fig:Means}, and much less power
is reflected.  The standard deviation of the net power of the first
$39$ modes, shown with the red dashed line in Figure \ref{fig:TransPT},
is smaller than its mean. Thus,
$
\sum_{j=1}^{39} P_j(\om,z) \approx \sum_{j=1}^{39}
\big<P_j(\om,z)\big>,
$
with less than $10\%$  relative error (i.e., random fluctuations).

\begin{figure}[t]
\begin{center}
\includegraphics[width=0.5\textwidth]{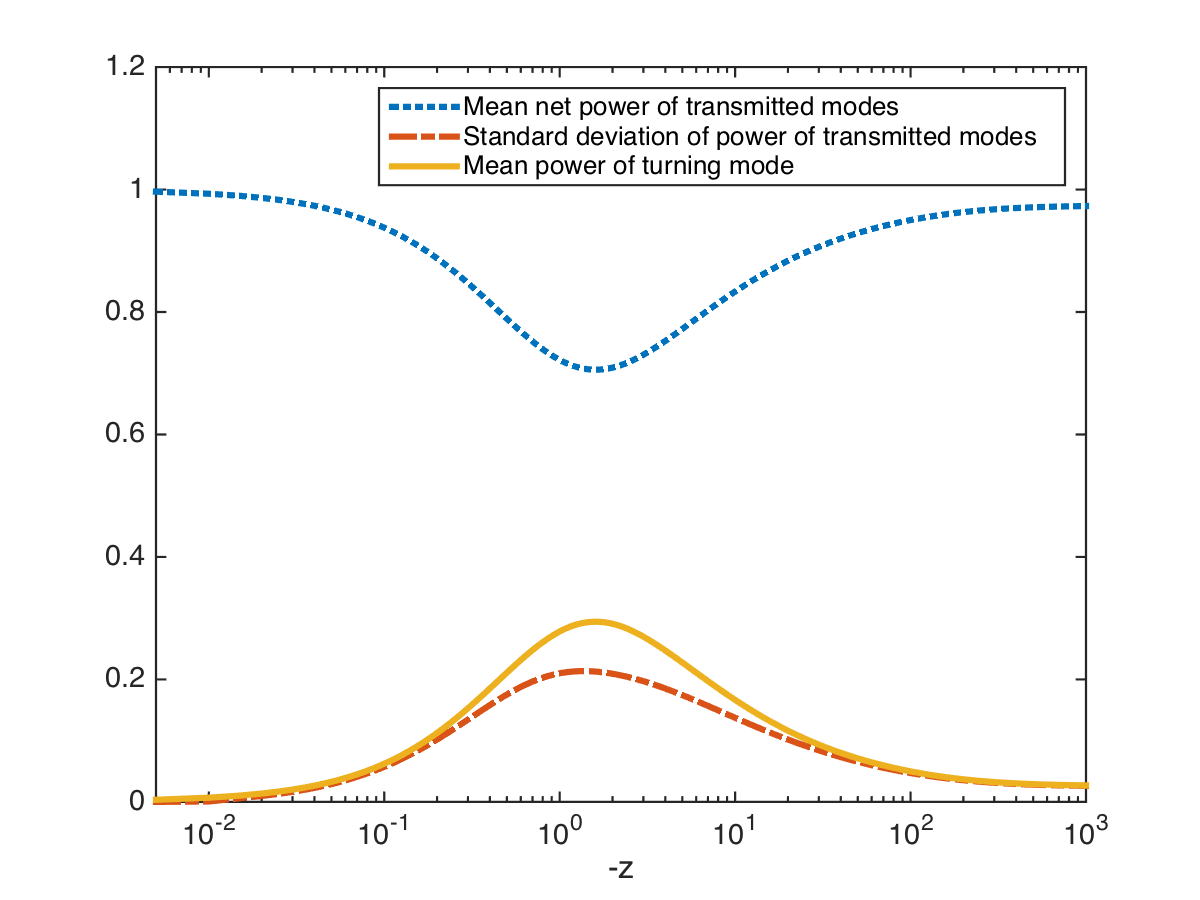}
\hspace{-0.1in}\includegraphics[width=0.5\textwidth]{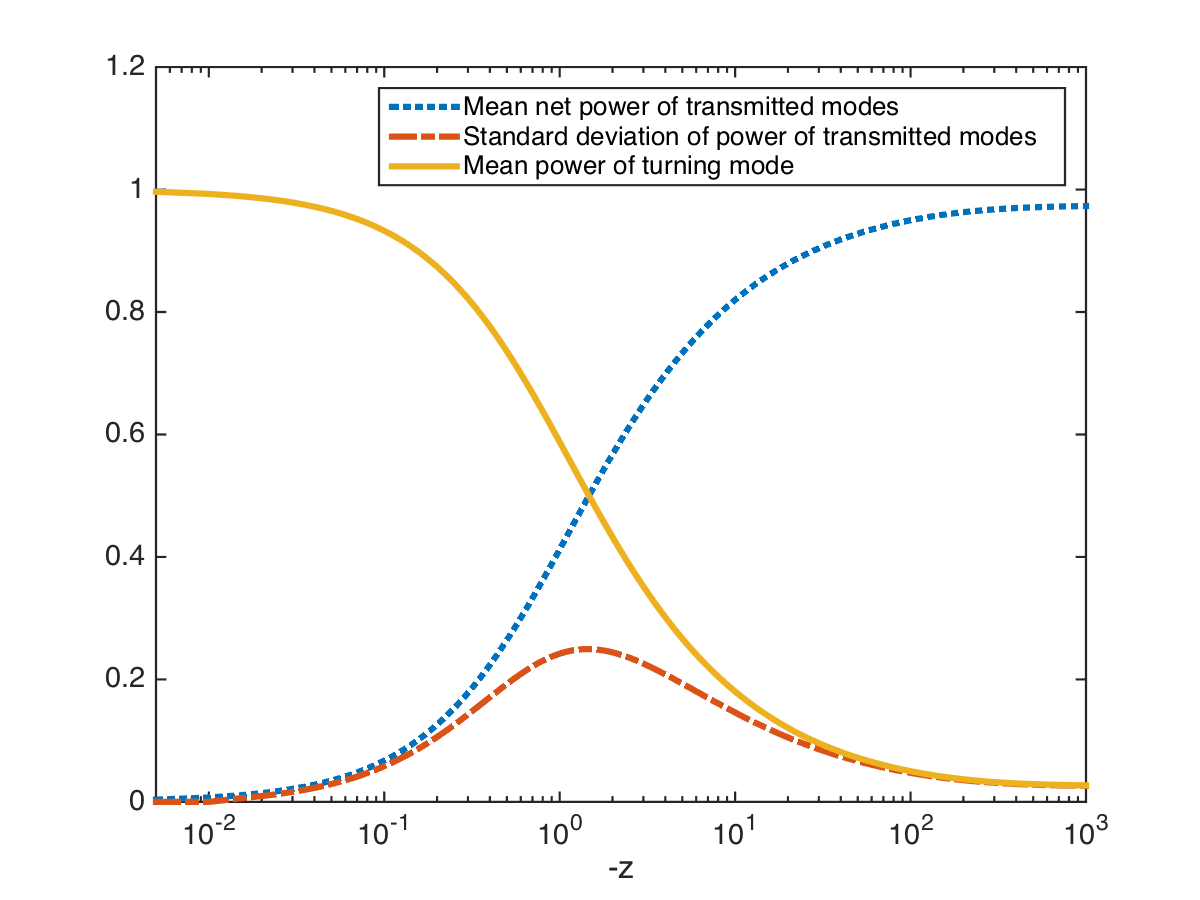}
\end{center}
\vspace{-0.1in}
\caption{ Display of the mean net power of the transmitted modes in
  dashed blue line, of the standard deviation of this power in dashed
  red line, and the mean power of the turning mode, indexed by $j =
  40$.  The abscissa is the arc length in units of $\la$ (in
  logarithmic scale). Only one mode was excited initially, the $39$-th
  one in the left plot and the $40$-th one in the right plot. }
\label{fig:Means1}
\end{figure}
The last illustration, in Figure \ref{fig:Means1}, shows the mean and
standard deviation of $\sum_{j=1}^{39} P_j(\om,z)$, and the mean power
$\big< P_{40}(\om,z)\big>$ of the turning mode, as functions of
$z$, for initial excitations of a single mode. In the left plot the
$39$-th mode is excited, and in the right plot the $40$-th
mode is excited. In the absence of the random fluctuations, these
initial conditions would determine the transmitted power at the
turning point.  Specifically, in the first case the power would stay
in the $39$-th mode and would propagate through, whereas in the second
case the power of the $40$-th mode would be totally reflected. The
cumulative scattering in the random waveguide distributes the power
among the modes, and we note in the left plot of Figure
\ref{fig:Means1} that slightly less power is transmitted, due to the
power transfer to the turning mode, whereas in the right plot, most of
the power is transmitted, due to the transfer of power from the
turning mode to the other modes.

\subsection{Universal transmission properties for strong scattering}
In case of strong scattering, the mean transmitted power through the left part of the waveguide becomes universal and equal to ${\cal P}_0 N_{\rm min}/N^{(0)} $,
where ${\cal P}_0 = \sum_{j=1}^{N^{(0)}} |b_{j,o}(\om)|^2$ is the power transmitted to the left by the source. 
More exactly, if scattering is so strong that equipartition is reached in each section between two turning points, 
in the sense that $z_-^{(t-1)} - z_-^{(t)} > L_{\rm eq}^t$ for all $t=1,\ldots,t_{M}^-$ (where $L_{\rm eq}^t$ is the equipartition distance in the section $(z_-^{(t)} , z_-^{(t-1)} ) $), then
 the fraction of mean power transmitted through the $t$-th turning point $z_-^{(t)}$ is $1-1/N^{(t-1)}_-$, 
 because the $N^{(t-1)}_-$-th mode carrying a fraction $1/N^{(t-1)}_-$ of the mean power is reflected. 
 By denoting $\left< {\cal P}^{(t-1)}_{\rm trans} \right>$ the net transmitted power in the $t$-th section $(z_-^{(t)},z_-^{(t-1)})$, 
 we get the recursive relation  
 \begin{equation}
 \left< {\cal P}^{(t)}_{\rm trans} \right>= \left< {\cal P}^{(t-1)}_{\rm trans} \right> (N^{(t-1)}_--1) / N^{(t-1)}_- , \quad 
 t = 1,\ldots,t_{M}^-,
 \end{equation} 
 which gives that the mean transmitted power at $-Z_{M}$ is ${\cal P}_0 N_{\rm min}/N^{(0)}$.
 
\section{Diffusion approximation theorem}
\label{app:adif}
In this section we state and prove the diffusion approximation theorem
used to obtain the asymptotic limit of the mode amplitudes in section
\ref{sect:results}.  We state the theorem for a general system of
random differential equations
 \begin{equation}
\frac{d {\bm{X}}^\ep(z)}{dz} = \frac{1}{\sqrt{\ep}} {\itbf F} \big(
     {\bm{X}}^\ep(z),q^\ep(z),\btheta^\ep(z),z\big) , \quad z > 0,
\label{eq:DIFF0}
\end{equation}
with unknown vector ${\bm{X}}^\ep \in \mathbb{R}^d$ satisfying the
initial condition ${\bm{X}}^\ep(0) = {\itbf x}_o$, and right hand side
defined by a function of the form
\begin{equation}
\label{eq:DIFF1}
  {\itbf F} \big( {\bm{X}},q,\btheta,z\big) = \sum_{j=1}^p {\itbf
    F}^{(j)} \big( {\bm{X}} ,q ,\theta_j ,z\big), \quad \mbox{for}~ ~
  \bm{\theta} = \big(\theta_j\big)_{j=1}^p \in \mathbb{R}^p.
\end{equation}

The second argument of ${\itbf F}$ is defined by $q^\ep(z) =
q(z/\ep)$, where $q(z)$ is a stationary and ergodic Markov
process
%\footnote{\textcolor{red}{ Does this result generalize to non
%    Markovian? Do we need to assume our fluctuations $\nu$ are
%    Markovian? Perhaps we can have a comment here?}} 
taking values in a space $E$, with generator $Q$ and stationary
distribution~$\pi_q$.  We assume that  $Q$ satisfies the
Fredholm alternative, which holds true for many different classes of
Markov processes \cite[section 6.3.3]{fouque07}.
Note that the Markovian assumption on the driving process $q$ is convenient for the proof, 
but the statement of the diffusion approximation theorem \ref{prop.diflim} 
generalizes to  a process $q$ that is not Markovian, but  $\phi$-mixing 
with $\phi \in L^{1/2}$ \cite[Sec. 4.6.2]{kushner}.

The third argument of ${\itbf F}$ is the vector valued function
$\bm{\theta}^\ep(z) $ taking values in $\mathbb{R}^p$, with components
satisfying the equation
\[
\frac{d \theta_j^\ep}{dz} = \frac{1}{\ep} \beta_j(z), \qquad j = 1,
\ldots, p,
\]
where $\beta_j(z)$ is a $\RR$-valued smooth function, bounded as  $C
\leq \beta_j(z) \leq 1/C$ for some constant $C>0$.

We assume that the components ${\itbf F}^{(j)}$ in \eqref{eq:DIFF1}
satisfy the following conditions, for all $j = 1, \ldots, p$:

\vspace{0.06in}
\begin{enumerate}
\itemsep 0.015in
\item The mappings $({\itbf x},z) \in \RR^d\times \RR\to {\itbf F}^{(j)}
  ({\itbf x},q,\theta_j,z)\in \RR^d $ are smooth for all $q\in E$
  and $\theta_j \in \RR$.
\item The mappings $ q \in E \to {\itbf F}^{(j)}({\itbf x},q,\theta_j,z) $
  are centered with respect to the stationary distribution $\pi_q$,
\[
\EE [ {\itbf F}^{(j)}({\itbf x},q(0),\theta,z) ] =\int_E {\itbf F}^{(j)}
({\itbf x},q,\theta_j,z) \pi_q(dq) =0 ,
\]
for any ${\itbf x} \in \RR^d$, $\theta_j \in \RR$ and $z \in \RR$. 
\item The mappings $ \theta_j \in \RR \to {\itbf F}^{(j)}({\itbf
  x},q,\theta_j,z) $ are periodic with period $1$ for all ${\itbf x} \in
  \RR^d$ and  $q\in E$.
\end{enumerate}

\vspace{0.06in}
\begin{theorem}
\label{prop.diflim}
Let ${\bm{X}}^\ep(z)$ be the solution of \eqref{eq:DIFF0}, with
right hand side $\bm{F}$ defined in terms of the functions
$\bm{F}^{(j)}$as in \eqref{eq:DIFF1}, and $\bm{F}^{(j)}$ satisfying
the three properties above.  In the limit $\ep \to 0$, the continuous
processes $({\bm{X}}^\ep(z))_{z \geq 0}$ converge in distribution to
the Markov diffusion process $({\bm{X}}(z))_{z \geq 0}$ with the
inhomogeneous generator
\begin{align}
\label{def:genX1}
{\cal L}_z f({\itbf x}) &= \sum_{m=1}^d {\textsl{h}}_m({\itbf x},z)
\partial_{x_m} f({\itbf x}) + \sum_{m,n=1}^d {\textsl{a}}_{m,n}({\itbf
  x},z) \partial_{x_m x_n}^2 f({\itbf x}) ,\\ {\textsl h}_m({\itbf x},z)
&= \sum_{n=1}^d \big< \int_{0}^{\infty} \hspace{-0.05in}\EE \left[ {F}_n ({\itbf
    x},q(0),\cdot ,z) \partial_{x_n} {F}_m({\itbf x},q(\zeta),\cdot
  +\bbeta(z) \zeta ,z ) \right] d\zeta \big>_{\hspace{-0.0in}\bbeta(z)} ,
 \label{def:genX4}
 \\ {\textsl a}_{m,n}({\itbf x},z) &=\big<
 \int_{0}^{\infty}\hspace{-0.05in} \EE \left[ {F}_n ({\itbf
     x},q(0),\cdot ,z) {F}_m({\itbf x},q(\zeta),\cdot +\bbeta(z) \zeta
   ,z ) \right] d\zeta \big>_{\hspace{-0.0in}\bbeta(z)} ,
\label{def:genX5}
\end{align}
where $\big< . \big>_\bbeta $ is the mean value for almost
periodic functions, 
$$ \big< {\itbf H}(\cdot) \big>_\bbeta = \lim_{S \to \infty}
\frac{1}{S} \int_0^S {\itbf H}(\btheta+\bbeta s) ds.
$$
\end{theorem}

Note that the mean values for the terms involved
in~(\ref{def:genX4}-\ref{def:genX5}) exist and are independent
of~$\btheta$, since the functions
\begin{align*}
G_{n,m}(s) &= \int_{0}^{\infty}\hspace{-0.05in} \EE \left[ {F}_n ({\itbf
    x},q(0),\btheta+\bbeta s ,z) {F}_m({\itbf x},q(\zeta), \btheta
  +\bbeta s +\bbeta \zeta ,z ) \right] d\zeta ,\\ \tilde{G}_{n,m}(s)
&= \int_{0}^{\infty} \hspace{-0.05in}\EE \left[ {F}_n ({\itbf
    x},q(0),\btheta+\bbeta s ,z) \partial_{x_n} {F}_m({\itbf
    x},q(\zeta), \btheta +\bbeta s +\bbeta \zeta ,z ) \right] d\zeta ,
\end{align*}
are periodic or almost periodic in $s$, for any fixed ${\itbf x}$ and $q$.

\vspace{0.06in}
\noindent \debproof Let us define the projection on the torus
$\mathbb{S} \simeq \RR/\mathbb{Z}$:
$$ \theta \in \RR \to \dot{\theta}=\theta \mbox{ mod } 1 \in
\mathbb{S},
$$ and observe that if a function $f(\theta)$ is periodic with
period $1$, then $f(\theta) = f(\dot{\theta})$.  We also define
$\dot{\btheta}^\ep(z) = \btheta^\ep(z)$ mod $1$, and $Z(z)=z$.  
The process $({\bm{X}}^\ep(z),q^\ep(z),\dot{\btheta}^\ep(z),Z(z))_{z
  \geq 0}$ is a Markov process with values in $\RR^d \times E \times
\mathbb{S}^p \times \RR$ and infinitesimal generator
\begin{equation}
\label{def:Leps}
{\cal L}^\ep = \frac{1}{\ep} \big( Q + \bbeta(Z) \cdot
\nabla_{\dot{\btheta}}\big) +\frac{1}{\sqrt{\ep}} {\itbf
  F}({\bm{X}},q,\dot{\btheta},Z) \cdot \nabla_{\bm{X}}+ {\partial_Z}
\, .
\end{equation}

One can show by the perturbed test function method \cite[Section
  6.3.2]{fouque07} and Lemma \ref{lem:2} that the continuous processes
$({\bm{X}}^\ep(z),Z(z))_{z \geq 0}$ converge in distribution to the
Markov diffusion process $({\bm{X}}(z),Z(z))_{z \geq 0}$ with the
homogeneous generator:
\begin{align}
&{\cal L} f({\itbf x},Z) ={\partial_Z} f({\itbf x},Z) \label{eq:Lepsb}
  \\ &\hspace{0.05in}+ \big< \int_{0}^{\infty} \hspace{-0.05in} \EE
  \left[ {\itbf F}({\itbf x},q(0),\cdot,Z ) \cdot \nabla_{\itbf x} \left(
    {\itbf F}({\itbf x},q(\zeta),\cdot +\bbeta(Z) \zeta ,Z) \cdot
    \nabla_{\itbf x} f({\itbf x},Z) \right) \right] d\zeta
  \big>_{\hspace{-0.0in}\bbeta(z)}.  \nonumber
\end{align}
Since $(Z(z))_{z \geq 0}$ is deterministic, we conclude that
$({\bm{X}}(z) )_{z \geq 0}$ is a Markov process and its inhomogeneous
generator is
\begin{equation}
\hspace{-0.1in}{\cal L}_z f({\itbf x}) =\big< \int_{0}^{\infty} \hspace{-0.05in}\EE
\left[ {\itbf F} ({\itbf x},q(0),\cdot ,z) \cdot \nabla_{\itbf x} \left(
  {\itbf F}({\itbf x},q(\zeta),\cdot +\bbeta(z) \zeta ,z ) \cdot
  \nabla_{\itbf x} f({\itbf x}) \right) \right] d\zeta \big>_{\hspace{-0.0in}\bbeta(z)}
\label{eq:L}
\end{equation}
or equivalently (\ref{def:genX1}).
\finproof

\vspace{0.1in}
\begin{lemma}
\label{lem:1}%
 We have the following two statements:

\vspace{0.06in}
1. Let $\beta \in \RR \backslash \{0\}$. Let $g(q,\theta)$ be a bounded
function, periodic in $\theta \in \RR$ with period $1$, such that
$$
 \EE [ g(q(0),\theta)]   =0 \mbox{ for all } \theta \in \RR \, .
$$ 
The Poisson equation
$$
 \big( Q + \beta  \partial_{\dot{\theta}}\big) f = g 
$$ has a unique solution $f$, periodic in $\theta$, up to an additive
 constant. The solution with mean zero is
\begin{equation}\label{eq:formf1}
f(q,\dot{\theta}) = - \int_0^\infty \hspace{-0.05in}\EE [ g(q(\zeta),
  \dot{\theta}+\beta \zeta) | q(0)=q] d\zeta \, .
\end{equation}

\vspace{0.06in}
2. Let $\bbeta \in \RR^2$ with non-zero entries.  Let $g(q,\btheta)$ be a
bounded function, periodic in $\btheta \in \RR^2$ with period $1$,
such that
$$ \EE [ g(q(0),\btheta)] =0 \mbox{ for all } \btheta \in \RR^2 \, .
$$ 
The Poisson equation
$$ \big( Q + \bbeta \cdot \nabla_{\dot{\btheta}}\big) f = g
$$ has a unique solution $f$, periodic in $\btheta$, up to an additive
constant.  The solution with mean zero is
\begin{equation}
\label{eq:formf2}
f(q,\dot{\btheta}) = - \int_0^\infty \hspace{-0.05in}\EE [ g(q(\zeta),
  \dot{\btheta}+\bbeta \zeta) | q(0)=q] d\zeta \, .
\end{equation}
\end{lemma}

\vspace{0.06in} Note that in the second item of Lemma \ref{lem:1} it
is important to assume that $ \EE [ g(q(0),\btheta)] =0 \mbox{ for all
} \btheta \in \RR^2$, and not only that $\int_{\mathbb{S}^2} \EE [
  g(q(0),\dot{\btheta} )] d\dot{\btheta} =0$.  The latter weaker
hypothesis ensures the desired result only when $\beta_1 / \beta_2$ is
irrational.

\vspace{0.06in}
\noindent
\debproof To prove statement 1. let $\beta \in \RR$ be fixed. We denote by
$\theta_\beta(\zeta)$ the solution to
$\frac{d\theta_\beta}{d\zeta}=\beta$ and by $\dot{\theta}_\beta(\zeta)
= \theta_\beta(\zeta)\mbox{ mod } 1$.  The process $(q(\zeta),
\dot{\theta}_\beta(\zeta))_{\zeta \geq 0}$ is a Markov process with
values in $E \times \mathbb{S}$ and with generator $Q + \beta
\partial_{\dot{\theta}}$.  It is a stationary process with the
stationary distribution $\pi_q \otimes \nu_\mathbb{S}$ where
$\nu_\mathbb{S}$ is the uniform distribution over the torus
$\mathbb{S}$.  It is also an ergodic process with respect to the
stationary distribution $\pi_q \otimes \nu_\mathbb{S}$.  Since $g$
satisfies $\int g(q,\dot{\theta}) \pi_q(dq)=0$ for all $\dot{\theta}$,
it a fortiori satisfies $\iint g(q,\dot{\theta}) \pi_q(dq)
\nu_\mathbb{S}(d\dot{\theta})=0$, and the result then follows from
standard arguments \cite[section 6.5.2]{fouque07}:
$$ f(q_0,\dot{\theta}_0) = - \int_0^\infty \hspace{-0.05in}\EE \big[
  g( q(\zeta), \dot{\theta}_\beta(\zeta) ) | q(0)=q_0 ,
  \dot{\theta}_\beta(0)=\dot{\theta}_0 \big] d\zeta \, ,
$$
which gives (\ref{eq:formf1}).

To prove statement 2. let $\bbeta \in \RR^2$ be fixed.  We denote by
$\btheta_\bbeta(\zeta)$ the solution to
$\frac{d\btheta_\bbeta}{d\zeta}=\bbeta$ and by
$\dot{\btheta}_\bbeta(\zeta) = \btheta_\bbeta(\zeta)\mbox{ mod } 1$.
The process $(q(\zeta), \dot{\btheta}_\bbeta(\zeta))_{\zeta \geq 0}$
is a Markov process with values in $E \times \mathbb{S}^2$ and with
generator $Q + \bbeta \cdot \nabla_{\dot{\btheta}}$.

If the ratio $\beta_1/\beta_2$ of the entries of $\beta_2$ of $\bbeta$
is irrational, the process $(q(\zeta),
\dot{\btheta}_\bbeta(\zeta))_{\zeta \geq 0}$ is stationary and
ergodic, with the stationary distribution $\pi_q \otimes
\nu_{\mathbb{S}^2}$, where $\nu_{\mathbb{S}^2}$ is the uniform
distribution over the torus $\mathbb{S}^2$.  Since $g$ satisfies $\int
g(q,\dot{\btheta}) \pi_q(dq)=0$ for all $\dot{\btheta}$, it a fortiori
satisfies $\iint g(q,\dot{\btheta}) \pi_q(dq)
\nu_{\mathbb{S}^2}(d\dot{\btheta})=0$, and the result then follows
from standard arguments \cite[section 6.5.2]{fouque07}:
$$ f(q_0,\dot{\btheta}_0) = - \int_0^\infty \hspace{-0.05in} \EE [ g(
  q(\zeta), \dot{\btheta}_\bbeta(\zeta) ) | q(0)=q_0 ,
  \dot{\btheta}_\bbeta(0)=\dot{\btheta}_0 ] d\zeta \, ,
$$
which gives (\ref{eq:formf2}).

If the ratio $\beta_1/\beta_2$ of the entries of  $\bbeta$
is rational, that is to say, if there exist nonzero integers $n_1,n_2$
such that $n_1 \beta_1=n_2\beta_2$, then $(
\dot{\btheta}_\bbeta(\zeta) )_{\zeta \geq 0}$ is not ergodic over the
torus $\mathbb{S}^2$.  However, for a given starting point
$\dot{\btheta}_0$, it satisfies the ergodic theorem over the compact
manifold $\mathbb{S}^1_{\dot{\btheta}_0} := \{ \dot{\btheta}_0+\bbeta
s \mbox{ mod } 1, \ s\in \RR\}$, with the uniform distribution over
the manifold $\mathbb{S}^1_{\dot{\btheta}_0}$.  Since $g$ satisfies
$\int g(q,\dot{\btheta}) \pi_q(dq)=0$ for all $\dot{\btheta}$, it a
fortiori satisfies $\iint g(q,\dot{\btheta}) \pi_q(dq)
\nu_{\mathbb{S}^1_{\dot{\btheta}_0}}(d\dot{\btheta})=0$.  We can then
define
$$ f(q_0,\dot{\btheta}_0) = - \int_0^\infty \EE [ g( q(\zeta),
  \dot{\btheta}_\bbeta(\zeta) ) | q(0)=q_0 ,
  \dot{\btheta}_\bbeta(0)=\dot{\btheta}_0 ] d\zeta \, ,
$$
which gives (\ref{eq:formf2}).
\finproof

\vspace{0.06in}
We can now state the lemma used in the proof of Theorem \ref{propdiff}:

\vspace{0.06in}
\begin{lemma}
\label{lem:2}
For all $ f \in {\cal C}_b^\infty(\RR^d \times \RR ,\RR)$, and all
compact sets $ K $ of $\RR^d \times \RR$, there exists a family
$f^\ep$ such that:
\begin{eqnarray}
\label{pertf0}
&&\sup_{({\itbf x},Z )\in K,q \in E,\dot{\btheta}\in \mathbb{S}^p}
|f^\ep({\itbf x},q,\dot{\btheta},Z) - f({\itbf x},Z)| \stackrel{\ep
  \rightarrow 0}{\longrightarrow} 0 , \\ && \sup_{({\itbf x},Z ) \in
  K,q\in E,\dot{\btheta} \in \mathbb{S}^p} |{\cal L}^\ep f^\ep({\itbf
  x},q,\dot{\btheta},Z) - {\cal L} f({\itbf x},Z)| \stackrel{\ep
  \rightarrow 0}{\longrightarrow} 0 ,
\label{pertf0bis}
\end{eqnarray}
where ${\cal L}^\ep$ is the generator (\ref{def:Leps}) and ${\cal L}$
is the generator (\ref{eq:Lepsb}).
\end{lemma}

\vspace{0.06in} \noindent \debproof Let $f \in {\cal C}_b^\infty(\RR^d
\times \RR ,\RR)$, and define
\begin{equation}
\label{eq:deff}
 f^\ep({\itbf x},q,\dot{\btheta},Z ) = f({\itbf x},Z) + \sqrt{\ep}
 f_{1}({\itbf x},q,\dot{\btheta},Z) + \ep f_{2}({\itbf
   x},q,\dot{\btheta},Z) + \ep f_{3}^\ep ({\itbf x},\dot{\btheta},Z ) ,
\end{equation} 
where $f_1$, $f_2$, and $f_3^\ep$ will be specified later on.
Applying ${\mathcal L}^\ep$ to $f^\ep$, we get
\begin{align}
{\mathcal L}^\ep f^\ep({\itbf x},q,\dot{\btheta},Z ) &=
\frac{1}{\sqrt{\ep}} \left( \big( Q + \bbeta(Z) \cdot
\nabla_{\dot{\btheta}}\big) f_1 + {\itbf F}({\itbf
  x},q,\dot{\btheta},Z)\cdot\nabla_{\itbf x} f({\itbf x},Z) \right)
\nonumber\\ & + \left( \big( Q + \bbeta(Z) \cdot
\nabla_{\dot{\btheta}}\big) f_2 + {\itbf F}({\itbf x},q,\dot{\btheta},Z)
\cdot \nabla_{\itbf x} f_1 ({\itbf x},q,\dot{\btheta},Z) \right)
\nonumber\\ & + \bbeta(Z) \cdot \nabla_{\dot{\btheta}} f_3^\ep ({\itbf
  x},\dot{\btheta},Z) + {\partial_Z} f({\itbf x},Z)+ O(\sqrt{\ep})
. \label{eq:Lep1}
\end{align}

Now let us define the correction $f_1$ as
\begin{align}
f_1({\itbf x},q,\dot{\btheta} ,Z) = \sum_{j=1}^p f_{1}^{(j)} ({\itbf
  x},q,\dot{\theta}_j ,Z) , \label{eq:f1}
\end{align}
where 
\begin{align*}
f_{1}^{(j)} ({\itbf x},q,\dot{\theta}_j ,Z) = - \big( Q + \beta_j(Z)
\partial_{\dot{\theta}_j}\big)^{-1} \left( {\itbf F}^{(j)} ({\itbf
  x},q,\dot{\theta}_j,Z) \cdot \nabla_{\itbf x} f({\itbf x},Z) \right) .
\end{align*}
These functions are well-defined and admit the representation
$$ f_{1}^{(j)}({\itbf x},q,\dot{\theta}_j,Z ) =
\int_{0}^{\infty} \hspace{-0.05in} \EE \big[{\itbf F}^{(j)}({\itbf
    x},q(\zeta),\dot{\theta}_j+\beta_j(Z)\zeta ,Z) \cdot \nabla_{\itbf
    x} f({\itbf x},Z) | q(0)=q \big] d \zeta ,
$$
by Lemma \ref{lem:1}.

The second correction $f_2$ is defined by
\begin{equation}
f_{2}({\itbf x},q,\dot{\btheta} ,Z) = \sum_{j,l=1}^p f_{2}^{(jl)} ({\itbf
  x},q,\dot{\theta}_j,\dot{\theta}_l ,Z) ,\label{eq:f2}
\end{equation}
where
\begin{align*}
&f_{2}^{(jl)}({\itbf x},q,\dot{\theta}_j,\dot{\theta}_l ,Z) = - \big( Q
+ \beta_j(Z) \partial_{\dot{\theta}_j} + \beta_l(Z)
\partial_{\dot{\theta}_l}\big)^{-1} \\  &\times\hspace{-0.02in} \Big( {\itbf F}^{(j)}({\itbf
  x},q,\dot{\theta}_j,Z) \cdot \nabla_{\itbf x} f_{1}^{(l)}({\itbf
  x},q,\dot{\theta}_l,Z)  - \EE \big[ {\itbf F}^{(j)} ({\itbf
    x},q(0),\dot{\theta}_j,Z) \cdot \nabla_{\itbf x} f_{1}^{(l)} ({\itbf
    x},q(0),\dot{\theta}_l,Z) \big] \Big) .
\end{align*}
These functions  are well defined by Lemma \ref{lem:1} since
the argument of the operator $ \big( Q + \beta_j(Z) \partial_{\dot{\theta}_j} +
\beta_l(Z) \partial_{\dot{\theta}_l}\big)^{-1}$ has mean zero for all
$\btheta$.  

Substituting \eqref{eq:f1} and \eqref{eq:f2} in \eqref{eq:Lep1} we obtain
\begin{align}
{\mathcal L}^\ep f^\ep({\itbf x},q,\dot{\btheta},Z ) =\sum_{j,l=1}^p
g_{3}^{(jl)} ({\itbf x}, \dot{\theta}_j, \dot{\theta}_l ,Z) + \bbeta(Z)
\cdot \nabla_{\dot{\btheta}} f_3^\ep ({\itbf x},\dot{\btheta},Z)\nonumber \\ +
      {\partial_Z} f({\itbf x},Z)+ O(\sqrt{\ep}) ,\label{eq:Lep2}
\end{align} 
with
\begin{equation} 
 g_{3}^{(jl)} ({\itbf x},\dot{\theta}_j, \dot{\theta}_l,Z) = \EE \big[
   {\itbf F}^{(j)} ({\itbf x},q(0),\dot{\theta}_j,Z) \cdot \nabla_{\itbf
     x} f_{1}^{(l)} ({\itbf x},q(0),\dot{\theta}_l,Z ) \big] .
\label{eq:defgep}
 \end{equation} 
We now define the third correction function
\begin{equation}
f_{3}^\ep({\itbf x},\dot{\btheta} ,Z ) = \sum_{j,l=1}^p
f_{3}^{(jl),\ep} ({\itbf x},\dot{\theta}_j,\dot{\theta}_l ,Z ) ,
\end{equation}
with terms
\begin{equation*}
f_{3}^{(jl),\ep}({\itbf x},\dot{\theta}_j, \dot{\theta}_l,Z) =
\int_0^\infty e^{-\sqrt{\ep} s} \,\widetilde{g}_{3}^{(jl)} ({\itbf
  x},\dot{\theta}_j+\beta_j (Z) s , \dot{\theta}_l +\beta_l (Z)s ,Z)
ds,
\end{equation*}
defined by 
\begin{align*}\widetilde{g}_{3}^{(jl)} ({\itbf x},\dot{\theta}_j,
\dot{\theta}_l,Z) &= g_{3}^{(jl)} ({\itbf x},\dot{\theta}_j,
\dot{\theta}_l,Z) - \mathcal{G}_{3}^{(jl)} ({\itbf x}, Z),
\end{align*}
where 
\begin{align} 
\mathcal{G}_3^{(jl)} ({\itbf
  x}, Z) &= \lim_{S \to \infty} \frac{1}{S} \int_0^S g_{3}^{(jl)}
({\itbf x},\dot{\theta}_j+\beta_j (Z) s , \dot{\theta}_l +\beta_l
(Z)s,Z) ds .\label{eq:defGep}
 \end{align}
These are are well defined because 
$s \mapsto g_{3}^{(jl)} ({\itbf x},\dot{\theta}_j+\beta_j (Z) s ,
\dot{\theta}_l +\beta_l (Z)s,Z)$ are almost periodic mappings.  

Note that $\sqrt{\ep} f_{3}^{(jl),\ep}$ is uniformly bounded because
$\widetilde{g}_{3}^{(jl)} $ is bounded. This and definitions
\eqref{eq:f1}, \eqref{eq:f2} of the corrections $f_1$ and $f_2$ used
in equation \eqref{eq:deff} imply that $f^\ep$ satisfies
(\ref{pertf0}).  Note also that $\sqrt{\ep} f_{3}^{(jl),\ep}$ goes to
zero as $\ep \to 0$, because the mapping $s \mapsto
\widetilde{g}_{3}^{(jl)} ({\itbf x},\dot{\theta}_j+\beta_j (Z) s ,
\dot{\theta}_l +\beta_l (Z)s ,Z)$ is almost periodic and with mean
zero. Moreover, using the chain rule and integration by parts, we
obtain
\begin{align*}
 \big( \beta_j(Z) \partial_{\dot{\theta}_j} + \beta_l(Z)
 \partial_{\dot{\theta}_l}\big) &f_{3}^{(jl),\ep}({\itbf
   x},\dot{\theta}_j, \dot{\theta}_l,Z) \\ &= \int_0^\infty e^{-\sqrt{\ep} s}
 \partial_s \big[ \widetilde{g}_{3}^{(jl)} ({\itbf
     x},\dot{\theta}_j+\beta_j (Z) s , \dot{\theta}_l +\beta_l (Z)s
   ,Z) \big] ds\\ &= - \widetilde{g}_{3}^{(jl)} ({\itbf
   x},\dot{\theta}_j, \dot{\theta}_l,Z) + \sqrt{\ep}
 f_{3}^{(jl),\ep}({\itbf x},\dot{\theta}_j, \dot{\theta}_l,Z) .
\end{align*}

Gathering the results, equation \eqref{eq:Lep2} becomes 
\begin{eqnarray*}
{\mathcal L}^\ep f^\ep &=&\sum_{j,l=1}^p \mathcal{G}_{3}^{(jl)} ({\itbf
  x} ,Z) + {\partial_Z} f({\itbf x},Z) +\sqrt{\ep} f_{3}^\ep ({\itbf x},
\dot{\btheta},Z) + O(\sqrt{\ep}) .
\end{eqnarray*} 
The result \eqref{pertf0bis} follows from this equation and definitions
\eqref{eq:defgep}, \eqref{eq:defGep} and \eqref{eq:f1}, because $\sqrt{\ep}
f_3^\ep$ goes to zero as $\ep \to 0$.  \finproof

\section{Summary}
\label{sect:sum}
We studied the transmission and reflection of time-harmonic sound
waves emitted by a point source in a two-dimensional random waveguide
with turning points. The waveguide has sound soft boundaries, a slowly
bending axis and variable cross-section. The variation consists of a
slow and monotone change of the opening $D$ of the waveguide, and small
amplitude random fluctuations of the boundary.  The slow variations
are on a long scale with respect to the wavelength $\la$, whereas the
random fluctuations are on a scale comparable to $\la$. The  wavelength $\la$ is chosen smaller than $D$,
so that the wave field is a
superposition of multiple propagating modes, and infinitely many
evanescent modes. The turning points are the locations along the axis
of the waveguide where the number of propagating modes decreases by
$1$ in the direction of decrease of $D$, or increases by $1$ in the
direction of increase of $D$. The change in the number of propagating
modes means that there are modes that transition from propagating to
evanescent. Due to energy conservation, the incoming such waves are
turned back i.e., they are reflected at the turning points.

We analyzed the transmitted and reflected
propagating modes in the waveguide and quantified their interaction
with the random boundary. This interaction is called cumulative
scattering and it manifests as mode coupling which causes
randomization of the wave field and exchange of power between the
modes. We analyzed these effects from first principles, 
starting from the wave equation, using stochastic asymptotic analysis.
We focused attention on the transport of power in the
waveguide and showed that cumulative scattering may increase or decrease
the transmitted power depending on the source.

% ---------------------------------------
\section*{Acknowledgements}
The research of LB and DW was supported in part by NSF grant
DMS1510429.  LB also acknowledges support from AFOSR grant
FA9550-15-1-0118.

% ----------------
\appendix
\section{Transformation to curvilinear coordinates}
\label{ap:CurvCoord}
The Frenet-Serret formulas give
\begin{align*}
\partial_z \xpar(z) &= \btau\Big(\frac{z}{L}\Big),
\quad  \partial_z \btau\Big(\frac{z}{L}\Big) &= \frac{1}{L} \kappa
\Big(\frac{z}{L}\Big) \bn \Big(\frac{z}{L}\Big), 
\quad \partial_z \bn\Big(\frac{z}{L}\Big) &= -\frac{1}{L} \kappa
\Big(\frac{z}{L}\Big) \btau \Big(\frac{z}{L}\Big),
\end{align*}
and from \eqref{eq:coord1} we obtain that the vectors 
$
\partial_r \bx = \bn\Big(\frac{z}{L}\Big)$ and $\partial_z \bx =
\left[ 1 - \frac{r}{L} \kappa \Big( \frac{z}{L}\Big)\right]
\btau\Big(\frac{z}{L}\Big)$
are orthogonal. Their norm defines the Lam\'{e} coefficients 
$
h_r = |\partial_r \bx| = 1$ and $h_z = |\partial_z \bx| = \Big|1 -
\frac{r}{L} \kappa \Big(\frac{z}{L}\Big)\Big|,$
which in turn define the Laplacian operator  in curvilinear coordinates \cite{spiegel1959schaum}
$
\Delta = \frac{1}{h_rh_z} \left[ \partial_r \Big(\frac{h_z}{h_r}
  \partial_r \Big) + \partial_x \Big(\frac{h_r}{h_z} \partial_z \Big)\right].
$
This is written explicitly in the left hand-side of equation
\eqref{eq:we3}. For the right hand-side we used the formula 
$
\delta(\bx-\bx_\star) = \frac{1}{h_r h_z} \delta(z) \delta(r-r_\star).
$

\section{Derivation of the asymptotic model}
\label{ap:DerAs}
We begin with the region $|z| < Z_M$. 
The change of variables
\eqref{eq:AM3} and the chain rule give
\begin{align*}
\partial_r p &= \frac{\partial_\rho p^\ep(\om,\rho,z)}{1+
  \frac{\sqrt{\ep}}{2} \sigma \nu \big(\frac{z}{\ep}\big)}, \quad 
\partial_r^2 p = \frac{\partial_\rho^2 p^\ep(\om,\rho,z)}{\big[1+
    \frac{\sqrt{\ep}}{2} \sigma \nu
    \big(\frac{z}{\ep}\big)\big]^2},
\end{align*} 
and 
\begin{align*}
 \partial_z p &= \Big\{\partial_z -
\frac{\big[ [2 \rho + D(z)] \frac{\sigma}{\sqrt{\ep}} \nu'
    \big(\frac{z}{\ep}\big) + D'(z) \sqrt{\ep} \sigma \nu
    \big(\frac{z}{\ep}\big)\big]}{2\big[2+ \sqrt{\ep} \sigma
    \nu \big(\frac{z}{\ep}\big)\big]} \partial_\rho\Big\}
p^\ep(\om,\rho,z), \\ \partial_z^2 p &= \Big\{ \partial_z -
\frac{\big[ [2 \rho + D(z)] \frac{\sigma}{\sqrt{\ep}} \nu'
    \big(\frac{z}{\ep}\big) + D'(z) \sqrt{\ep} \sigma \nu
    \big(\frac{z}{\ep}\big)\big]}{2\big[2+ \sqrt{\ep} \sigma
    \nu \big(\frac{z}{\ep}\big)\big]}\partial_\rho \Big\}^2
p^\ep(\om,\rho,z).
\end{align*}
Substituting in \eqref{eq:AM1} we get 
\begin{align*} 
&\partial_z^2 p^\ep(\om,\rho,z) + \frac{\big[ [2 \rho +
          D(z)]\frac{\sigma}{\sqrt{\ep}} \nu' \big(\frac{z}{\ep}\big) +
      {D'(z)} \sqrt{\ep} \sigma \nu
      \big(\frac{z}{\ep}\big)\big]^2}{4\big[2+ \sqrt{\ep}
      \sigma \nu \big(\frac{z}{\ep}\big)\big]^2} \partial_\rho^2
  p^\ep(\om,z) \\ &- \frac{\big[[2 \rho +
          D(z)]\frac{\sigma}{\sqrt{\ep}} \nu' \big(\frac{z}{\ep}\big) +
      {D'(z)} \sqrt{\ep} \sigma \nu
      \big(\frac{z}{\ep}\big)\big]}{\big[2+ \sqrt{\ep} \sigma
      \nu \big(\frac{z}{\ep}\big)\big]} \partial_{\rho z}^2 p^\ep(\om,z)
  \\ &+ \frac{\Big\{ [2 \rho + D(z)]\frac{\sigma^2}{\ep}{\nu'}^{\,2}
        \big(\frac{z}{\ep}\big) + {D'(z)}\sigma^2 \nu'
      \big(\frac{z}{\ep}\big)\nu \big(\frac{z}{\ep}\big)\Big\}}{\big[2+
      \sqrt{\ep} \sigma \nu \big(\frac{z}{\ep}\big)\big]^2}
  \partial_\rho p^\ep(\om,\rho,z) \\&- \frac{[2 \rho +
        D(z)]\frac{\sigma}{\ep^{3/2}} \nu''\big(\frac{z}{\ep}\big) +
    2 D'(z) \frac{\sigma}{\sqrt{\ep}} \nu'\big(\frac{z}{\ep}\big) +
    D''(z) \sqrt{\ep} \sigma
    \nu\big(\frac{z}{\ep}\big)}{2\big[2+ \sqrt{\ep} \sigma \nu
    \big(\frac{z}{\ep}\big)\big]} \partial_\rho p^\ep(\om,\rho,z) \\ & +
  \frac{\Big\{ 1 - \ep \kappa(z)\big[ \rho + \frac{[2 \rho +
          D(z)]}{4}\sqrt{\ep} \sigma \nu
      \big(\frac{z}{\ep}\big)\big]\Big\}^2}{\ep^2} \Big\{
  \frac{\partial_\rho^2 p^\ep(\om,\rho,z)}{\big[1+ \frac{\sqrt{\ep}}{2}
      \sigma \nu \big(\frac{z}{\ep}\big)\big]^2} + k^2
  p^\ep(\om,\rho,z)\Big\}\\ &-\frac{\kappa(z) \Big\{ 1 - \ep \kappa(z) \big[
      \rho + \frac{[2 \rho + D(z)]}{4}\sqrt{\ep} \sigma \nu
      \big(\frac{z}{\ep}\big)\big]\big\}}{\ep[1+
      \frac{\sqrt{\ep}}{2} \sigma \nu \big(\frac{z}{\ep}\big)]}
  \partial_\rho p^\ep(\om,\rho,z) \\ &+\frac{\ep \kappa'(z) \big[ \rho
      + \frac{[2 \rho + D(z)]}{4}\sqrt{\ep} \sigma \nu
      \big(\frac{z}{\ep}\big)\big]}{\big\{1 - \ep \kappa(z)\big[ \rho + \frac{[2
            \rho + D(z)]}{4}\sqrt{\ep} \sigma \nu
        \big(\frac{z}{\ep}\big)\big]\big\}} \Big\{\partial_z
  p^\ep(\om,\rho,z) \\&\hspace{1.3in}-\frac{\big[ [2 \rho + D(z)]
        \frac{\sigma}{\sqrt{\ep}} \nu' \big(\frac{z}{\ep}\big) +
      D'(z) \sqrt{\ep} \sigma \nu
      \big(\frac{z}{\ep}\big)\big]}{2\big[2+ \sqrt{\ep} \sigma
      \nu \big(\frac{z}{\ep}\big)\big]} \partial_\rho
  p^\ep(\om,\rho,z)\Big\} \\ &= \frac{f(\om)\Big\{ 1 -
    \ep \Big[\rho_\star + \frac{[2 \rho + D(0)]}{4}\sqrt{\ep} \sigma
      \nu(0)\Big] \Big\} }{\ep \Big[1 + \frac{\sqrt{\ep}}{2} \sigma
        \nu(0)\Big]} \delta(\rho-\rho_\star) \delta(z).
\end{align*}
By assumption $\nu$, $\nu'$ and $\nu''$ are bounded
almost surely. Moreover, $\kappa'$ and $D'$, $D''$ are bounded
uniformly in $\mathbb{R}$. Thus, we can expand the coefficients of the
differential operator in powers of $\ep$ and obtain after multiplying
through by $\ep$,
\begin{align}
&\frac{1}{\ep}\Big[(\ep \partial_z)^2 +\partial_\rho^2 + k^2
    \Big]p^\ep(\om,\rho,z) - 2 \rho \kappa(z)\big[1 +
    O(\sqrt{\ep})\big] \big(\partial_\rho^2 + k^2) \nonumber \\
&- \kappa(z) \big[1 +
    O(\sqrt{\ep})\big]
  \partial_\rho p^\ep(\om,\rho,z)-\ep \kappa'(z)\big[1 +
      O(\sqrt{\ep})\big] (\ep \partial_z)p^\ep(\om,\rho,z)\nonumber\\ &-
    \frac{[2 \rho + D(z)]}{2} \Big[ \frac{\sigma}{\sqrt{\ep}}
      \nu'\Big(\frac{z}{\ep}\Big) - \frac{\sigma^2}{2}
      \nu'\Big(\frac{z}{\ep}\Big) \nu\Big(\frac{z}{\ep}\Big) +
      O(\sqrt{\ep}) \Big] \ep \partial_{\rho z}^2 p^\ep(\om,\rho,z)
    \nonumber\\ &- \Big\{ \frac{\sigma}{\sqrt{\ep}}
    \nu\Big(\frac{z}{\ep}\Big) - \frac{3 \sigma^2}{4}
    \nu^2\Big(\frac{z}{\ep}\Big) - \frac{[2 \rho +
        D(z)]^2\sigma^2}{16}{\nu'}^{\,2} \Big(\frac{z}{\ep}\Big)
    + O(\sqrt{\ep}) \Big\}\partial_\rho^2
    p^\ep(\om,\rho,z) \nonumber\\ &-\frac{\big[2 \rho + D(z)\big]}{4} \Big\{
    \frac{\sigma}{\sqrt{\ep}} \nu''\Big(\frac{z}{\ep}\Big) -
    \frac{\sigma^2}{2}
    \nu''\Big(\frac{z}{\ep}\Big)\nu\Big(\frac{z}{\ep}\Big) -
    \sigma^2{\nu'}^{\,2}\Big(\frac{z}{\ep}\Big) + O(\sqrt{\ep})
    \Big\} \partial_\rho p^\ep(\om,\rho,z)\nonumber\\&= f(\om)\Big[1 +
      O(\sqrt{\ep})\Big] \delta(\rho-\rho_\star) \delta(z), \label{eq:APC1}
\end{align}
for $|z| < Z_M$. This is the asymptotic series in \eqref{eq:AM5}.

The result simplifies at $|z| > Z_M$, where the waveguide has no variations,
as stated at the end of section \ref{sect:form3}.

\section{Derivation of the mode coupling equations}
\label{ap:ModeCoup}
Substituting \eqref{eq:DC1} in \eqref{eq:AM5}, taking the inner
product with $y_j(\rho,z)$ and using the identities
\eqref{eq:orthog}-\eqref{eq:D6} we obtain the following system of
equations for the modes
\begin{align}
&\frac{1}{\ep} \Big[ \big(\ep \partial_z)^2 + k^2 - \mu_j^2(z)\Big]
  p_j^\ep(\om,z) +\frac{\sigma}{\sqrt{\ep}} \Big[ \mu_j^2(z) \nu
    \Big(\frac{z}{\ep}\Big) + \frac{1}{4} \nu''\Big(\frac{z}{\ep}\Big)
    + \frac{1}{2} \nu'\Big(\frac{z}{\ep}\Big) \ep \partial_z\Big]
  p_j^\ep(\om,z) \nonumber\\ &\hspace{0.76in}- \frac{\sigma^2}{4}
  \Big\{ 3 \mu_j^2(z) \nu^2\Big(\frac{z}{\ep}\Big) + \Big[\frac{(\pi
      j)^2}{3} + \frac{1}{2}\Big]
    {\nu'}^{\,2}\Big(\frac{z}{\ep}\Big) + \frac{1}{2}
  \nu\Big(\frac{z}{\ep}\Big) \nu''\Big(\frac{z}{\ep}\Big)\Big]
    p_j^\ep(\om,z)\nonumber\\ &\hspace{0.76in}-\frac{\sigma^2}{4}
    \nu\Big(\frac{z}{\ep}\Big) \nu'\Big(\frac{z}{\ep}\Big)\ep
    \partial_z p_j^\ep(\om,z) \approx C_j^\ep(\om,z) + f(\om)
    y_j(\rho_\star,0) \delta(z), \label{eq:MC1}
\end{align}
at $|z| < Z_M$, where the approximation is because we neglect the
$O(\sqrt{\ep})$ terms that vanish in the limit $\ep \to 0$. The
coupling term is
\begin{align}
C_j^\ep(\om,z) = \sum_{q=1, q \ne j}^\infty & \Big\{ \frac{2 j q
  (-1)^{j+q}}{(q^2-j^2)} \Big[ \frac{\sigma}{\sqrt{\ep}}
  \nu'\Big(\frac{z}{\ep}\Big) - \frac{\sigma^2}{2}
  \nu\Big(\frac{z}{\ep}\Big)\nu'\Big(\frac{z}{\ep}\Big)\Big]\ep
\partial_z p_q^\ep(\om,z) \nonumber \\& + \frac{D'(z)}{D(z)} \,
\frac{2jq[1+(-1)^{j+q}]}{(q^2-j^2)} \ep \partial_z p_q^\ep(\om,z)
\nonumber\\&+\frac{j q (-1)^{j+q}}{(q^2-j^2)}
\Big[\frac{\sigma}{\sqrt{\ep}} \nu''\Big(\frac{z}{\ep}\Big) -
  \frac{\sigma^2}{2} \nu\Big(\frac{z}{\ep}\Big)
  \nu''\Big(\frac{z}{\ep}\Big) \Big] p_q^\ep(\om,z)
\nonumber\\ &+\frac{jq (j^2+q^2)(-1)^{j+q}}{(q^2-j^2)^2} \sigma^2
            {\nu'}^{\,2}\Big(\frac{z}{\ep}\Big)
            p_q^\ep(\om,z)\nonumber \\ &+\frac{\kappa(z)}{D(z)}
            \frac{2 j q [1-(-1)^{j+q}]\big[ j^2 + 3 q^2 - 4
                \big(\frac{k D(z)}{\pi}\big)^2 \big]}{(q^2-j^2)^2}
            p_q^\ep(\om,z)\Big\}, \label{eq:MC2}
\end{align}
where we obtained from \eqref{eq:Coef1}-\eqref{eq:Coef5} that
\begin{align*}
\lin \big(2 \rho + D\big) y_j,\partial_\rho y_q\rin &=\frac{4 j q
  (-1)^{j+q}}{q^2-j^2}, \\ \lin y_j, \partial_z y_q\rin &=
\frac{D'(z)}{D(z)} \,\frac{jq[1+(-1)^{j+q}]}{j^2-q^2},
\\ \frac{\mu_q^2(z) \lin (2\rho + D)^2 y_j,y_q \rin}{16} - \frac{\lin
  \big(2 \rho + D\big) y_j,\partial_\rho y_q\rin}{4} & = \frac{jq
  (j^2+q^2)(-1)^{j+q}}{(q^2-j^2)^2}, \\ \big(k^2-\mu_q^2(z)\big) \lin
(2 \rho +D) y_j,y_q\rin + \lin y_j,\partial_\rho y_q \rin &= \frac{2 j
  q [1-(-1)^{j+q}]\big[ j^2 + 3 q^2 - 4
  \big(\frac{k D(z)}{\pi}\big)^2 \big]}{D(z)(q^2-j^2)^2} .
\end{align*}

We now use integrating factors to simplify equations \eqref{eq:MC1}.
Specifically, we define 
\begin{equation}
u_j^\ep(\om,z) = p_j^\ep(\om,z) \exp \Big[ \frac{\sigma \sqrt{\ep}}{4}
  \nu \Big(\frac{z}{\ep}\Big) - \frac{\sigma^2 \ep}{16}
  \nu^2\Big(\frac{z}{\ep}\Big) \Big] =
p_j^\ep(\om,z)[1+O(\sqrt{\ep})],
\label{eq:MC3}
\end{equation}
and obtain after substituting in \eqref{eq:MC1} that 
\begin{align}
\frac{1}{\ep} \Big[ (\ep \partial_z)^2 + k^2 -
  \mu_j^2(z)\Big]u_j^\ep(\om,z) + \frac{\sigma}{\sqrt{\ep}} \mu_j^2(z)
\nu\Big(\frac{z}{\ep}\Big)u_j^\ep(\om,z) \nonumber \\+ \sigma^2
\Big\{-\frac{3}{4} \mu_j^2(z) \nu^2\Big(\frac{z}{\ep}\Big) - \Big[
  \frac{(\pi j)^2}{12} + \frac{1}{16}\Big]
            {\nu'}^{\,2}\Big(\frac{z}{\ep}\Big)\Big\} u_j^\ep(\om,z)
            \nonumber \\ \approx \mathcal{C}_j^\ep(\om,z) + f(\om)
            y_j(\rho_\star,0)\delta(z),
\end{align}
with coupling term 
\begin{align}
\mathcal{C}_j^\ep(\om,z) = \sum_{q=1, q \ne j}^\infty & \Big\{ \frac{2 j q
  (-1)^{j+q}}{(q^2-j^2)} \Big[ \frac{\sigma}{\sqrt{\ep}}
  \nu'\Big(\frac{z}{\ep}\Big) - \frac{\sigma^2}{2}
  \nu\Big(\frac{z}{\ep}\Big)\nu'\Big(\frac{z}{\ep}\Big)\Big]\ep
\partial_z u_q^\ep(\om,z) \nonumber \\& + \frac{D'(z)}{D(z)}
\frac{2jq[1+(-1)^{j+q}]}{(q^2-j^2)} \ep \partial_z u_q^\ep(\om,z)
\nonumber\\&+\frac{j q (-1)^{j+q}}{(q^2-j^2)}
\Big[\frac{\sigma}{\sqrt{\ep}} \nu''\Big(\frac{z}{\ep}\Big) -
  \frac{\sigma^2}{2} \nu\Big(\frac{z}{\ep}\Big)
  \nu''\Big(\frac{z}{\ep}\Big) \Big] u_q^\ep(\om,z)
\nonumber\\ &+\frac{jq (3j^2+q^2)(-1)^{j+q}}{2(q^2-j^2)^2} \sigma^2
            {\nu'}^{\,2}\Big(\frac{z}{\ep}\Big)
            u_q^\ep(\om,z)\nonumber \\ &+\frac{\kappa(z)}{D(z)}
            \frac{2 j q [1-(-1)^{j+q}]\big[ j^2 + 3 q^2 - 4
                \big(\frac{k D(z)}{\pi}\big)^2 \big]}{(q^2-j^2)^2}
            u_q^\ep(\om,z)\Big\}. \label{eq:MC4}
\end{align}
This is the expression \eqref{eq:DC5} and the leading coupling
coefficients are 
\begin{align}
\Gamma_{jq} = \frac{jq(-1)^{j+q}}{(q^2-j^2)} , \quad \Theta_{jq} =
\frac{2jq(-1)^{j+q}}{(q^2-j^2)}.\label{eq:DC6}
\end{align}
The second order coefficients, due to the random fluctuations, are 
\begin{align}
\gamma_{jq}\Big(\frac{z}{\ep}\Big) &= \frac{jq(-1)^{j+q}}{2(q^2-j^2)}\Big[ \frac{(3
    j^2+q^2)}{(q^2-j^2)} {\nu'}^{\,2}\Big(\frac{z}{\ep}\Big) -
  \nu\Big(\frac{z}{\ep}\Big)\nu''\Big(\frac{z}{\ep}\Big)\Big],
\label{eq:DC7}\\
\theta_{jq}\Big(\frac{z}{\ep}\Big) &=-\frac{jq(-1)^{j+q}}{(q^2-j^2)}
\nu\Big(\frac{z}{\ep}\Big)\nu'\Big(\frac{z}{\ep}\Big),\label{eq:DC8}
\end{align}
and those due to the slow changes in the waveguide are
\begin{align}
\gamma_{jq}^o(z) &= \frac{\kappa(z)}{D(z)}
            \frac{2 j q [1-(-1)^{j+q}]\big[ j^2 + 3 q^2 - 4
                \big(\frac{k D(z)}{\pi}\big)^2 \big]}{(q^2-j^2)^2},
\label{eq:DC9} \\
\theta_{jq}^o(z) &=\frac{D'(z)}{D(z)} \frac{2 j q[1+(-1)^{j+q}]}{(q^2-j^2)}.
\label{eq:DC10}
\end{align}

\section{Useful identities}
\label{ap:identities}
Here we give a few identities satisfied by the eigenfunctions
\eqref{eq:eigf}, for all $z \in \mathbb{R}$. The first identity
is just the statement that the eigenfunctions are orthonormal 
\begin{equation}
\int_{-D(z)/2}^{D(z)/2} d \rho \, y_j(\rho,z) y_q(\rho,z) =
\delta_{jq},
\label{eq:orthog}
\end{equation}
where $\delta_{jq}$ is the Kronecker delta symbol. The second identity
\begin{align}
\int_{-D(z)/2}^{D(z)/2} d \rho \, \rho y_j^2(\rho,z)&= 0,
\label{eq:odd}
\end{align}
is because the integrand is odd. The third identity
follows from the fundamental theorem of calculus,
\begin{align}
\int_{-D(z)/2}^{D(z)/2} d \rho\,y_j(\rho,z) \partial_\rho
y_j(\rho,z) &= \frac{1}{2}\int_{-D(z)/2}^{D(z)/2} d \rho \,
\partial_\rho y_j^2(\rho,z) = 0,
\end{align}
because the eigenfunctions vanish at $\rho = \pm D(z)/2$. 
The fourth identity is 
\begin{align}
\int_{-D(z)/2}^{D(z)/2} d \rho \, [2\rho + D(z)] y_j(\rho,z)
\partial_\rho y_j(\rho,z) = \int_{-D(z)/2}^{D(z)/2} d \rho \, \rho
\partial_\rho y_j^2(\rho,z) \nonumber \\ =\int_{-D(z)/2}^{D(z)/2} d
\rho \left\{ \partial_\rho \left[\rho y_j^2(\rho,z)\right] -
y_j^2(\rho,z)\right\} = -1,
\end{align}
where we used integration by parts. The fifth identity is 
\begin{equation}
\int_{-D(z)/2}^{D(z)/2} d \rho \, y_j(\rho,z) \partial_z y_j(\rho,z) = 0.
\end{equation}
To derive it, we take the $z$ derivative  in
\eqref{eq:orthog}, for $q=j$, and obtain that
\begin{align*}
0 =& \partial_z \int_{-D(z)/2}^{D(z)/2} d \rho \, y_j^2(\rho,z) = 2
\int_{-D(z)/2}^{D(z)/2} d \rho \, y_j(\rho,z) \partial_z
y_j(\rho,z) \\&+ \frac{D'(z)}{2} \left[ y_j^2(D(z)/2,z)-
  y_j^2(-D(z)/2,z)\right] = 2
\int_{-D(z)/2}^{D(z)/2} d \rho \, y_j(\rho,z) \partial_z
y_j(\rho,z).
\end{align*}
We also have  from \eqref{eq:orthog}, \eqref{eq:odd}, and
definition \eqref{eq:eigf} that 
\begin{align}
\int_{-D(z)/2}^{D(z)/2} d \rho [2\rho + D(z)]^2 y_j^2(\rho,z) &=
D^2(z) + \frac{8}{D(z)} \int_{-D(z)/2}^{D(z)/2} d \rho\, \rho^2 \sin^2
\left[ \left(\frac{\rho}{D(z)} + \frac{1}{2} \right) \pi j\right]
\nonumber \\&= D^2(z) \left[\frac{4}{3}-\frac{2}{(\pi j)^2}\right]
\label{eq:D6}.
\end{align}

For $j \ne q$ we  have from definition \eqref{eq:eigf} of the
eigenfunctions that 
\begin{align}
\int_{-D(z)/2}^{D(z)/2} d \rho\, [2\rho + D(z)] y_j(\rho,z)
\partial_\rho y_q(\rho,z) = 2 \pi q \int_{-D(z)/2}^{D(z)/2} \frac{d
  \rho}{D(z)} \Big[\frac{2\rho}{D(z)} + 1\Big] \nonumber \\
\times \sin\left[
  \left(\frac{\rho}{D(z)} + \frac{1}{2} \right) \pi j\right]
\cos\left[ \left(\frac{\rho}{D(z)} + \frac{1}{2} \right) \pi q\right] = -
\frac{4 j q(-1)^{j+q}}{j^2-q^2}.
\label{eq:Coef1}
\end{align}
Similarly, we obtain after taking the derivative with respect to $z$ of 
$y_q(\rho,z)$ and substituting in the integral below that 
\begin{align}
\int_{-D(z)/2}^{D(z)/2} d \rho \, y_j(\rho,z)
\partial_z y_q(\rho,z) = \frac{D'(z)}{D(z)} \, \frac{jq
  \Big[(-1)^{j+q}+1\Big]}{j^2-q^2}.
\label{eq:Coef2}
\end{align}
We also calculate using the expression \eqref{eq:eigf} that
\begin{align}
\int_{-D(z)/2}^{D(z)/2} d \rho\, [2\rho + D(z)]^2 y_j(\rho,z)
y_q(\rho,z) = \frac{32 D^2(z)}{\pi^2} \,\frac{j q (-1)^{j+q}}{(j^2-q^2)},
\label{eq:Coef3}
\end{align}
and 
\begin{align}
\int_{-D(z)/2}^{D(z)/2} d \rho \, y_j(\rho,z) \partial_\rho y_q(\rho,z) = 
\frac{2 j q [1-(-1)^{j+q}]}{D(z)(j^2-q^2)},
\label{eq:Coef4}
\end{align}
and 
\begin{align}
\int_{-D(z)/2}^{D(z)/2} d \rho \, (2\rho + D(z)) y_j(\rho,z)y_q(\rho,z) = 
-\frac{8 D(z) j q [1-(-1)^{j+q}]}{\pi^2 (j^2-q^2)^2}.
\label{eq:Coef5}
\end{align}
% --------------------
\section{The evanescent waves}
\label{ap:evanesc}
Let us begin by rewriting equation \eqref{eq:DC3} in first order system form, 
for the unknown vector with components $u_j^\ep(\om,z)$ and 
\begin{equation}
v_j^\ep(\om,z) = \frac{\ep}{\beta_j(\om,z)} \partial_z u_j^\ep(\om,z),
\label{eq:E1}
\end{equation}
where $j > \mathcal{N}$ and $z \in \big( z_{-}^{(t)},z_{-}^{(t-1)}
\big)$. The mode wave number $\beta_j$ is defined in \eqref{eq:BetaEv},
and the system is
\begin{align}
 &\left\{ \partial_z  - \frac{\beta_j(\om,z)}{\ep}
\begin{pmatrix} 0 & 1 \\ 1 & 0 \end{pmatrix}  + \Big[
   \frac{\sigma \mu_j^2(z)}{\sqrt{\ep}\beta_j(\om,z)} \nu
   \big(\frac{z}{\ep}\big) +
   \frac{\sigma^2g_j^\ep(\om,z)}{\beta_j(\om,z)} \Big] \begin{pmatrix}
  0 & 0 \\ 1 & 0 \end{pmatrix} \right\} \begin{pmatrix} u_j^\ep(\om,z)
  \\ v_j^\ep(\om,z) \end{pmatrix} \nonumber \\&\hspace{2.5in}=
\frac{\mathcal{C}_j^\ep(\om,z)}{\beta_j(\om,z)} \begin{pmatrix} 0
  \\ 1 \end{pmatrix}. \label{eq:E2}
\end{align}
The matrix $\begin{pmatrix} 0 & 1 \\ 1 & 0 \end{pmatrix}$ in the
leading term has the eigenvalues $\pm 1$, and the orthonormal
eigenfunctions $\frac{1}{\sqrt{2}} \begin{pmatrix} 1 \\ \pm
  1 \end{pmatrix}$. Expanding the solution in the basis of these 
eigenfunctions 
\begin{equation}
\begin{pmatrix} u_j^\ep(\om,z)
  \\ v_j^\ep(\om,z) \end{pmatrix} =
\frac{\alpha_j^+(\om,z)}{\sqrt{2}} \begin{pmatrix} 1 \\ 1 \end{pmatrix}
+ \frac{\alpha_j^-(\om,z)}{\sqrt{2}} \begin{pmatrix} 1 \\ -
  1 \end{pmatrix},
\label{eq:E3}
\end{equation}
and substituting in \eqref{eq:E2} gives the following equations for
the coefficients
\begin{align}
\left[ \partial_z \mp \frac{\beta_j(\om,z)}{\ep}\right]
\alpha_j^\pm(\om,z) &= \pm\frac{\mathcal{C}_j^\ep(\om,z)}{\sqrt{2}
  \beta_j(\om,z)} \mp \frac{\big[ \alpha_j^+(\om,z) +
    \alpha_j^-(\om,z)\big]}{2 \beta_j(\om,z)} \nonumber \\ & \times \Big[
  \frac{\sigma}{\sqrt{\ep}}\mu_j^2(z) \nu \Big(\frac{z}{\ep}\Big) +
  \sigma^2 g_j^\ep(\om,z)\Big].
\label{eq:E4}
\end{align}
These are complemented with the boundary conditions
\begin{equation}
\alpha_j^+\big(\om, z_{-}^{(t-1)}\big) = \sqrt{2}c_j^{(t)+}, \quad
\alpha_j^-\big(\om, z_{-}^{(t)} \big) = 0,
\label{eq:E5}
\end{equation}
with constant $c_j^{(t)}$ to be determined later, indexed by $t$ to
remind us that we work in the sector $z \in
\big(z_{-}^{(t)},z_{-}^{(t-1)}\big)$. In \eqref{eq:E5} we set to zero
the component $\alpha_j^{-}$ at the farther end $z_{-}^{(t)}$ from the
source, to suppress the growing part of the solution.

We obtain after integration of \eqref{eq:E4} that
\begin{align}
 \alpha_j^+(\om,z) &= \sqrt{2}c_j^{(t)} \exp \Big[- \frac{1}{\ep}
   \int_{z}^{z_{-}^{(t-1)}} d \zeta \,\beta_j(\om,\zeta)\Big] -
 \int_z^{z_{-}^{(t-1)}} \hspace{-0.1in}d \zeta \frac{\exp \Big[
    - \frac{1}{\ep} \int_z^{\zeta} ds \, \beta_j(\om,s) \Big]}{\sqrt{2}
   \beta_j(\om,\zeta)} \nonumber \\ &\hspace{-0.1in} \times
 \left\{\mathcal{C}_j^\ep(\om,\zeta) - \frac{\big[ \alpha_j^+(\om,\zeta) +
     \alpha_j^-(\om,\zeta)\big]}{\sqrt{2}} \Big[
   \frac{\sigma}{\sqrt{\ep}}\mu_j^2(\zeta) \nu
   \Big(\frac{\zeta}{\ep}\Big) + \sigma^2
   g_j^\ep(\om,\zeta)\Big]\right\}, \label{eq:E6}
\end{align}
and 
\begin{align}
\alpha_j^-(\om,z) = - \int_{z_{-}^{(t)}}^z &d \zeta \frac{\exp \Big[
   - \frac{1}{\ep} \int_{\zeta}^z ds \, \beta_j(\om,s) \Big]}{\sqrt{2}
  \beta_j(\om,\zeta)} \left\{\mathcal{C}_j^\ep(\om,\zeta) - \frac{\big[
    \alpha_j^+(\om,\zeta) + \alpha_j^-(\om,\zeta)\big]}{\sqrt{2}}
\right. \nonumber \\ &\left. \times \Big[
  \frac{\sigma}{\sqrt{\ep}}\mu_j^2(\zeta) \nu
  \Big(\frac{\zeta}{\ep}\Big) + \sigma^2
  g_j^\ep(\om,\zeta)\Big]\right\}. \label{eq:E7}
\end{align}
All the exponential terms in these equations are decaying in $z$, so we can
change the variable of integration as $ \zeta = z+ \ep \xi, $ and note
that only $\xi = O(1)$ contributes to the result.  Equation
\eqref{eq:E6} becomes
\begin{align}
&\alpha_j^+(\om,z) \approx \sqrt{2}c_j^{(t)} \exp \Big[ - \frac{1}{\ep}
    \int_z^{z_{-}^{(t-1)}} d \zeta \, \beta_j(\om,\zeta)\Big] -
  \frac{\ep}{\sqrt{2} \beta_j(\om,z)} \int_0^\infty d \xi \, e^{- \xi
    \beta_j(\om,z)}\nonumber \\& \hspace{0.2in} \times
  \left\{\mathcal{C}_j^\ep(\om,z + \ep \xi) - u_j^\ep(\om,z+\ep \xi)
  \Big[ \frac{\sigma}{\sqrt{\ep}}\mu_j^2(z) \nu \Big(\frac{z}{\ep}+
    \xi\Big) + \sigma^2 g_j^\ep(\om,z + \ep \xi
    )\Big]\right\}, \label{eq:E8}
\end{align}
where we used \eqref{eq:E3} in the integrand, and the approximation
means that we neglect terms that vanish in the limit $\ep \to
0$. Similarly, equation \eqref{eq:E7} becomes
\begin{align}
\alpha_j^-(\om,z) \approx - \frac{\ep}{\sqrt{2} \beta_j(\om,z)}
\int_{-\infty}^0 d \xi \,&e^{\xi \beta_j(\om,z)}
\left\{\mathcal{C}_j^\ep(\om,z + \ep \xi) - u_j^\ep(\om,z+\ep \xi)
\right. \nonumber \\ & \times \left. \Big[
  \frac{\sigma}{\sqrt{\ep}}\mu_j^2(z) \nu \Big(\frac{z}{\ep}+ \xi\Big)
  + \sigma^2 g_j^\ep(\om,z + \ep \xi )\Big]\right\}. \label{eq:E9}
\end{align}
The expression of $u_j^\ep$ follows from these equations and
\eqref{eq:E3},
\begin{align}
&u_j^\ep(\om,z) \approx c_j^{(t)}(\om) \exp \Big[ - \frac{1}{\ep}
    \int_z^{z_{-}^{(t-1)}} d \zeta \, \beta_j(\om,\zeta)\Big] -
  \frac{\ep}{2\beta_j(\om,z)} \int_{-\infty}^\infty d \xi \, e^{- |\xi|
    \beta_j(\om,z)}\nonumber \\& \hspace{0.1in} \times
  \left\{\mathcal{C}_j^\ep(\om,z + \ep \xi) - u_j^\ep(\om,z+\ep \xi)
  \Big[ \frac{\sigma}{\sqrt{\ep}}\mu_j^2(z) \nu \Big(\frac{z}{\ep}+
    \xi\Big) + \sigma^2 g_j^\ep(\om,z + \ep \xi
    )\Big]\right\}. \label{eq:E10}
\end{align}
The derivative $\ep \partial_z u_j^\ep$ is obtained from
\eqref{eq:E1}, \eqref{eq:E3}, \eqref{eq:E8}-\eqref{eq:E9} and
integration by parts
\begin{align}
&\ep \partial_z u_j^\ep(\om,z) \approx \beta_j(\om,z) c_j^{(t)}(\om)
  \exp \Big[\frac{1}{\ep} \int_{z_{-}^{(t-1)}}^z d \zeta \,
    \beta_j(\om,\zeta)\Big] - \frac{\ep}{2\beta_j(\om,z)}
  \int_{-\infty}^\infty d \xi \, e^{- |\xi| \beta_j(\om,z)}\nonumber
  \\& \times \ep \partial_z \left\{
  \mathcal{C}_j^\ep(\om,z + \ep \xi) -
  u_j^\ep(\om,z+\ep \xi) \Big[ \frac{\sigma}{\sqrt{\ep}}\mu_j^2(z) \nu
    \Big(\frac{z}{\ep}+ \xi\Big) + \sigma^2 g_j^\ep(\om,z + \ep \xi
    )\Big]\right\}. \label{eq:E10p}
\end{align}

Now let us recall the expression \eqref{eq:DC5} of
$\mathcal{C}_j^\ep(\om,z)$, which models the coupling with the other
modes, and write it as the sum of two terms:
\begin{equation}
\mathcal{C}_j^\ep(\om,z) = \mathcal{C}_j^{\ep (p)}(\om,z) +
\mathcal{C}_j^{\ep(e)}(\om,z).
\label{eq:E11}
\end{equation}
The first term is the coupling with the propagating modes, and is
given by restricting the sum in \eqref{eq:DC5} to $ q \le
\mathcal{N}$. The second term is the remaining series, with terms
indexed by $q > \mathcal{N}$, and $q \ne j$. Each term in this series
involves $u_q^\ep(\om,z)$ and $\ep \partial_z u_q^\ep(\om,z)$ that
have expressions like \eqref{eq:E10}-\eqref{eq:E10p}.  Stringing all
the unknowns in the infinite-dimensional vector $\bm{U} = \big(
\bm{U}_{\mathcal{N}+1}, \bm{U}_{\mathcal{N}+2}, \ldots, \big)$ where
$\bm{U}_j = (u_j^\ep,\ep \partial_z u_j^\ep)$, for $j > \mathcal{N}$,
we can write equations \eqref{eq:E10}-\eqref{eq:E10p} in compact form
as
\begin{equation}
\Big({\bf I}-\ep {\bf K}\Big) \bm{U}(\om,z) = \bm{F}(\om,z),
\label{eq:E12} 
\end{equation}
with right hand side given by the concatenation of 
\begin{align}
\bm{F}_j(\om,z) = &\begin{pmatrix} 1 \\ \beta_j(\om,z) \end{pmatrix}
c_j^{(t)}(\om) \exp \Big[- \frac{1}{\ep} \int_z^{z_{-}^{(t-1)}} d \zeta
  \, \beta_j(\om,\zeta)\Big] \nonumber \\ &-
\frac{\ep}{2\beta_j(\om,z)} \int_{-\infty}^\infty d \xi e^{- |\xi|
  \beta_j(\om,z)}  \begin{pmatrix} 1 \\ \ep
  \partial_z \end{pmatrix} \mathcal{C}_j^{\ep(p)}(\om,z + \ep \xi),
\label{eq:E13}
\end{align}
for $j \ge \mathcal{N}$. In the left hand side of \eqref{eq:E12} we
have the perturbation of the identity ${\bf I}$ by the integral
operator ${\bf K}$, whose kernel follows easily from the $(u_q^\ep)_{q
  > \mathcal{N}}$ dependent terms in the integrand in
\eqref{eq:E10}-\eqref{eq:E10p}, including those in $\mathcal{C}_j^{\ep
  (e)}$. This integral operator is basically the same as that analyzed
in \cite[Lemma 3.1]{alonso2011wave}, and it is bounded with respect to
an appropriate norm. This means that we can solve \eqref{eq:E12} using
Neumann series and obtain
\begin{equation}
\bm{U}(\om,z) = \bm{F}(\om,z) + O(\ep).
\label{eq:E14}
\end{equation}

The first term in \eqref{eq:E13} matters only in the $O(\ep)$ vicinity
of $z_{-}^{(t-1)}$, over which the mode coupling is negligible. The
constant $c_j^{(t)}$ is determined by continuity conditions at
$z_{-}^{(t-1)}$ as follows: If $t = 1$, so that $z_{-}^{(t-1)} = 0$,
$c_j^{(1)}$ is determined by the source excitation, and it equals the
coefficient of the $j$-th evanescent mode in the perfect waveguide
with width $D(0)$. If $t > 1$, then $c_{j}^{(t)}$ is determined by
continuity of the wave field at the turning point $z_{-}^{(t-1)}$.

Assuming that $z_{-}^{(t-1)} - z \gg \ep$, so we can neglect the
first term in \eqref{eq:E13}, we have
\begin{align}
\begin{pmatrix} u_j^\ep(\om,z) \\
\ep \partial_z u_j^\ep(\om,z) \end{pmatrix} \approx -
\frac{\ep}{2\beta_j(\om,z)} \int_{-\infty}^\infty d \xi e^{- |\xi|
  \beta_j(\om,z)} \begin{pmatrix} 1 \\ \ep \partial_z \end{pmatrix}
\mathcal{C}_j^{\ep(p)}(\om,z + \ep \xi),
\label{eq:CE}
\end{align}
with
\begin{align}
\hspace{-0.1in}\ep \mathcal{C}_j^{\ep(p)}(\om,z + \ep \xi) \approx &
    {\sigma}\sqrt{\ep} \sum_{q=1}^\mathcal{N} \Big[\Gamma_{jq} \nu''
      \Big(\frac{z}{\ep} + \xi\Big) + \Theta_{jq} \nu'
      \Big(\frac{z}{\ep} + \xi\Big) \ep \partial_z\Big]
    u_q^\ep(\om,z+\ep \xi), \label{eq:CE1}
\end{align}
obtained from \eqref{eq:DC5}.  Here the modes $u_q^\ep$ and their
derivative $\ep \partial_z u_q^\ep(\om,z)$ are given in
\eqref{eq:FM10}-\eqref{eq:FM11}, and the constant coefficients
$\Gamma_{jq}$ and $\Theta_{jq}$ are defined in \eqref{eq:DC6}.
Substituting in \eqref{eq:CE1} and then \eqref{eq:CE}, and using that
the derivatives of the mode amplitudes given in \eqref{eq:FM12} are at
most $O(\ep^{-1/2})$, we obtain equation \eqref{eq:CE2}. The
derivative in the integrand in \eqref{eq:CE} is
\begin{align}
\ep \partial_z \Big[ \ep \mathcal{C}_j^{\ep (p)}(\om,z
    + \ep \xi) \Big] &= \sigma \sqrt{\ep} \sum_{q=1}^\mathcal{N} \Big[
    \Gamma_{jq} \nu'''\Big(\frac{z}{\ep} + \xi\Big) + \big(\Gamma_{jq}
    +\Theta_{jq}\big) \nu''\Big(\frac{z}{\ep} + \xi\Big) \ep
    \partial_z \nonumber \\&- \beta_q^2(\om,z) \Theta_{jq} \nu'
    \Big(\frac{z}{\ep} + \xi \Big) \Big]u_q^\ep(\om,z + \ep \xi)
  , \label{eq:CE3}
\end{align}
where we used equation \eqref{eq:DC3} for $(\ep \partial_z)^2
u_q^\ep$. Substituting \eqref{eq:FM10}-\eqref{eq:FM11} in \eqref{eq:CE3}
and then in \eqref{eq:CE}, we obtain \eqref{eq:CE4}.

\bibliographystyle{siam} \bibliography{TURNING}

\end{document}